\documentclass[preprint,3p,authoryear]{elsarticle}

\usepackage[pdftex,backref=section,hypertexnames=true,plainpages=false,
naturalnames,colorlinks=true,linkcolor=blue,bookmarks]{hyperref}
\usepackage{amsfonts,amssymb,amsmath,amsopn,float}
\usepackage[nameinlink,capitalise]{cleveref}
\usepackage{xifthen}
\usepackage{rotating,graphicx}
\graphicspath{{./files/}}
\usepackage{multicol}
\usepackage{tikz,pgfplots}
\pgfplotsset{compat=1.13}
\usepackage{epstopdf}
\usepackage{url}
\usepackage{subcaption}

\usepackage{algorithm}
\usepackage{algcompatible}

\usepackage{siunitx}
\usepackage{xfrac}

\usepackage{color}
\definecolor{mygray}{RGB}{47,79,79}

\usepackage{comment}

\newtheorem{theorem}{Theorem}
\newtheorem{lem}[theorem]{Lemma}

\newtheorem{rem}[theorem]{Remark}
\crefname{rem}{Remark}{Remarks}
\newtheorem{definition}[theorem]{Definition}
\crefname{definition}{Definition}{Definitions}
\newcommand{\sref}[2]{\hyperref[#2]{#1 \ref*{#2}}}
\newcommand{\sautoref}[2]{\hyperref[#2]{#1 \ref*{#2}}}

\newcommand{\bcolon}{\boldsymbol{:}}
\newcommand{\bbR}{\mathbb{R}}

\newcommand{\bzero}{\mathbf{0}}

\newcommand{\bxi}{{\boldsymbol{\xi}}}
\newcommand{\btau}{{\boldsymbol{\tau}}}

\newcommand{\Ltwo}[1]{%
	\ifthenelse{\equal{#1}{}}{L^2}{L^2(#1)}%
}

\newcommand{\Ltwoz}[1]{%
	\ifthenelse{\equal{#1}{}}{L^2_0}{L^2_0(#1)}%
}

\newcommand{\Cone}[1]{%
	\ifthenelse{\equal{#1}{}}{C^{1}}{C^{1}(#1)}%
}

\newcommand{\Conez}[1]{%
	\ifthenelse{\equal{#1}{}}{C^{1}_{0}}{C^{1}_{0}(#1)}%
}

\newcommand{\Ctwo}[1]{%
	\ifthenelse{\equal{#1}{}}{C^{2}}{C^2(#1)}%
}

\newcommand{\Ctwoz}[1]{%
	\ifthenelse{\equal{#1}{}}{C^{2}_{0}}{C^{2}_{0}(#1)}%
}

\newcommand{\Cholder}[1]{%
	\ifthenelse{\equal{#1}{}}{C^{0,\gamma}}{C^{0,\gamma}(#1)}%
}

\newcommand{\Cholderz}[1]{%
	\ifthenelse{\equal{#1}{}}{C^{0,\gamma}_{0}}{C^{0,\gamma}_{0}(#1)}%
}

\newcommand{\bolds}[1]{\boldsymbol{#1}}

\newcommand{\bb}{\bolds{b}}

\newcommand{\be}{\bolds{e}}
\newcommand{\bff}{\bolds{f}}

\newcommand{\br}{\bolds{r}}

\newcommand{\bu}{\bolds{u}}
\newcommand{\bv}{\bolds{v}}
\newcommand{\bw}{\bolds{w}}
\newcommand{\bx}{\bolds{x}}
\newcommand{\by}{\bolds{y}}
\newcommand{\bz}{\bolds{z}}

\newcommand{\bE}{\bolds{E}}
\newcommand{\bF}{\bolds{F}}

\newcommand{\bU}{\bolds{U}}

\newcommand{\bX}{\bolds{X}}

\numberwithin{theorem}{section}

\ifpdf
\DeclareGraphicsExtensions{.eps,.pdf,.png,.jpg}
\else
\DeclareGraphicsExtensions{.eps}
\fi

\begin{document}
	
	\begin{frontmatter}
		
		\title{Nodal finite element approximation of peridynamics}
		
		\author[aa]{Prashant K. Jha\corref{cor1}}
		\ead{pjha.sci@gmail.com}
		
		\author[cc,phys]{Patrick Diehl}
		\ead{pdiehl@cct.lsu.edu}
		
		\author[bb,cc]{Robert Lipton}
		\ead{lipton@lsu.edu}
		
		\address[aa]{Department of Mechanical Engineering, South Dakota School of Mines and Technology, Rapid City, USA}
		
		\address[bb]{Department of Mathematics, Louisiana State University, Baton Rouge, USA}
		
		\address[cc]{Center for Computation and Technology, Louisiana State University, Baton Rouge, USA}
		
		\address[phys]{Department of Physics and Astronomy, Louisiana State University, Baton Rouge, USA
			\\[4pt]
			{\bf\large
				\begin{center}
					$^\dagger$Published in Computer Methods in Applied Mechanics and Engineering DOI: \url{https://doi.org/10.1016/j.cma.2024.117519}
				\end{center}
			}
		}

		\cortext[cor1]{Corresponding author (pjha.sci@gmail.com)}
		
		
		\begin{abstract}
			This work considers the nodal finite element approximation of peridynamics, in which the nodal displacements satisfy the peridynamics equation at each mesh node. For the nonlinear bond-based peridynamics model, it is shown that, under the suitable assumptions on an exact solution, the discretized solution associated with the central-in-time and nodal finite element discretization converges to a solution in the $L^2$ norm at the rate $C_1 \Delta t + C_2 h^2/\epsilon^2$. Here, $\Delta t$, $h$, and $\epsilon$ are time step size, mesh size, and the size of the horizon or nonlocal length scale, respectively. Constants $C_1$ and $C_2$ are independent of $h$ and $\Delta t$ and depend on norms of the solution and nonlocal length scale. Several numerical examples involving pre-crack, void, and notch are considered, and the efficacy of the proposed nodal finite element discretization is analyzed. 
		\end{abstract}
		
		
	\end{frontmatter}
	
	\section{Introduction}
	
	Peridynamics is a reformulation of classical continuum mechanics introduced by Silling in~\cite{CMPer-Silling,States}. The strain inside the medium is expressed in terms of displacement differences as opposed to the displacement gradients, and the internal force at a material point is due to the sum of all pairwise interactions between a point and its neighboring points. The new formulation bypasses the difficulty incurred by displacement gradients and discontinuities, as in the case of classical fracture theories. The nonlocal fracture theory has been applied numerically to model the complex fracture phenomenon in materials, see, \emph{e.g.}, \cite{WeckAbe,SillBob,CMPer-Silling4,CMPer-Silling5,CMPer-Silling7,HaBobaru,CMPer-Agwai,BobaruHu,CMPer-Ghajari,CMPer-Lipton2,DuTaoTian,lipton2019complex,Jha2020peri,jha2020peridynamics}. \cite{diehl2019review} is referred for a comprehensive survey. In peridynamics, every point interacts with its neighbors inside a ball of fixed radius called the horizon. The size of the horizon sets the length scale of nonlocal interaction. When the forces between points are linear and when the nonlocal length scale tends to zero, it is seen that peridynamics converge to the classical linear elasticity, \cite{CMPer-Emmrich,CMPer-Silling4,AksoyluUnlu,CMPer-Mengesha2}. For nonlinear forces, in which the bond behaves like an elastic spring for small strains and softens with increasing strains, peridynamics converges in the small horizon limit to linear elastic fracture mechanics, where the material has a sharp crack, and away from a sharp crack the material is governed by linear elastodynamics, see \cite{CMPer-Lipton,CMPer-Lipton3,Jha2020peri,jha2018numerical2}. 
	For simulation of fracture using peridynamics, there are several choices, e.g., meshfree discretization \cite{silling2005meshfree, trask2019asymptotically} (see \cite{seleson2016consistency} for numerical convergence test and overview of existing meshfree methods for peridynamics), commonly used nodal-based discretization (similar to finite difference approximation of partial differential equations) \cite{CMPer-Silling, silling2003dynamic, silling2010peridynamic, jha2018numerical, jha2019numerical}, and finite element approximations and their variants have been used in works such as~\cite{CMPer-Richard,CMPer-Madenci,CMPer-Wildman2,chen2011continuous,diyaroglu2017peridynamic,de2017finite,anicode2022bond,ni2018peridynamic,huang2019finite,yang2019implementation}. To reduce the computational cost associated with the nonlocal interaction, coupled local (classical continuum mechanics) and peridynamics equations are also considered in which the local model is used away from cracking zone while the rest of the domain is modeled using peridynamics, see \cite{liu2012coupling, shojaei2017coupling}.
	
	Existing meshfree methods and methods similar to finite difference discretization where the nonlocal integral is approximated using the node-node interaction offer multiple advantages over standard finite element approximation of peridynamics, such as easy implementation and reduced computational cost. However, when compared to the finite element methods, meshfree methods suffer from poor numerical convergence, a lack of continuous representation of the displacement field, which could be crucial in post-processing, and difficulty in coupling peridynamics with other physics (e.g., heat equation and diffusion equation for corrosion).
	
	Motivated by the above arguments, this work considers nodal finite element approximation (or, in brief, NFEA) that overcomes some of the limitations of the meshfree method while retaining the key features of finite element approximation, such as continuous representation of displacement, convergence estimates, and suitability for combining peridynamics with multiphysics models. In the nodal finite element approximation, the equation for the discretized displacement field is written at each mesh node. In contrast, in the standard finite element approximation (FEA) of peridynamics, the approximate solution satisfies the variational form of the peridynamics equation. 
	Node-based calculations considered in this work are quite suitable for peridynamics/nonlocal equations, where a point nonlocally interacts with neighboring points at a distance larger than the mesh size. Classical finite element discretization of peridynamics, \emph{e.g.}, \cite{jha2020finite, jha2021finite}, involves computing interactions of a quadrature point with all neighboring quadrature points within a nonlocal neighborhood (typically a ball of radius greater than the mesh size). Thus, the computation cost is large and prohibitive if one chooses higher-order quadrature approximations. 
	In contrast, the nodal finite element discretization considered in this work applies a discretized equation at each node, and nonlocal interactions are computed between the mesh nodes.
	Comparing the discretized equations in NFEA and FEA, NFEA includes the appearance of an additional error in representing the peridynamics force; see~\cref{ss:compareFE}.  
	
	The main goal of this work is to perform an error analysis of the NFEA approximation and show a-priori convergence of numerical solutions. The convergence of the numerical approximation is established by combining our previous work on a-priori convergence of finite element approximation of peridynamics \cite{jha2021finite, jha2020finite} with new estimates that control the additional error introduced by nodal finite element approximation. For suitable initial conditions and boundary conditions, the NFEA solutions are shown to converge at a rate $C_1 \Delta t + C_2 h^2$, where $\Delta t$ gives the size of the time step and $h$ mesh size. Here, $C_1$ and $C_2$ are constants independent of $\Delta t$ and $h$ and depend on the nonlocal length scale $\epsilon$, the norm of the exact solution, choice of influence function, and the peridynamics force potential (anti-derivative) $\psi$. Several fracture problems involving pre-crack, void, and notch are presented. These problems not only highlight the efficacy of the NFEA but also show the utility of peridynamics in nucleation and crack propagation.
	
	The solution of the nonlocal problem is more regular than the local solution. For an evolution over the time interval $[0,t_F]$ a solution of the nonlocal problem is in $C^2[0,t_F];H^2(D)\cap H_0^1(D))$, see \cite{jha2021finite, jha2020finite}. On the other hand the $\epsilon=0$ limit solution with the limit taken with respect to the $C([0,t_F];L^2(D))$ topology  was shown to be in the space $C([0,t_F];SBD(\Omega))$, see \cite{CMPer-Lipton}. For a definition of $SBD(\Omega)$, see \cite{Ambrosio}. The space of Special Functions of Bounded Deformation $SBD(\Omega)$ as well as $SBV(\Omega)$ is used to describe displacement fields in fractured media, see \cite{AmbrosioBrades}.  In this work, the location of the large strain represents the crack. The simulations show that the regions of high strain are localized to thin regions of width on the order of the peridynamic horizon. 
	The theoretical error bound deteriorates with the horizon, reflected in the growth of constants appearing in the a-priori error as the horizon becomes smaller.
	Motivated by these considerations, an alternate nodal interpolation given by the Cl\'ement interpolant is introduced.  An improved a-priori convergence of nodal finite element approximation is proved.
	The implementation of Cl\'ement interpolations in NFEA and a-posteriori error estimates will be discussed in future work. It is seen that the nodal finite element approximation using Cl\'ement interpolation is asymptotically compatible under the assumption the $L^\infty$ norm of the solution is bounded uniformly with respect to the horizon. This assumption is strongly supported in the simulations, \cref{s:results}.

	\paragraph{Outline of the article.} In~\cref{s:model}, bond-based peridynamics theory is described, and the peridynamics equation of motion is presented. \cref{s:fespace} develops nodal finite element approximation, and it is compared with the standard finite element approximation. In~\cref{s:convergence}, a-priori convergence of an alternate nodal finite element approximation for the nonlinear bond-based model is stated and proved. In~\cref{sec:Clement Node}, an a-priori convergence of the Cl\'ement nodal finite element approximation is proved. In~\cref{sec:Asymp}, the method is shown to be asymptotically compatible.
	Numerical experiments involving pre-crack, void, and notch are presented in~\cref{s:results}. Conclusions are drawn in~\cref{s:concl}.
	
	\section{Bond-based peridynamics}\label{s:model}
	
	Let $D\subset \bbR^d$, for $d=2,3$, be the material domain and $\epsilon> 0$ denote the size of the horizon. In the peridynamics formulation, a material point $\bx \in D$ interacts with all the material points within a neighborhood of $\bx$. Neighborhood of a point $\bx$ is taken to be the ball of radius $\epsilon$ centered at $\bx$ and is denoted by $H_\epsilon(\bx) = \{ \by\in \bbR^d: |\by - \bx| < \epsilon\}$. In what follows, $\bx\in D$ denote the material point, $\bu (\bx, t)$ the displacement of $\bx$ at time $t$ for $t\in [0, t_F]$, and $\bz(\bx,t) = \bx + \bu(\bx,t)$ current coordinate of $\bx$. The bond strain (or bond stretch or pairwise strain) between material points $\bx$ and $\by$ is defined as
	\begin{equation}\label{eq:bondstrain}
		\tilde{S}(\by,\bx, t) = \frac{|\bz(\by,t) - \bz(\bx,t)| - |\by - \bx|}{|\by - \bx|}.
	\end{equation}
	For prototype microelastic brittle (PMB) material, the pairwise force between $\bx$ and $\by$ takes the form (see~\cite{CMPer-Silling,BobaruHu})
	\begin{equation}\label{eq:pmbbondforce}
		\tilde{\bff}_{pmb}(\by, \bx, t) = c J^\epsilon(|\by - \bx|) \tilde{S}(\by,\bx,t) \mu(\tilde{S}(\by,\bx, t),t) \frac{\bz(\by,t) - \bz(\bx,t)}{|\bz(\by,t) - \bz(\bx,t)|}.
	\end{equation}
	Here, $c$ is a constant that depends on the elastic strength of a material, $J^\epsilon = J^\epsilon(|\by - \bx|)$ the influence function, and $\mu(S, t)$ the bond-breaking function that models the breakage of the bond if the pairwise strain exceeds certain threshold strain:
	\begin{equation}\label{eq:pmbbondbreaking}
		\begin{split}
			\mu(\tilde{S}(\by,\bx,t),t) = \begin{cases}
				1, \qquad \text{if } \tilde{S}(\by, \bx, \tau) < S_c(\by,\bx) \quad \forall \tau \leq t, \\
				0, \qquad \text{otherwise}.
			\end{cases}
		\end{split}
	\end{equation}
	In the above, $S_c$ is the critical bond strain between the material points $\by, \bx$. In the PMB model, $S_c$ is independent of $\by,\bx$. In general, the value of the critical bond strain $S_c$ depends on the critical fracture energy $G_c$ and the elastic strength of the material. The total force at $\bx$ is given by the sum of the pairwise forces in the neighborhood of $\bx$, \emph{i.e.}, 
	\begin{equation}\label{eq:pmbtotalforce}
		\tilde{\bF}_{pmb}(\bx, t) = \int_{H_\epsilon(\bx)\cap D} \tilde{\bff}_{pmb}(\by, \bx, t) d\by.
	\end{equation}
	
	Under the small deformation assumption given by $|\bu(\bx)-\bu(\by)|<<1$, the bond strain $\tilde{S}(\by, \bx,t)$ can be approximated by linearizing $\tilde{S}$ as follows
	\begin{equation}\label{eq:bondstrainsmall}
		S(\by, \bx,t) = \frac{\bu(\by, t) - \bu(\bx, t)}{|\by - \bx|} \cdot \frac{\by - \bx}{|\by - \bx|} \approx \tilde{S}(\by, \bx,t)
	\end{equation}
	and the pairwise force taking the form
	\begin{equation}\label{eq:pmbbondforcesmall}
		\bff_{pmb}(\by, \bx, t) = c J^\epsilon(|\by - \bx|) S(\by,\bx,t) \mu(S(\by,\bx, t),t) \frac{\by - \bx}{|\by - \bx|},
	\end{equation}
	where the force now acts along the bond vector in the reference configuration, \emph{i.e.}, $\frac{\by - \bx}{|\by - \bx|}$. The total force at a material point is simply
	\begin{equation}\label{eq:pmbtotalforcesmall}
		\bF_{pmb}(\bx, t) = \int_{H_\epsilon(\bx)\cap D} \bff_{pmb}(\by, \bx, t) d\by.
	\end{equation}
	In this treatment, $\tilde{\,\cdot\,}$, \emph{e.g.}, pairwise strain $\tilde{S}$ and force $\tilde{\bF}_{pmb}$, indicates the quantity associated with the large deformation, whereas the notations without $\tilde{ }$, \emph{e.g.}, $S$ and $\bF_{pmb}$, correspond to the small deformation assumption.
	
\begin{figure}
		\centering
		\includegraphics[scale=0.14]{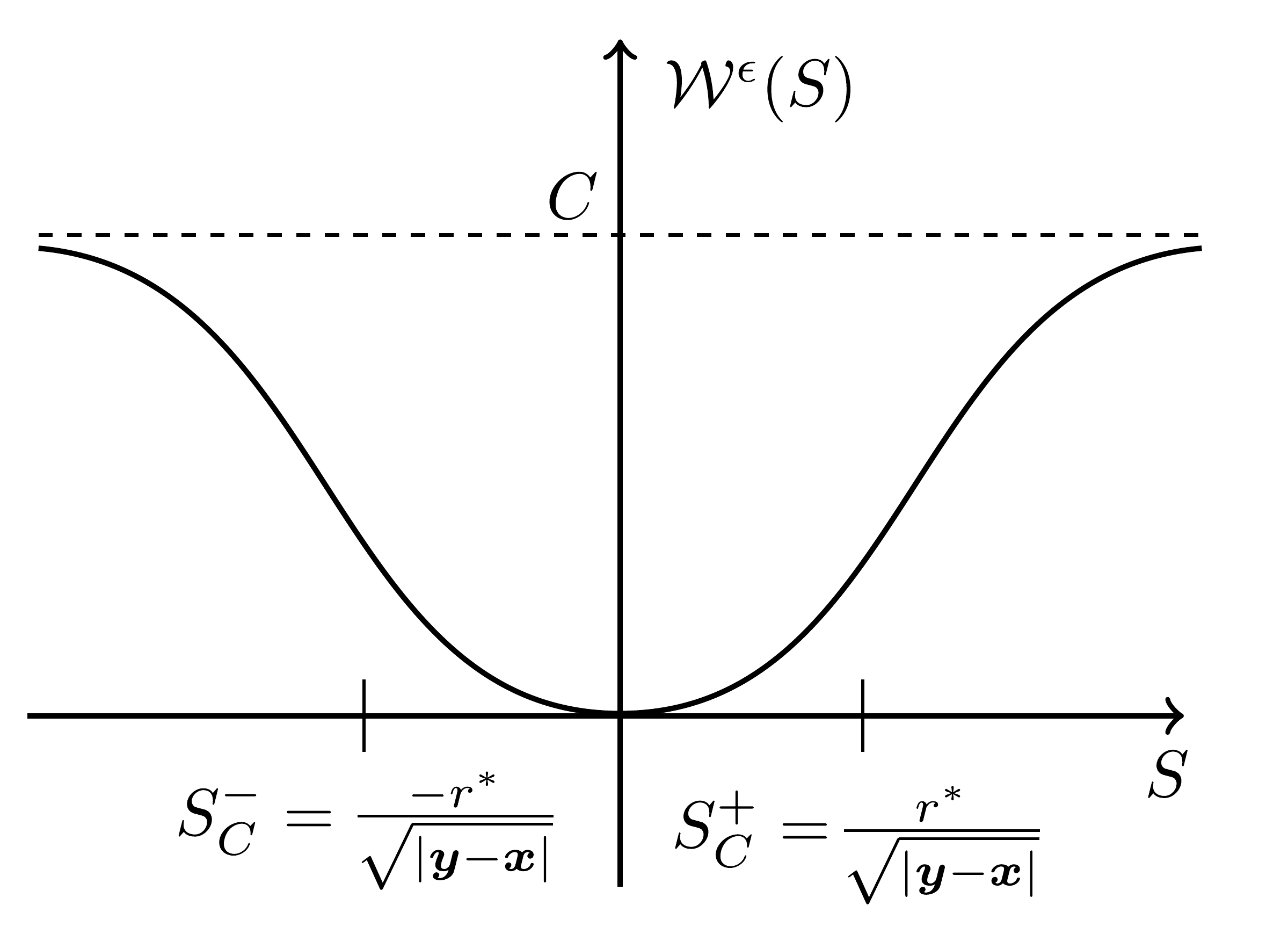}
		\caption{Profile of pairwise potential $\mathcal{W}^\epsilon(S)$ defined in \eqref{eq:rnppotential}. Here $C = \lim_{S \to \infty}\mathcal{W}^\epsilon(S)$.}\label{fig:rnppotential}
	\end{figure}

\begin{figure}[tb]
		\centering
		\begin{subfigure}[t]{0.45\textwidth}
			\centering
			\includegraphics[scale=0.18]{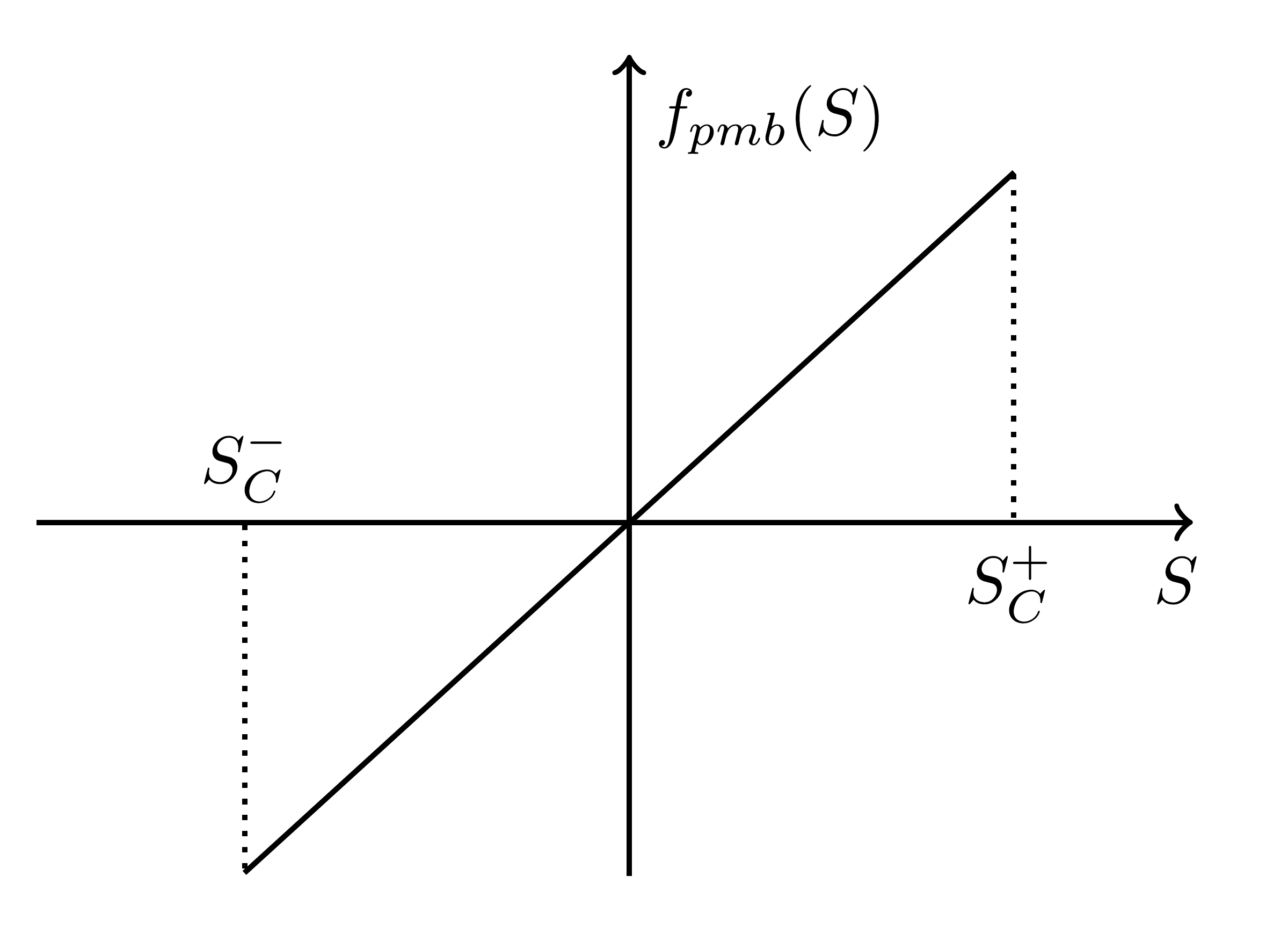}
			\caption{}
			\label{fig:pmb}
		\end{subfigure}
		\begin{subfigure}[t]{0.45\textwidth}
			\centering
			\includegraphics[scale=0.18]{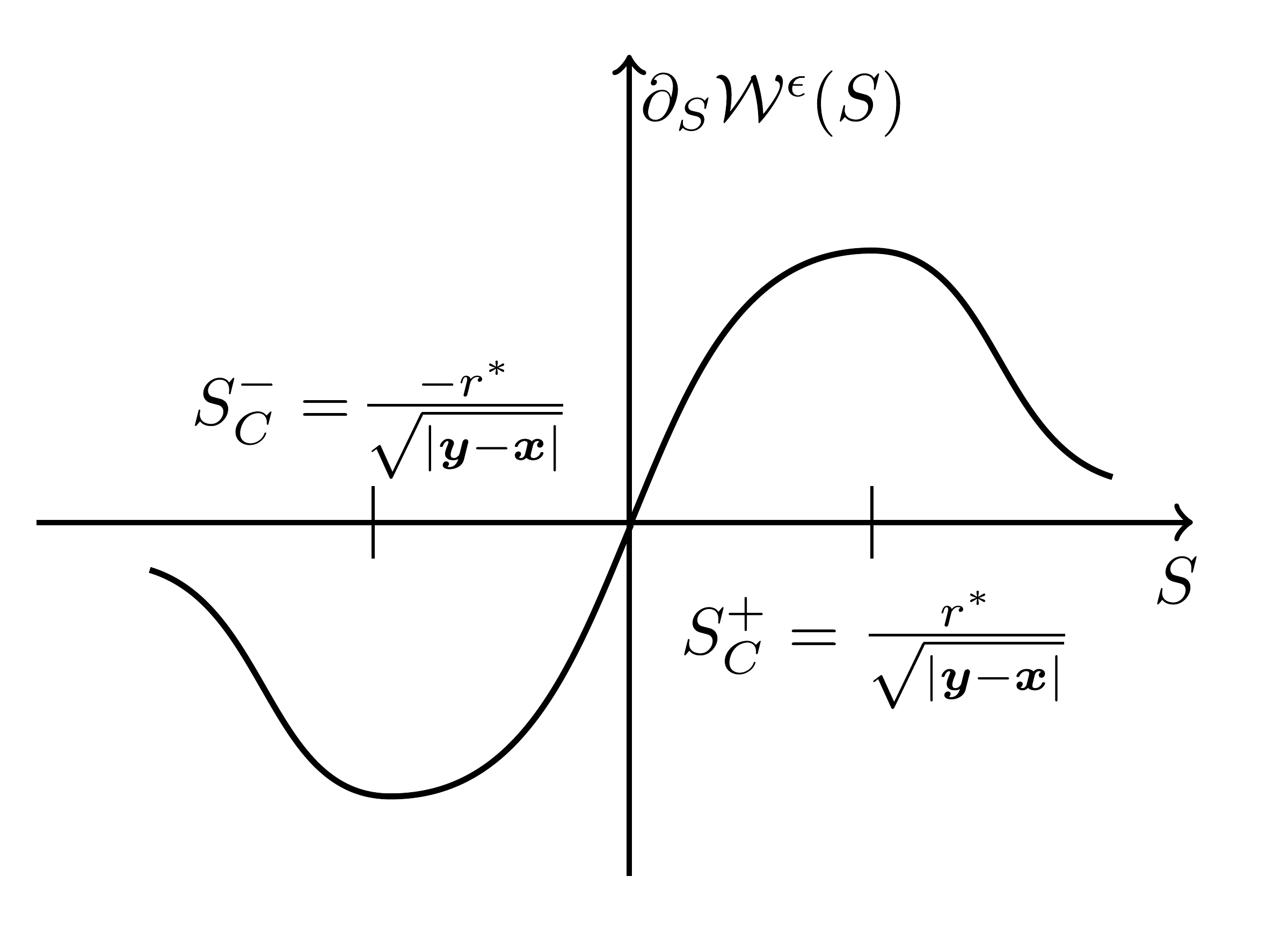}
			\caption{}
			\label{fig:rnp}
		\end{subfigure}
		\caption{(a) Prototype microelastic brittle (PMB) material. Here, $f_{pmb}(S)$ is a scalar such that the bond force is $\bff_{pmb}(\by, \bx, t;\bu) = f_{pmb}(S)\frac{\by - \bx}{|\by - \bx|}$; see \eqref{eq:pmbbondforcesmall}. (b) Regularized nonlinear peridynamics (RNP) material.}
	\end{figure}

	In the PMB model, the interaction between two material points comes to an abrupt stop as soon as the pairwise strain exceeds the critical strain. In contrast, pairwise force considered in~\cite{CMPer-Lipton, CMPer-Lipton3} regularizes the pairwise strain-force profile such that the bond under small strains behaves like a linear elastic material, and for larger strains, yields and softens with increasing strain, and eventually, the bond breaks for large strains. The force model introduced in~\cite{CMPer-Lipton, CMPer-Lipton3} is referred to as the regularized nonlinear peridynamics (RNP) material model. The pairwise potential -- force given by the derivative of the potential -- in the RNP model is defined by
	\begin{equation}\label{eq:rnppotential}
		\mathcal{W}^\epsilon(S(\by, \bx, t)) = \frac{1}{\mathrm{w}_d \epsilon^d}\omega(\bx) \omega(\by) \frac{J^\epsilon(|\by - \bx|)}{\epsilon|\by - \bx|} \psi(|\by - \bx| S(\by, \bx, t)^2).
	\end{equation}
	Here $\mathrm{w}_d = \vert\{\by\in \bbR^d: |\by| < 1\}\vert = \vert H_1(\bzero)\vert$ is the volume of a unit ball in the dimension $d$, i.e. $\mathrm{w}_d = \pi$ in 2-d and $\mathrm{w}_d = 4\pi/3$ in 3-d. $\omega : D \to [0,1]$ is a boundary function which takes the value $1$ for all $\bx \in D_i:= \{ \by\in D: \mathrm{dist}(\by, \partial D) > \epsilon \}$ and decays smoothly from $1$ to $0$ as $\bx$ approaches the boundary $\partial D$. The potential function $\psi: \bbR^+ \to \bbR$ is smooth, positive, and concave. For such a choice of $\psi$, the profile of potential $\mathcal{W}^\epsilon$ as a function of strain $S$ is shown in~\cref{fig:rnppotential}. The pairwise force is written as (see~\cite{CMPer-Lipton, CMPer-Lipton3})
	\begin{equation}\label{eq:rnpbondforce}
		\begin{split}
			\bff_{rnp}(\by, \bx, t) &= 2\,\partial_{S} \mathcal{W}^\epsilon(S(\by, \bx, t)) \,\frac{\by - \bx}{|\by - \bx|} \\
			&= \frac{4}{\mathrm{w}_d \epsilon^d} \omega(\bx) \omega(\by) \frac{J^\epsilon(|\by - \bx|)}{\epsilon} \psi'(|\by - \bx| S(\by, \bx, t)^2)  S(\by,\bx,t) \frac{\by - \bx}{|\by - \bx|}.
		\end{split}
	\end{equation}
	The critical strain $S_c$ depends on the distance between material points $\by$ and $\bx$, and it is given by $S_c(\by, \bx) = \pm r^\ast/{\sqrt{|\by - \bx|}}$. $r^\ast>0$ is the inflection point of function $r \mapsto \psi(r^2)$. The total force at $\bx$ is given by
	\begin{equation}\label{eq:rnptotalforce}
		\bF_{rnp}(\bx, t) = \int_{H_\epsilon(\bx)\cap D} \bff_{rnp}(\by, \bx, t) d\by.
	\end{equation}
	\cref{fig:pmb} and \cref{fig:rnp} shows PMB and RNP force profiles. The RNP model is amenable to a-priori convergence rate analysis and is investigated in this paper.
	
	\subsection{Peridynamics equation of motion using the RNP model}
	
	In the rest of the article, the pairwise strain $S$ defined in \eqref{eq:bondstrainsmall} is considered, and the RNP model is employed where the pairwise force is given by \eqref{eq:rnpbondforce}. The peridynamics equation of motion for the displacement field $\bu : D \times [0,t_F] \to \bbR^d$ is given by the Newton's second law as follows
	\begin{equation}\label{eq:equationofmotion}
		\rho \partial^2_{tt} \bu(\bx,t) = \bF(\bu)(\bx,t) + \bb(\bx,t),
	\end{equation}
	where, $\rho $ is the density, $\bF(\bu)(\bx,t)$ peridynamics force defined in \eqref{eq:rnptotalforce}, and $\bb(\bx,t)$ the body force per unit volume. Let $\partial D$ be the boundary of the material domain $D$. For the analysis, the Dirichlet boundary condition is assumed, \emph{i.e.},
	\begin{equation}\label{eq:bc}
		\bu(\bx,t) = 0, \qquad \forall \bx \in \partial D, \forall t\in [0,t_F].
	\end{equation}
	The initial conditions for displacement and velocities are 
	\begin{equation}\label{eq:initialvalues}
		\bu(\bx,0) = \bu_0(\bx) \qquad \text{and} \qquad \partial_t \bu(\bx,0)= \bv_0(\bx), \forall \bx \in D.
	\end{equation}
	In the rest of the article, density $\rho$ is assumed to be constant. 
	
	For the RNP model, the initial boundary value problem given by \eqref{eq:equationofmotion} with~\eqref{eq:bc} and~\eqref{eq:initialvalues} for $\bu_0, \bv_0 \in H^2(D) \cap H^1_0(D)$ and $\bb \in C^2([0,t_F]; H^2(D)\cap H^1_0(D))$, is shown to be well-posed in the space $C^2([0,t_F];H^2(D) \cap H^1_0(D))$; see \cite[Theorem 3.2]{jha2021finite}. Here, $H^1_0(D)$ is given by the space of functions in $H^1(D)$ taking value zero on the boundary $\partial D$. In what follows, $||\cdot ||$ and $||\cdot ||_n$ will denote the $L^2(D)$ and $H^n(D)$ norms, for $n=1,2$, respectively.

	%
	\section{Finite element approximation}\label{s:fespace}
	Consider a discretization $\mathcal{T}_h$ of the domain $D$ by triangular (in 2-d) or tetrahedral (in 3-d) elements, where $h$ denotes the size of mesh assuming that the elements are conforming and the mesh is shape regular. Let $\bar{V}_h$ and $V_h$, with $V_h \subset \bar{V}_h$, denote the spaces of functions spanned by continuous piecewise linear interpolation over mesh $\mathcal{T}_h$ such that $\bar{V}_h\subset H^1(D)$ and $V_h \subset H^1_0(D)$. It is further assumed that there exist constants $c_1,c_2> 0$ such that
	\begin{equation}\label{eq:l2L2relation}
		c_1 \sum_{i=1}^N |\bu(\bx_i)|^2 \leq ||\bu||^2 \leq c_2 \sum_{i=1}^N |\bu(\bx_i)|^2, \qquad \forall \bu \in \bar{V}_h ,
	\end{equation}
	where, $N$ is the total number of mesh nodes, and $\bx_i$ is the material coordinate of $i^\text{th}$ node.
	
	For a continuous function $\bu$ on $\bar{D}$, it's continuous piecewise linear interpolant on $\mathcal{T}_h$ is defined as
	\begin{equation}
		\mathcal{I}_h(\bu)\bigg\vert_{T} = \mathcal{I}_T(\bu) \qquad \forall T\in \mathcal{T}_h,
	\end{equation}
	where, $\mathcal{I}_T(\bu)$ is the local interpolant associated with the finite element $T$ such that
	\begin{equation}
		\label{eq:Interpolant}
		\mathcal{I}_T(\bu) = \sum_{i\in N_T} \bu(\bx_i^T)\phi_i^T.
	\end{equation}
	Here, $N_T$ is the list of nodes as a vertex of the element $T$, $\bx_i^T$ the position of $i^{\text{th}}$ vertex of the element $T$, and $\phi_i^T$ the linear interpolation function associated with the vertex $i$. 
	
	Application of Theorem 4.4.20 and Remark 4.4.27 in~\cite{MANa-Susanne} gives the bound on the interpolation error in $L^2$ norm as follows
	\begin{equation}\label{eq:interpolationerror}
		|| \bu - \mathcal{I}_h(\bu) || \leq c_3 h^2 || \bu ||_2, \qquad\qquad \forall \bu \in H^2(D)^2
	\end{equation}
	and in $L^\infty$ norm
	\begin{equation}\label{eq:interpolationerrorpointwise}
		\sup_{\bx}|\bu(\bx) - \mathcal{I}_h(\bu)(\bx)| \leq c_4 h^2 \sup_{\bx\in D} \left\vert\nabla^2 \bu(\bx)\right\vert, \qquad \qquad \forall \bu \in C^2(D)^2.
	\end{equation}
	Here, constants $c_3,c_4$ are independent of mesh size $h$.
	
	\paragraph{Projection onto $V_h$}Let $\br_h(\bu)$ denote the projection of $\bu\in H^1_0(D)$ on $V_h$ with respect to the $L^2$ norm. It is defined by
	\begin{equation}
		||\bu - \br_h(\bu)|| = \inf_{\tilde{\bu}\in V_h} ||\bu - \tilde{\bu}||  \label{eq:proju}
	\end{equation}
	and satisfies the orthogonality property
	\begin{equation}\label{eq:projorthogonalprop}
		(\br_h(\bu), \tilde{\bu}) = (\bu, \tilde{\bu}), \qquad \forall \tilde{\bu} \in V_h.
	\end{equation}
	
	Since $\mathcal{I}_h(\bu) \in V_h$, from \eqref{eq:interpolationerror} it follows that
	\begin{equation}
		||\bu - \br_h(\bu)|| \leq ||\bu - \mathcal{I}_h(\bu)|| \leq c_3 h^2 || \bu ||_2,\hbox{     } \qquad \forall \bu \in H^2(D). \label{eq:projerror}
	\end{equation}
	
	\paragraph{Cl\'ement interpolation}
	The Cl\'ement interpolant \cite{clement1975} is a linear operator $\mathcal{I}_h^C:\,L^2(\Omega)^2\rightarrow V_h$ and to define it, let $S_i$ denote the set of elements $T$ with a common vertex $i$. Let $\mathcal{P}_1(S_i)$ be the space of continuous piecewise
	linear finite elements on $S_i$. Denote $P_i:\, L^2(S_i)^2\rightarrow \mathcal{P}_1(S_i)$ as the $L^2$-projection. Given $\bu\in L^2(\Omega)^2$ 
	\begin{equation}
		\label{clement interpolant}
		\mathcal{I}_T^C(\bu) = \sum_{i\in N_T} \mathcal{P}_i(\bu)(\bx_i^T)\phi_i^T,
	\end{equation}
	the bound on the interpolation error is
	\begin{equation}\label{eq:Clemtnt interpolationerror}
		|| \bu - \mathcal{I}_h^C(\bu) ||_t \leq C h^{s-t} || \bu ||_s, \qquad\qquad \forall \bu \in H^s(D)^2,
	\end{equation}
	for $0\leq t\leq s\leq 2$ and $t\leq 1$, where $t=0$ corresponds to the $L^2$ norm.
	Here, the constant $C$ is independent of mesh size $h$. On writing $\mathcal{I}_h^C(\bu)+\bu-\bu$ together with the triangle inequality delivers the stability
	\begin{equation}\label{eq:Clemtnt stab}
		||\mathcal{I}_h^C(\bu) ||_t \leq C h^{s-t} || \bu ||_s, \qquad\qquad \forall \bu \in H^s(D)^2,
	\end{equation}

	\subsection{Nodal finite element approximation}\label{ss:nodalfe}
	Let $\Delta t$ be the size of the time step and $t^k = k \Delta t$ be the time at step $k$. Let $\{\bU_i^k\}_{1\leq i\leq N}$ be the set of approximate nodal displacements at time step $k$. Associated to the discrete set $\{\bU_i^k\}_{1\leq i\leq N}$, displacement field $\bu^k_h \in V_h$ can be defined as follows
	\begin{equation}\label{eq:dispstrong}
		\bu_h^k(\bx)\bigg\vert_{T} = \sum_{i\in N_T} \bU_i^k\phi_i(\bx),\qquad \bx\in T, \qquad \forall T\in \mathcal{T}_h.
	\end{equation}

	The discrete solution satisfies, for all $1\leq i\leq N$ and for $k\geq 1$,
	\begin{equation}\label{eq:centralstrong}
		\rho\dfrac{\bU^{k+1}_i - 2 \bU^k_i + \bU^{k-1}_i}{\Delta t^2} = \bF(\bu^k_h)(\bx_i) + \bb(\bx_i, t^k)
	\end{equation}
	and, for $k = 0$ (first time step),
	\begin{equation}\label{eq:centralzerostrong}
		\rho\dfrac{\bU^1_i - \bu_0(\bx_i)}{\Delta t^2} = \dfrac{1}{2}\bF(\bu^0_h)(\bx_i) + \dfrac{1}{\Delta t} \bv_0(\bx_i) + \dfrac{1}{2}\bb(\bx_i, 0).
	\end{equation}
	In the above, $\bu^0$ and $\bv^0$ are the initial conditions. 
	
	Fix $\bx\in T$, where $T\in \mathcal{T}_h$. By multiplying interpolation function $\phi_i(\bx)$ to both sides of~\eqref{eq:centralstrong} and summing over $i\in N_T$, it follows that
	\begin{equation}\label{eq:centralstrong2.1}
		\rho\dfrac{\bu^{k+1}_h(\bx) - 2 \bu^k_h(\bx) + \bu^{k-1}_h(\bx)}{\Delta t^2} \bigg\vert_T = \sum_{i\in N_T}\bF(\bu^k_h)(\bx_i)\phi_i(\bx) + \sum_{i\in N_T}\bb(\bx_i, t^k)\phi_i(\bx),
	\end{equation}
	Let $\bF_h(\bu^k_h)$ and $\bb^k_h$ be the continuous piecewise linear interpolation of $\bF(\bu^k_h)$ and $\bb(t^k)$, \emph{i.e.},
	\begin{equation}
		\begin{split}
			\bF_h(\bu^k_h)(\bx)\bigg\vert_{T} &= \sum_{i\in N_T} \bF(\bu^k_h)(\bx_i) \phi_i(\bx), \qquad \bx\in T, \quad \forall T\in \mathcal{T}_h, \notag \\
			\bb^k_h(\bx)\bigg\vert_{T} &= \sum_{i\in N_T} \bb(\bx_i, t^k) \phi_i(\bx),\qquad \bx\in T, \quad \forall T\in \mathcal{T}_h.
		\end{split}
	\end{equation}
	Then, \eqref{eq:centralstrong2.1} can be written compactly as
	\begin{equation}\label{eq:centralstrong2}
		\rho\dfrac{\bu^{k+1}_h - 2 \bu^k_h + \bu^{k-1}_h}{\Delta t^2} = \bF_h(\bu^k_h) + \bb^k_h.
	\end{equation}

	\subsection{Comparison of NFEA with the standard FEA}\label{ss:compareFE}
	Let $\hat{\bu}^k_h\in V_h$ be the standard FEA solution that satisfies (see \cite{jha2021finite}), for all test functions $\tilde{\bu}\in V_h$ and $k\geq 1$,
	\begin{equation}\label{eq:stdFEA}
		\left( \rho\dfrac{\hat{\bu}^{k+1}_h - 2 \hat{\bu}^k_h + \hat{\bu}^{k-1}_h}{\Delta t^2}, \tilde{\bu} \right) = (\bF(\hat{\bu}^k_h), \tilde{\bu} ) + (\bb(t^k), \tilde{\bu}).
	\end{equation}
	To see the difference between the above discretization and the nodal FEA, multiply~\eqref{eq:centralstrong2} with the test function $\tilde{\bu}\in V_h$ and integrate over a domain $D$ to have
	\begin{equation}\label{eq:centralstrong3}
		\left( \rho\dfrac{\bu^{k+1}_h - 2 \bu^k_h + \bu^{k-1}_h}{\Delta t^2}, \tilde{\bu} \right) = (\bF_h(\bu^k_h), \tilde{\bu} ) + (\bb^k_h, \tilde{\bu}).
	\end{equation} 
	Thus, in the NFEA, the exact peridynamics force $\bF$ and body force $\bb$ are replaced by their continuous piecewise linear interpolation $\bF_h$ and $\bb_h$, respectively. By doing so, NFEA reduces the computational complexity of computing the integral of the product of peridynamics force and test function in~\eqref{eq:stdFEA} but at the cost of an additional discretization error; compare $(\bF(\bu^k_h), \tilde{\bu} )$ and $(\bb^k, \tilde{\bu})$ in~\eqref{eq:stdFEA} with $(\bF_h(\bu^k_h), \tilde{\bu} )$ and $(\bb^k_h, \tilde{\bu})$ in~\eqref{eq:centralstrong2}, respectively. 
	
	Next, a-priori convergence of NFEA solution $\bu^k_h$ to the exact solution $\bu(t^k)$ in the limit mesh size, $h$, and time step, $\Delta t$, tending to zero is shown. 
	
	%
	\section{A-priori convergence of nodal FEA for nonlinear peridynamics models}\label{s:convergence}
	This section establishes the convergence of the NFEA approximation to the exact peridynamics solution. The error analysis is focused on nonlinear peridynamics force (RNP), see \eqref{eq:rnptotalforce}, and exact solutions (displacement and velocity) are assumed to be $\bu, \bv \in C^2\left([0,t_F]; H^2(D)\cap H^1_0(D)\right)$. Before the main result is presented, equations for errors are obtained, and the consistency of the numerical discretization is shown. 
	
	Equation \eqref{eq:centralstrong3} is equivalent to \eqref{eq:centralstrong2} and it can be decoupled into two equations by introducing $\bv^k_h \in V_h$, $\bv^k_h$ being the approximate velocity field at time $t^k$, as follows
	\begin{equation}\label{eq:forward}
		\begin{split}
			\left( \dfrac{\bu^{k+1}_h - \bu^k_h}{\Delta t}, \tilde{\bu} \right) &= (\bv^{k+1}_h, \tilde{\bu}), \qquad \qquad \forall \tilde{\bu} \in V_h, \\
			\left( \dfrac{\bv^{k+1}_h - \bv^k_h}{\Delta t}, \tilde{\bu} \right) &= (\bF_{rnp, h}(\bu^k_h), \tilde{\bu} ) + (\bb^k_h, \tilde{\bu}), \qquad \qquad \forall \tilde{\bu} \in V_h .
		\end{split}
	\end{equation}
	
	Similar to \cite[Section 5]{jha2021finite}, $L^2$ projections of the exact displacement and velocity into $V_h$ are compared with the approximate displacement and velocity, and the errors are defined as
	\begin{equation*}
		\be^k_{u,h} = \bu^k_h - \br_h(\bu^k), \qquad \be^k_{v,h} = \bv^k_h -\br_h(\bv^k),
	\end{equation*}
	where $\bu^k$ is the exact solution at time $t^k = k\Delta t$, $\bv^k = \partial_t \bu(t^k)$, and $\br_h(\bu) \in V_h$ is the $L^2$ projection of $\bu$ defined in~\eqref{eq:proju}. Using the peridynamics equation of motion~\eqref{eq:equationofmotion}, \eqref{eq:forward}, and property~\eqref{eq:projorthogonalprop} of projection $\br_h$, it can be shown that
	\begin{align}
		( \be^{k+1}_{u,h}, \tilde{\bu}) &= (\be^{k}_{u,h}, \tilde{\bu}) + \Delta t (\be^{k+1}_{v,h}, \tilde{\bu} ) + \Delta t (\btau^k_{u,h},\tilde{\bu}) , \label{eq:eku} \\
		(\be^{k+1}_{v,h}, \tilde{\bu}) &= (\be^k_{v,h}, \tilde{\bu}) + \Delta t (\bF_{rnp}(\bu^k) - \bF_{rnp,h}(\bu^k_h), \tilde{\bu}) + \Delta t (\btau^k_{v,h}, \tilde{\bu}) + \Delta t (\bb^k_h - \bb(t^k), \tilde{\bu}), \label{eq:ekv}
	\end{align}
	where, $\btau^k_{u,h}, \btau^k_{v,h}$ are consistency errors and take the form
	\begin{equation}\label{eq:deftauk}
		\begin{split}
			\btau^k_{u,h} &= \dfrac{\partial \bu^{k+1}}{\partial t} - \dfrac{\bu^{k+1} - \bu^k}{\Delta t}, \\
			\btau^k_{v,h} &= \dfrac{\partial \bv^k}{\partial t} - \dfrac{\bv^{k+1} - \bv^k}{\Delta t}. 
		\end{split}
	\end{equation}
	
	\subsection{Key estimates}\label{ss:keyest}
	This section estimates the error terms in \eqref{eq:eku} and \eqref{eq:ekv}. In this direction, note that, if $\bu, \bv \in C^2([0,t_F];H^2(D))$, then 
	\begin{equation}\label{eq:consistenctytime}
		||\btau^k_{u,h}|| + ||\btau^k_{v,h}|| \leq C_t \Delta t, \qquad\qquad C_t = \sup_t ||\partial^2_{tt}\bu(t)|| + \sup_t ||\partial^3_{ttt}{\bu}(t)||.
	\end{equation}
	Further, if $\bb \in C([0,t_F]; H^2(D))$ then noting that $\bb^k_h$ is a linear interpolation of $\bb(t^k)$ it can be easily shown using \eqref{eq:interpolationerror} that
	\begin{equation}\label{eq:consistencybodyforce}
		||\bb^k_h - \bb(t^k)|| \leq c_3h^2 \sup_t ||\bb(t)||_2.
	\end{equation}
	
	Focusing on the remaining consistency error term in \eqref{eq:ekv}, $\bF_{rnp}(\bu^k) - \bF_{rnp,h}(\bu^k_h)$, using the triangle inequality, the error can be shown to be bounded by the sum of the four terms as follows:
	\begin{equation}\label{eq:consistency force nodal}
		\begin{split}
			&||\bF_{rnp}(\bu^k) - \bF_{rnp,h}(\bu^k_h)||  \\
			&\leq \underbrace{||\bF_{rnp}(\bu^k) - \bF_{rnp,h}(\bu^k)||}_{=: I_1} + \underbrace{||\bF_{rnp,h}(\bu^k) - \bF_{rnp,h}(\mathcal{I}_h(\bu^k))||}_{=:I_2} \\
			&\quad + \underbrace{||\bF_{rnp,h}(\mathcal{I}_h(\bu^k)) - \bF_{rnp,h}(\br_h(\bu^k))||}_{=:I_3} + \underbrace{||\bF_{rnp,h}(\br_h(\bu^k)) - \bF_{rnp,h}(\bu^k_h)||}_{=:I_4} .
		\end{split}
	\end{equation}
	To bound the above terms, the following property of nonlinear peridynamics force $\bF_{rnp}$ is crucial. 
	
	\begin{rem}
		Assuming that the domain $D$ is a $C^1$ domain, the boundary function $\omega \in C^2(D)$, and the peridynamics potential $\psi$ (see~\eqref{eq:rnppotential} or~\eqref{eq:rnpbondforce}) is smooth with up to 4$^\text{th}$ order bounded derivatives, from \cite[Section 3]{jha2021finite}, it holds that
		\begin{equation}\label{eq:bddrnpforceH2}
			||\bF_{rnp}(\bu)||_2 \leq \frac{L_2 ||\bu||_2 + L_3||\bu||_2^2}{\epsilon^{5/2}}, \qquad \forall \bu \in H^2(D).
		\end{equation}
		Further, the peridynamics force satisfies the following Lipschitz continuity condition in the $L^2$ norm
		\begin{equation}\label{eq:lipschitzL2}
			||\bF_{rnp}(\bu) - \bF_{rnp}(\bv)|| \leq \frac{L_1}{\epsilon^2} ||\bu - \bv|| \qquad \forall \bu,\bv \in L^2(D).
		\end{equation}
		Here, constants $L_1, L_2, L_3$ are independent of $\bu,\bv$ and depend on the influence function $J$ and peridynamics force potential $\psi$. For future reference, $L_1$, from \cite[Section 3]{jha2021finite}, is given by
		\begin{equation}\label{eq:defL1}
			L_1 := 4 \underbrace{\left(\frac{1}{\mathrm{w}_d} \int_{H_1(\bzero)} \frac{J(|\bxi|)}{|\bxi|} d\bxi\right)}_{=: \bar{J}_1} \;  \underbrace{ \left(\sup_r \bigg|\frac{d^2}{dr^2} \psi(r^2) \bigg| \right)}_{=: C_2} = 4 \bar{J}_1 C_2.
		\end{equation}
	\end{rem}
	
	The lemma below collects the bounds on the errors $I_n$, $n=1,2,3,4$.
	\begin{lem}\label{lemma:nodal consistency}
		\textbf{Consistency of the peridynamics force}\\
		For $\bu^k$ in $H^2(D)\cap C^2(D)$, the following estimates hold
		\begin{equation}\label{eq:consistencyspatial}
			\begin{split}
				I_1 &\leq \left[c_3 \frac{L_2 ||\bu^k||_2 + L_3 ||\bu^k||_2^2}{\epsilon^{5/2}}\right] h^2,  \qquad I_2 \leq \left[\frac{L_1 c_4 \sqrt{\frac{c_2|D|}{c_1}}}{\epsilon^2} \sup_{\bx\in D} \left\vert\nabla^2 \bu^k(\bx)\right\vert \right] h^2, \\
				I_3 &\leq \left[\frac{2 c_3 L_1 \bar{n} \sqrt{\frac{c_2}{c_1}}}{\epsilon^2} ||\bu^k||_2  \right] h^2, \quad\qquad \;
				I_4 \leq \left[\frac{L_1 \bar{n} \sqrt{\frac{c_2}{c_1}}}{\epsilon^2}   \right] ||\be^k_{u,h}||.
			\end{split}
		\end{equation}
		Here, $c_i$, $i=1,2,3,4$, are constants depending only on the triangulation $\mathcal{T}_h$, see~\eqref{eq:interpolationerror}, \eqref{eq:interpolationerrorpointwise}, \eqref{eq:l2L2relation}. Moreover, $L_i$, $i=1,2,3$, are constants that only depend on the influence function $J$ and the peridynamics force potential function $\psi$. Finally, the constant $\bar{n}$ is given by
		\begin{equation}
			\bar{n} = \max_{T\in \mathcal{T}_h} \{ \text{number of vertices of element }T \}.
		\end{equation} 
		
	\end{lem}
	
	\paragraph{Proof}
	Let us consider $I_1$ first. Since $\bu^k \in H^2(D) \cap C^2(D)$, note that $\bF_{rnp,h}(\bu^k) = \mathcal{I}_h(\bF_{rnp}(\bu^k))$, where $\mathcal{I}_h$ is the continuous piecewise linear interpolant. Using \eqref{eq:interpolationerror} and \eqref{eq:bddrnpforceH2}, it can be shown that
	\begin{equation}
		\begin{split}
			I_1 &= ||\bF_{rnp}(\bu^k) - \mathcal{I}_h(\bF_{rnp}(\bu^k))|| \leq c_3 h^2 ||\bF_{rnp}(\bu^k)||_2 \notag \\
			&\leq \left[c_3\frac{L_2 ||\bu^k||_2 + L_3 ||\bu^k||_2^2}{\epsilon^{5/2}} \right] h^2.
		\end{split}
	\end{equation}
	
	Next, $I_3$ and $I_4$ are bounded from above. Let $\bw_1, \bw_2 \in V_h$, then both $I_3$ and $I_4$ are of the form $||\bF_{rnp,h}(\bw_1) - \bF_{rnp,h}(\bw_2)||$. Now, using the definition of $\bF_{rnp}$ in~\eqref{eq:rnptotalforce}, it follows that
	\begin{equation}\label{eq:diffrnpforce}
		\begin{split}
			&\bF_{rnp}(\bw_1)(\bx_i) - \bF_{rnp}(\bw_2)(\bx_i) \\
			&= \frac{4}{\mathrm{w}_d \epsilon^{d+1}} \int_{H_\epsilon(\bx) \cap D} \omega(\bx_i) \omega(\by) J^\epsilon(|\by - \bx_i|) \\
			&\quad  [\psi'(|\by - \bx_i| S(\by, \bx_i; \bw_1)^2)S(\by, \bx_i; \bw_1) - \psi'(|\by - \bx_i| S(\by, \bx_i; \bw_2)^2)S(\by, \bx_i; \bw_2)] \frac{\by - \bx_i}{|\by - \bx_i|} d\by.
		\end{split}
	\end{equation}
	Let $\Psi(r) := \psi(r^2)$ then $\Psi'(r) = 2r \psi'(r^2)$ and $|\Psi'(r_1) - \Psi'(r_2)| \leq \sup_r |\Psi''(r)| |r_1 - r_2|$. Since $\psi$ is smooth and has up to 4 bounded derivatives, $\sup_r |\Psi''(r)| = \sup_r |\psi(r^2)| = C_2 < \infty$. Using the constant $C_2$, it holds that
	\begin{equation*}
		\begin{split}
			&2|\psi'(|\by - \bx_i| S(\by, \bx_i; \bw_1)^2)S(\by, \bx_i; \bw_1) - \psi'(|\by - \bx_i| S(\by, \bx_i; \bw_2)^2)S(\by, \bx_i; \bw_2)| \\
			&\leq C_2 \sqrt{|\by-\bx_i|} \, |S(\by, \bx_i; \bw_1) - S(\by, \bx_i; \bw_2)| \\
			&= C_2 \sqrt{|\by-\bx_i|}\, \bigg|\frac{\bw_1(\by) - \bw_1(\bx_i)}{|\by - \bx_i|} \cdot \frac{\by - \bx_i}{|\by - \bx_i|} - \frac{\bw_2(\by) - \bw_2(\bx_i)}{|\by - \bx_i|} \cdot \frac{\by - \bx_i}{|\by - \bx_i|} \bigg| \\
			&\leq C_2 \frac{|\bw_1(\by) - \bw_2(\by)| + |\bw_1(\bx_i) - \bw_2(\bx_i)|}{\sqrt{|\by-\bx_i|}}.
		\end{split}
	\end{equation*}
	Using the above bound and change in variable $\bxi = (\by - \bx_i)/\epsilon \in H_1(\bzero)$, from~\eqref{eq:diffrnpforce}, one can show that
	\begin{equation}\label{eq:diffrnpest1}
		\begin{split}
			&|\bF_{rnp}(\bw_1)(\bx_i) - \bF_{rnp}(\bw_2)(\bx_i)| \\
			&\leq \frac{2C_2}{\epsilon^2 \mathrm{w}_d} \int_{H_1(\bzero)} \frac{J(|\bxi|)}{\sqrt{|\bxi|}} (|\bw_1(\bx_i + \epsilon \bxi) - \bw_2(\bx_i + \epsilon \bxi)|+|\bw_1(\bx_i) - \bw_2(\bx_i)|) d\bxi
		\end{split}
	\end{equation}
	and 
	\begin{equation}\label{eq:diffrnpest2}
		\begin{split}
			&|\bF_{rnp}(\bw_1)(\bx_i) - \bF_{rnp}(\bw_2)(\bx_i)|^2 \\
			&\leq 2\left(\frac{2C_2}{\epsilon^2 \mathrm{w}_d}\right)^2 \frac{\bar{J}_1}{\mathrm{w}_d} \int_{H_1(\bzero)} \frac{J(|\bxi|)}{|\bxi|} (|\bw_1(\bx_i + \epsilon \bxi) - \bw_2(\bx_i + \epsilon \bxi)|^2+|\bw_1(\bx_i) - \bw_2(\bx_i)|^2) d\bxi,
		\end{split}
	\end{equation}
	where, $\bar{J}_1 := \frac{1}{\mathrm{w}_d} \int_{H_1(\bzero)} J(|\bxi|)/|\bxi| d\bxi$. 
	
	Next, using the property of the finite element function space $V_h$ that relates $L^2$ norm to discrete $l^2$ norm in~\eqref{eq:l2L2relation}, it can be shown that
	\begin{equation}\label{eq:diffrnpest3}
		\begin{split}
			&||\bF_{rnp,h}(\bw_1) - \bF_{rnp,h}(\bw_2)||^2 \\
			&\leq c_2 \sum_{i=1}^N |\bF_{rnp}(\bw_1)(\bx_i) - \bF_{rnp}(\bw_2)(\bx_i)|^2  \\ 
			&\leq c_2 2\left(\frac{2C_2}{\epsilon^2 \mathrm{w}_d}\right)^2 \frac{\bar{J}_1}{\mathrm{w}_d} \int_{H_1(\bzero)} \frac{J(|\bxi|)}{|\bxi|} \left[ \sum_{i=1}^N (|\bw_1(\bx_i + \epsilon \bxi) - \bw_2(\bx_i + \epsilon \bxi)|^2+|\bw_1(\bx_i) - \bw_2(\bx_i)|^2) \right] d\bxi.
		\end{split}
	\end{equation}
	Since $\bw_1, \bw_2 \in V_h$, using \eqref{eq:l2L2relation}, it holds that
	\begin{equation}\label{eq:diffrnpest4}
		\sum_{i=1}^N |\bw_1(\bx_i) - \bw_2(\bx_i)|^2 \leq \frac{1}{c_1} ||\bw_1 - \bw_2||^2.
	\end{equation}
	Now, to estimate
	\begin{equation*}
		\sum_{i=1}^N |\bw_1(\bx_i + \epsilon \bxi) - \bw_2(\bx_i+\epsilon \bxi)|^2,
	\end{equation*} 
	consider any point $\by \in T$, where $T \in \mathcal{T}_h$. Denoting the set of vertices of an element $T$ as $N_T$, it follows that
	\begin{equation*}
		\begin{split}
			|\bw_1(\by) - \bw_2(\by)|^2 &= |\sum_{i\in N_T} (\bw_1(\bx_i) - \bw_2(\bx_i)) \phi_i(\by)|^2 \\
			&\leq |N_T| \sum_{i\in N_T} |\bw_1(\bx_i) - \bw_2(\bx_i)|^2 |\phi_i(\by)|^2 \leq |N_T|\sum_{i\in N_T} |\bw_1(\bx_i) - \bw_2(\bx_i)|^2,
		\end{split}
	\end{equation*}
	where, in the above, the property of the interpolation function $\phi_i \leq 1$ is used, and $|N_T|$ gives the size of set $N_T$. Let $\bar{n} = \max_{T\in \mathcal{T}_h} |N_T|$, and define the map which returns the element $T$ that contains the point $\by$ by $\Pi(\by)$, i.e., $\Pi(\by) = T$ such that $\by \in \bar{T}$, $\bar{T} = T \cup \partial T$ is the closure of the set $T$. It is assumed that $\Pi$ returns a unique element for all $\by$. Note that for $\by$ on the boundary of an element $T$, $\by$ could belong to more than one element. In such cases, $\Pi$ is assumed to pick one element out of multiple elements randomly or through some selection scheme. It is easy to see now that
	\begin{equation*}
		\sum_{i=1}^N |\bw_1(\bx_i + \epsilon\bxi) - \bw_2(\bx_i + \epsilon\bxi)|^2 \leq  \sum_{i=1}^N \left[ \bar{n} \sum_{j\in N_{\Pi(\bx_i + \epsilon \bxi)}} |\bw_1(\bx_j) - \bw_2(\bx_j)|^2 \right].
	\end{equation*}
	In above double summation, each $|\bw_1(\bx_l) - \bw_2(\bx_l)|^2$ for $l=1,...,N$ will be counted at max $\bar{n}$ times, so
	\begin{align*}
		\sum_{i=1}^N |\bw_1(\bx_i + \epsilon\bxi) - \bw_2(\bx_i + \epsilon\bxi)|^2 &\leq  \sum_{i=1}^N \left[ \bar{n} \sum_{j\in N_{\Pi(\bx_i + \epsilon \bxi)}} |\bw_1(\bx_j) - \bw_2(\bx_j)|^2 \right] \\
		&\leq \bar{n}^2 \sum_{i=1}^N |\bw_1(\bx_i) - \bw_2(\bx_i)|^2.
	\end{align*}
	Combining the above inequality with~\eqref{eq:diffrnpest4}, the following holds, for any $\bxi \in H_1(\bzero)$,
	\begin{equation}\label{eq:diffrnpest5}
		\sum_{i=1}^N |\bw_1(\bx_i + \epsilon\bxi) - \bw_2(\bx_i + \epsilon\bxi)|^2 \leq \frac{\bar{n}^2}{c_1} ||\bw_1 - \bw_2||^2.
	\end{equation}
	
	By combining~\eqref{eq:diffrnpest4} and \eqref{eq:diffrnpest5} with~\eqref{eq:diffrnpest3}, it can be shown that
	\begin{equation*}
		\begin{split}
			&||\bF_{rnp,h}(\bw_1) - \bF_{rnp,h}(\bw_2)||^2 \\
			&\leq 2c_2 \left(\frac{2C_2}{\epsilon^2 \mathrm{w}_d}\right)^2 \frac{\bar{J}_1}{\mathrm{w}_d} \int_{H_1(\bzero)} \frac{J(|\bxi|)}{|\bxi|} \left[\frac{1+\bar{n}^2}{c_1}||\bw_1 - \bw_2||^2 \right] d\bxi \\
			&= \frac{(1+\bar{n}^2) (c_2/c_1)  8C_2^2 \bar{J}_1^2}{\epsilon^4} ||\bw_1 - \bw_2||^2 \\
			&\leq \left[ \frac{\bar{n} \sqrt{c_2/c_1}\, L_1}{\epsilon^2} ||\bw_1 - \bw_2||\right]^2,
		\end{split}
	\end{equation*}
	where, $L_1 = 4C_2 \bar{J}_1$ (see~\eqref{eq:defL1}). Using the above bound that holds for any $\bw_1, \bw_2 \in V_h$, one readily obtains
	\begin{equation}
		\begin{split}
			I_3 &= ||\bF_{rnp,h}(\mathcal{I}_h(\bu^k)) - \bF_{rnp,h}(\br_h(\bu^k))|| \leq \frac{\bar{n} \sqrt{c_2/c_1}\, L_1}{\epsilon^2} ||\mathcal{I}_h(\bu^k) - \br_h(\bu^k)|| \\
			&\leq \frac{\bar{n} \sqrt{c_2/c_1}\, L_1}{\epsilon^2} \left[||\mathcal{I}_h(\bu^k) - \bu^k|| + ||\bu^k - \br_h(\bu^k)|| \right] \leq \frac{\bar{n} \sqrt{c_2/c_1}\, L_1}{\epsilon^2} 2c_3 h^2 ||\bu^k||_2,
		\end{split}
	\end{equation}
	where, \eqref{eq:interpolationerror} and \eqref{eq:projerror} are utilized in the last step. Similarly, $I_4$ can be bounded from the above as follows
	\begin{equation}\label{eq: previous timesteperror}
		I_4 = ||\bF_{rnp,h}(\br_h(\bu^k)) - \bF_{rnp,h}(\bu^k_h)|| \leq  \frac{\bar{n} \sqrt{c_2/c_1}\, L_1}{\epsilon^2} ||\br_h(\bu^k) - \bu^k_h||   = \frac{\bar{n} \sqrt{c_2/c_1}\, L_1}{\epsilon^2} ||\be^k_{u,h}||,
	\end{equation}
	where, the definition of the error $\be^k_{u,h}$ is used in the last step. 
	
	Next, $I_2$ is bounded from the above. Bounds established so far only used the fact that $\bu^k \in H^2(D)$. However, to bound $I_2$, additional regularity of $\bu^k$, $\bu^k \in H^2(D) \cap C^2(D)$, will be utilized. 
	
	\begin{equation}
		\begin{split}
			I_2^2 = &||\bF_{rnp,h}(\bu^k) - \bF_{rnp,h}(\mathcal{I}_h(\bu^k))||^2 
			\leq 2c_2 \left(\frac{2C_2}{\epsilon^2 \mathrm{w}_d}\right)^2 \frac{\bar{J}_1}{\mathrm{w}_d} \\
			\quad \int_{H_1(\bzero)} \frac{J(|\bxi|)}{|\bxi|}
			&\left[ \sum_{i=1}^N (|\bu^k(\bx_i + \epsilon \bxi) - \mathcal{I}_h(\bu^k)(\bx_i + \epsilon \bxi)|^2+|\bu^k(\bx_i) - \mathcal{I}_h(\bu^k)(\bx_i)|^2) \right] d\bxi.
		\end{split}
	\end{equation}
	Using the pointwise bound on the interpolant error, see~\eqref{eq:interpolationerrorpointwise}, for $\bu^k \in C^2(D)$, it follows
	\begin{equation}\label{eq:est1}
		I_2^2 \leq c_2 2\left(\frac{2C_2}{\epsilon^2 \mathrm{w}_d}\right)^2 \frac{\bar{J}_1}{\mathrm{w}_d} \int_{H_1(\bzero)} \frac{J(|\bxi|)}{|\bxi|} \left[ N c^2_4 h^4 \left(\sup_{\bx \in D} \left\vert\nabla^2 \bu^k(\bx)\right\vert \right)^2 \right] d\bxi,
	\end{equation}
	where recall that $N$ is the number of mesh nodes. Consider $\bar{\bv}_h = (\bar{v}_{1,h}, \bar{v}_{2,h}, ..., \bar{v}_{d,h}) \in \bar{V}_h$, $d = 2, 3$ being the spatial dimension, such that $\bar{v}_{n,h} = 0$ for $n\geq 2$ and $\bar{v}_{1,h} = 1$. Then, from~\eqref{eq:l2L2relation}, it holds that
	\begin{equation}
		c_1 N \leq ||\bar{\bv}_h||^2 = |D| \qquad \Rightarrow \qquad  N \leq \frac{|D|}{c_1}.
	\end{equation}
	Using above in~\eqref{eq:est1}, it follows that
	\begin{equation}
		\begin{split}
			I_2^2 &\leq 2c_2 \frac{c_4^2 h^4\,|D|}{c_1} \left(\sup_{\by \in D} |\nabla^2 \bu^k(\by)|\right)^2 \left(\frac{2C_2}{\epsilon^2 \mathrm{w}_d}\right)^2 \frac{\bar{J}_1}{\mathrm{w}_d} \int_{H_1(\bzero)} \frac{J(|\bxi|)}{|\bxi|}  d\bxi \\
			&\leq \left[ \frac{2h^2 c_4 \, \sqrt{c_2|D|/c_1}\, 2C_2\bar{J}_1 }{\epsilon^2} \sup_{\bx \in D} \left\vert\nabla^2 \bu^k(\bx)\right\vert \right]^2 \\
			&\leq \left[ \frac{h^2L_1 c_4 \, \sqrt{c_2|D|/c_1}}{\epsilon^2} \sup_{\bx \in D} \left\vert\nabla^2 \bu^k(\bx)\right\vert \right]^2,
		\end{split}
	\end{equation}
	where the definition of $L_1$ is used in the last step. This completes the proof of lemma. \qed
	
	\subsection{A-priori convergence}
	Let the discretization error $E^k$ at the $k^{\text{th}}$ step be given by
	\begin{equation}
		E^k = ||\bu_h^k - \bu(t^k)|| + ||\bv_h^k - \bv(t^k)||.
	\end{equation}
	Then, the application of triangle inequality and \eqref{eq:projerror} gives
	\begin{equation}\label{eq:totalerror}
		\begin{split}
			E^k &\leq ||\bu_h^k - \br_h(\bu(t^k))|| + ||\bv_h^k - \br_h(\bv(t^k))|| + ||\br_h(\bu^k) - \bu(t^k)|| + ||\br_h(\bv^k) - \bv(t^k)|| \\
			&= ||\be^k_{u,h}|| + ||\be^k_{v,h}|| + C_p h^2,
		\end{split}
	\end{equation}
	where
	\begin{equation}\label{eq:Cp}
		C_p = c_3\left[ \sup_t ||\bu(t)||_2 + \sup_t ||\partial_t{\bu}(t)||_2\right].
	\end{equation}
	
	The main result is as follows.
	\begin{theorem}\label{thm:convergence}
		\textbf{A-priori convergence of NFEA} \\
		If the solution $(\bu, \bv = \partial_t{\bu})$ of the  peridynamics equation~\eqref{eq:equationofmotion} is such that $\bu, \bv \in C^2([0,t_F]; H^2(D)\cap H^1_0(D) \cap C^2(D))$ then the scheme is consistent and the total error $E^k$ satisfies the following bound
		\begin{equation}\label{eq:totalErrorBound}
			\begin{split}
				&\sup_{k\leq t_F/\Delta t} E^k \\
				&\leq C_p h^2 + \exp[t_F\frac{(1+L_1\bar{n} \sqrt{c_2/c_1}/\epsilon^2)}{1-\Delta t}] \left[||\be^0_{u,h}|| + ||\be^0_{v,h}|| + \left(\frac{t_F}{1-\Delta t}\right) \left( C_t \Delta t + C_s \frac{h^2}{\epsilon^2} \right) \right],
			\end{split}
		\end{equation}
		where, constants $C_p$ and $C_t$ are defined in \eqref{eq:Cp} and \eqref{eq:consistenctytime}, respectively, and the constant $C_s$ is given by
		\begin{equation}
			\begin{split}
				C_s = &\left[\frac{c_3}{\epsilon^{5/2}} \left(L_2 \sup_t||\bu(t)||_2 + L_3(\sup_t ||\bu (t)||_2)^2\right)\right] \\
				&+ \left[\frac{L_1 c_4 \sqrt{\frac{c_2|D|}{c_1}}}{\epsilon^2} \sup_t \sup_{\bx\in D} \left\vert\nabla^2 \bu(\bx,t)\right\vert \right] + \left[\frac{2 c_3 L_1 \bar{n} \sqrt{\frac{c_2}{c_1}}}{\epsilon^2} \sup_t ||\bu(t)||_2  \right].
			\end{split}
		\end{equation}
	\end{theorem}
	
	The proof is similar to the proof of Theorem 5.1 in~\cite{jha2021finite} and relies on the estimates shown in~\cref{ss:keyest}.
	
	
	\section{An alternate nodal formulation based on the Cl\'ement interpolation}
	\label{sec:Clement Node}
	{An improved a priori convergence result is observed when the  Cl\'ement interpolant (see \cref{s:fespace}) is used in a nodal finite element formulation. More specifically, the Cl\'ement interpolation is used for the peridynamic force. 
		
		{Let $\bF_h^C(\bu^k_h)$ be the Cl\'ement interpolation of $\bF(\bu^k_h)$, \emph{i.e.},
			\begin{equation}
				\begin{split}
					\bF_h^C(\bu^k_h)(\bx)\bigg\vert_{T} &= \sum_{i\in N_T} \mathcal{P}_i(\bF(\bu^k_h))(\bx_i) \phi_i(\bx), \qquad \bx\in T, \quad \forall T\in \mathcal{T}_h, \label{cinterpF}
				\end{split}
			\end{equation}
			where, recall from \cref{s:fespace} that $P_i : L^2(S_i)^2 \to \mathcal{P}_1(S_i)$ is the $L^2$ projection, $S_i$ being the list of elements with $i$ as the vertex and $\mathcal{P}_1(S_i)$ the space of continuous piecewise linear finite elements on $S_i$.
			Then, \eqref{eq:centralstrong2.1} is written as
			\begin{equation}\label{eq:centralstrong99}
				\rho\dfrac{\bu^{k+1}_h - 2 \bu^k_h + \bu^{k-1}_h}{\Delta t^2} = \bF_h^C(\bu^k_h) + \bb^k_h.
			\end{equation}
			
			Previous analysis can be used to obtain the a-priori error estimates except for the peridynamics force. For the control of the error in peridynamics force, using the splitting of the error in \eqref{eq:consistency force nodal}, next results similar to \cref{lemma:nodal consistency} are obtained bounding terms $I_k$, $k=1,2,3,4$.

			First, an upper bound on $I_2$ is obtained. 
			From \eqref{eq:interpolationerror}, \eqref{eq:lipschitzL2}, \eqref{eq:Clemtnt stab}, it follows that
			\begin{equation}
				\label{I_2 second}
				\begin{split}
					I_2 = &||\bF_{h}^C(\bu^k) - \bF_{h}^C(\mathcal{I}_h(\bu^k))||
					= ||\mathcal{I}_h^C[\bF(\bu^k) - \bF(\mathcal{I}_h(\bu^k))]|| \\
					\leq &C||\bF(\bu^k) - \bF(\mathcal{I}_h(\bu^k))||\leq \frac{CL_1}{\epsilon^2} ||\bu^k - \mathcal{I}_h(\bu^k)|| \\
					\leq & \frac{c_3CL_1 h^2}{\epsilon^2}||\bu^k||_2.
				\end{split}
			\end{equation}
			To bound $I_3$, combining \eqref{eq:interpolationerror}, \eqref{eq:lipschitzL2}, \eqref{eq:Clemtnt stab} to get 
			\begin{equation}
				\label{I_3 second}
				\begin{split}
					I_3 = &||\bF_{h}^C(\mathcal{I}_h(\bu^k)) - \bF_{h}^C(r_h(\bu^k))||
					= ||\mathcal{I}_h^C[\bF(\mathcal{I}_h(\bu^k)) - \bF(r_h(\bu^k))]|| \\
					\leq &C||\bF(\mathcal{I}_h(\bu^k)) - \bF(r_h(\bu^k))||\leq \frac{CL_1}{\epsilon^2} ||\mathcal{I}_h(\bu^k) - r_h(\bu^k)|| \\
					\leq &\frac{CL_1}{\epsilon^2} \{||\mathcal{I}_h(\bu^k)-\bu^k||+||\bu^k - r_h(\bu^k)|| \}\leq2\frac{c_3CL_1 h^2}{\epsilon^2}||\bu^k||_2.
				\end{split}
			\end{equation}
			Applying $L^2$ stability of the Cl\'ement interpolant and arguing as in \eqref{eq: previous timesteperror} deliver
			\begin{equation}
				\label{I_4 second}
				\begin{split}
					I_4 = &||\bF_{h}^C(r_h(\bu^k)) - \bF_{h}^C(\bu_h^k)|| \leq \frac{\bar{n} \sqrt{c_2/c_1}\, L_1}{\epsilon^2} ||\be^k_{u,h}||.
				\end{split}
			\end{equation}
			Lastly, term $I_1$ is bounded by applying \eqref{eq:Clemtnt stab} with $t=0$ and $s=2$ and \eqref{eq:bddrnpforceH2} as follows 
			\begin{equation}
				\label{I_1 second}
				\begin{split}
					I_1 = &||\bF(\bu^k) - \bF_{h}^C(\bu^k)|| \leq Ch^2|| \bF(\bu^k)||_2 \leq \left[c_3C\frac{L_2 ||\bu^k||_2 + L_3 ||\bu^k||_2^2}{\epsilon^{5/2}} \right] h^2 \,.
				\end{split}
			\end{equation}
			
			Collecting the estimates \eqref{I_2 second}, \eqref{I_3 second}, \eqref{I_4 second}, and \eqref{I_1 second} and arguing again as in the proof of Theorem 5.1 in~\cite{jha2021finite} gives
			\begin{theorem}
				\textbf{A-priori convergence of NFEA with Cl\'ement interpolation of peridynamics force} \\
				\label{Clem}
				If the solution $(\bu, \bv = \partial_t{\bu})$ of the  peridynamics equation~\eqref{eq:equationofmotion} is such that $\bu, \bv \in C^2([0,t_F]; H^2(D)\cap H^1_0(D))$ then the scheme is consistent and the total error $E^k$ satisfies the following bound
				\begin{equation}
					\begin{split}
						&\sup_{k\leq t_F/\Delta t} E^k \\
						&\leq C_p h^2 + \exp[t_F\frac{(1+L_1\bar{n} \sqrt{c_2/c_1}/\epsilon^2)}{1-\Delta t}] \left[||\be^0_{u,h}|| + ||\be^0_{v,h}|| + \left(\frac{t_F}{1-\Delta t}\right) \left( C_t \Delta t + C_s \frac{h^2}{\epsilon^2} \right) \right],
					\end{split}
				\end{equation}
				where, constants $C_p$ and $C_t$ are defined in \eqref{eq:Cp} and \eqref{eq:consistenctytime}, respectively, and the constant $C_s$ is given by
				\begin{equation}
					\begin{split}
						C_s = &\left[\frac{c_3}{\epsilon^{5/2}} \left(L_2 \sup_t||\bu(t)||_2 + L_3(\sup_t ||\bu (t)||_2)^2\right)\right] \\
						&+ \left[\frac{3c_3CL_1 }{\epsilon^2}\sup_t{||\bu(t)||_2}\right]].
					\end{split}
				\end{equation}
			\end{theorem}
			
			This section is concluded with the following key remark.
			\begin{rem}
				The a-priori convergence rate for the alternate NFEA is an improvement over the rate given in Section~\cref{thm:convergence} as the constant $C_s$ depends only on $\sup_t||\bu(t)||_2$. 
				However, it is more expensive to implement. Future work will investigate the efficiency of the alternate NFEA scheme. Additionally, the Cl\'ement interpolation will be used to design a-posteriori estimates for use in adaptive schemes for mesh refinement.
			\end{rem}

			\section{Asymptotic compatibility of Cl\'ement NFEA}
			\label{sec:Asymp}
			In this section, the numerical error in initial values is assumed to be zero. 
			The asymptotic compatibility is established using the $C([0,T],L^2(D))$ compactness of sequences of solutions to the RNP model associated with vanishing horizon together with the fact that the limit displacement $\bu^0$ lies inside a dense subspace of $C([0,T],L^2(D))$.
			
			To begin, write the peridynamic solution for the RNP model given in~\cref{s:model} associated with horizon size $\epsilon$ as $\bu^{\epsilon}(t)$. Motivated by simulations, see~\cref{fig:compareHorizons} where the $L^\infty$ norms of displacement are bounded for three different horizons, 
			peridynamics solution is assumed to satisfy $\sup_{\epsilon>0}||{\bu}^\epsilon||_\infty<\infty$ for all $\epsilon$. With this hypothesis, multiplying~\cref{eq:equationofmotion} by $\dot{\bu^\epsilon}$, integrating by parts, and applying Gr\"onwall's inequality to find as in~\cite{CMPer-Lipton} that there is a sub-sequence denoted by $\bu^{\epsilon_n}$ converging strongly to a limit function $\bu^0$ in $C([0,t_F];L^2(D))$, i.e.,
			\begin{equation}\label{strongL2}
				\lim_{n\rightarrow\infty}\max_{t\in[0,t_F]}\Vert\bu^{\epsilon_n}(t)-\bu^0(t) \Vert_{L^2(D)}=0,
			\end{equation}
			where, $\bu^0$ belongs to SBD for every $t\in[0,t_F]$. Furthermore, there exists a constant $C$ depending only on $t_F$ bounding the Griffith (LEFM) energy, 
			\begin{eqnarray}
				\int_{D}\,\mu |\mathcal{E} \bu^0(t)|^2+\frac{\lambda}{2} |{\rm div}\,\bu^0(t)|^2\,dx+\mathcal{G}_c\mathcal{H}^{d-1}(J_{\bu^0(t)})\leq C, \hbox{ $d=2,3$,}
				\label{LEFMbound}
			\end{eqnarray}
			for $0\leq t\leq t_F$, where, $J_{u^0(t)}$ denotes the evolving fracture surface and $\mathcal{H}^{d-1}(J_{\bu^0(t)})$ is its $d-1$ dimensional Hausdorff measure at time $t$. Here the shear moduli $\mu$, and L\'ame moduli $\lambda$ are given by explicit formulas expressed in terms of $\psi'$ (derivative of the function $\psi$ in the RNP force model), see~\cite{CMPer-Lipton}. An analogous observation for H\"older continuous solutions is made in~\cite{jha2018numerical}.
			
			Defining as before the discrete times $t^{k_n}=k\Delta t_n$, $k=1,\ldots,N$, where $N=t_F/\Delta t_n$, let the approximation of $\bu^\epsilon(t^{k_n})$ based on the Cl\'ement NFEA is denoted by $({\mathbf{u}}^{\epsilon})^{k_n}_{h_n}$. Using the piecewise constant interpolation, the discrete solutions in times are extended to be defined at all times as follows:
			\begin{equation}\label{piece wise time}
				({\mathbf{u}}^{\epsilon})^{k_n}_{h_n}(t):=({\mathbf{u}}^{\epsilon})^{k_n}_{h_n} \hbox{ for $t$ in the interval $((k-1)\Delta t_n,k\Delta t_n)$}\,.
			\end{equation}
			
			\begin{theorem}
				\textbf{Asymptotic compatibility of NFEA with Cl\'ement interpolation of peridynamics force} \\
				\label{Asymtote}
				Let $({\mathbf{u}}^{\epsilon_n})^{k_n}_{h_n}$ be the Cl\'ement Nodal finite element approximation at time $t^{k_n}$ to $\mathbf{u}^{\epsilon_n}$, then there are sequences $\Delta t_n$ $h_n,\epsilon_n$, with $\Delta t_n\rightarrow 0$, $h_n\rightarrow 0$ and $\epsilon_n\rightarrow 0$ for $n=1,2,\ldots$ such that 
				$$\lim_{n\rightarrow\infty}\max_{t\in[0,t_f]}||({\mathbf{u}}^{\epsilon_n})^{k_n}_{h_n}(t)-{\mathbf{u}}^0(t)||_{L^2(D)}=0\,.$$
			\end{theorem}
			\begin{rem}\label{nouniform}
				It is noted that the method is not uniformly convergent with horizon $\epsilon_n$ as the constants multiplying the ratio $\frac{h_n^2}{\epsilon_n^2}$ go to $\infty$ as $\epsilon_n\rightarrow 0$. 
			\end{rem}
			
			To establish~\cref{Asymtote}, define
			\begin{equation}\label{piece wise time u}
				\tilde{\mathbf{u}}^{\epsilon_n,k_n}(t):={\mathbf{u}}^{\epsilon_n}(t^{k_n}) \hbox{ for $t$ in the interval $((k-1)\Delta t_n,k\Delta t_n)$},
			\end{equation}
			and apply the triangle inequality 
			\begin{align}\label{3ways}
				\max_{t\in[0,t_f]}||({\mathbf{u}}^{\epsilon_n})^{k_n}_{h_n}(t)-{\mathbf{u}}^0(t)||_{L^2(D)}\leq & \max_{t\in[0,t_f]}||({\mathbf{u}}^{\epsilon_n})^{k_n}_{h_n}(t)-\tilde{\mathbf{u}}^{\epsilon_n,k_n}(t)||_{L^2(D)}+\max_{t\in[0,t_f]}||\tilde{\mathbf{u}}^{\epsilon_n,k_n}(t)-\bu^{\epsilon_n}(t)||_{L^2(D)}\nonumber\\
				&+\max_{t\in[0,t_f]}||\bu^{\epsilon_n}(t)-\bu^0(t)||_{L^2(D)}:=A_1+A_2+A_3\,.
			\end{align}
			From~\cref{Clem}
			\begin{align}\label{Ione}
				A_1\leq\max_{t^{k_n}}\Vert{(\mathbf{u}}^{\epsilon})^{k_n}_{h_n}-\bu^{\epsilon_n}(t^{k_n})\Vert_{L^2(D)}\leq\frac{C^{\epsilon_n}}{\epsilon_n^2}(C_t\Delta t_n+C_s h_n^2)\,,
			\end{align}
			where, $C^{\epsilon_n}\rightarrow\infty$ as $\epsilon_n\rightarrow 0$. Set $O_{k_n}=((k-1)\Delta t_n,k\Delta t_n)$, and from Theorem 2.3 of \cite{CMPer-Lipton} one has
			\begin{align}\label{Itwo}
				A_2\leq\max_{k_n}\max_{t\in O_{k_n}}\Vert\bu^{\epsilon_n}(t^{k_n})-\bu^{\epsilon_n}(t)\Vert_{L^2(D)}\leq C\Delta t_n\,,
			\end{align}
			where, the constant $C$ is independent of $\epsilon_n$, $\bu^{\epsilon_n}$ and $\Delta t_n$. From the strong convergence~\cref{strongL2}, it is clear $A_3$ converges to $0$ for $\epsilon_n\rightarrow 0$. Thus, given a tolerance $\tau>0$, $\epsilon_n$ can be selected to be sufficiently small so that $A_3<\tau/3$. From-\cref{Itwo}, the size of time steps $\Delta t_n$ can be chosen sufficiently small so that $A_2\leq\tau/3$. Lastly, given $\epsilon_n$ from~\cref{Ione}, mesh size $h_n$ and $\Delta t_n$ can be selected sufficiently small to have $A_1\leq \tau/3$, and combining $A_1+A_2+A_3<\tau$ proving the theorem.
			
			\section{Numerical results}\label{s:results}
			This section presents results involving fracture evolution under different loading conditions and geometries. First, the procedure to numerically compute the peridynamics force is detailed, and the implementation of the NFEA method is briefly presented. In the implementation, integration over a horizon in the nodal peridynamics force is approximated by discrete summation involving node-node interaction; this is similar to the commonly used meshfree method and allows the NFEA method to be computationally efficient. Next, the material properties for numerical examples and calibration of the parameters in the peridynamics constitutive law are detailed. The remaining subsections are devoted to the numerical results. The first example concerns a simple non-fracture problem involving a square domain subjected to displacement-controlled loading, and the mesh convergence of the NFEA method is analyzed. The results show convergence rates depend on the two meshes used in the rate calculation, and the rates are close to or above 1.75. 
			In the second example, a Mode-I crack propagation problem is taken up. Using this example, several results are obtained: first, the convergence rate with mesh refinement; second, the solution is compared with the solution from the meshfree discretization based on \cite{silling2005meshfree, lipton2019complex, jha2019numerical}; and, third, the localization of damage zone is shown with refinement of horizon.  
			The third example involves a square specimen with a circular hole under displacement-controlled axial loading. This example shows the nucleation of the crack from the two points in the boundary of a hole. The fourth example is about the bending loading of the V-notch structure. This example also shows the crack nucleation. The last problem motivated by \cite[Figure 18]{CMPer-DaiBFEM} includes a rectangle specimen with a hole and pre-crack. This example shows the effect of stress concentration near the hole on crack path and propagation. 
			The section ends with a discussion of the crack speeds for the four problems. 
			
			Numerical results were obtained using a code similar to the C\texttt{++} code NLMech~\cite{diehl2020asynchronous,jha2021nlmech}. In all results, the mesh consisted of linear triangle elements. The second-order quadrature scheme calculates the integration over a finite element (triangle elements); see next subsection. The strain field -- symmetric gradient of the displacement -- is constant over each element and computed at the element's center. For the triangulation of a domain with a void and notch, an open-source library Gmsh~\cite{Gmsh} is utilized, and Paraview~\cite{paraview} is used to visualize the results. 
			
			\begin{figure}
				\centering
				\includegraphics[scale=0.18]{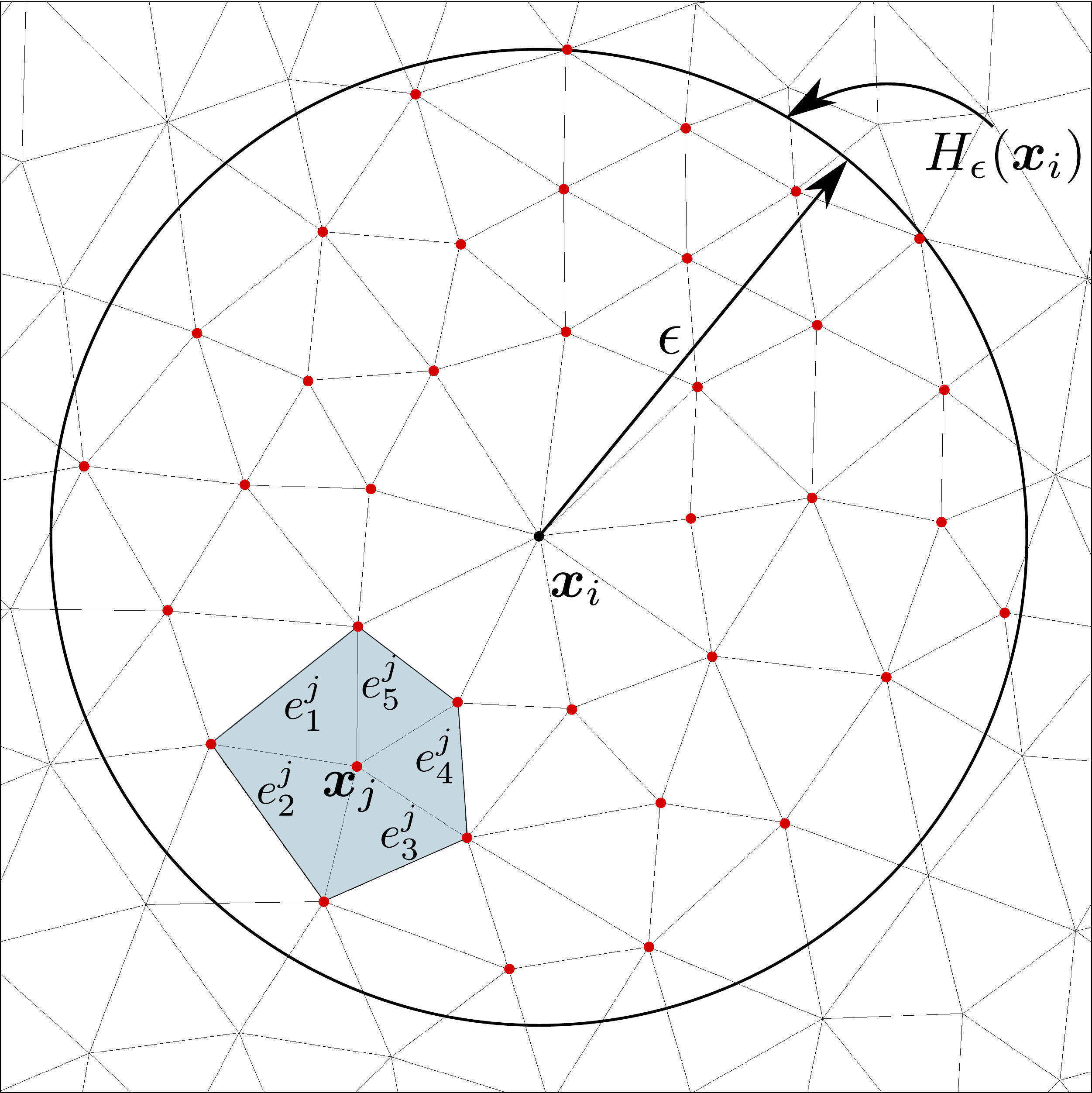}
				\caption{Typical mesh node $\bx_i$ and one of the neighboring nodes $\bx_j$ in an example 2-d finite element mesh. All the red nodes contribute to the force at $\bx_i$. The set $E_j = \{e^j_k\}_{k=1}^5$ of elements with the node $\bx_j$ as the vertex is shown in grey.}
				\label{fig:pdforcecalc}
			\end{figure}
			
			\subsection{Computation of a peridynamics force in NFEA}\label{ss:strongfeaimplementation}
			This section develops a computationally efficient approximation of the peridynamics forces at nodes based on the discrete summation of the node-node interaction. As a result of the approximation, the NFEA and meshfree methods will have the same computational complexity. The downside of the approximation, however, is the loss of accuracy. In the a-priori error analysis, the discretization error of nodal peridynamics forces is not considered.
			
			Let $\bu^k_h, \bv^k_h \in V_h$ be the finite element displacement and velocity functions, $V_h$ being the finite element space (see \eqref{eq:dispstrong} and \cref{ss:nodalfe}). Corresponding to $\bu^k_h$ and $\bv^K_h$, suppose $\bolds{U}^k, \bolds{V}^k$ are nodal displacement and velocity vectors, respectively, i.e., $\bolds{U}^k = (\bolds{U}^k_1, ..., \bolds{U}^k_N), \bolds{V}^k = (\bolds{V}^k_1, ..., \bolds{V}^k_N)$. Velocity is given by 
			\begin{equation*}
				\bolds{V}^k_i = (\bolds{U}^k_i - \bolds{U}^{k-1}_i)/\Delta t    
			\end{equation*}
			when $k\geq 1$ and $\bolds{V}^k_i = \bv_0(\bx_i)$ when $k=0$, where $\bv_0$ is the prescribed initial condition for the velocity. From \eqref{eq:centralstrong}, $\bolds{U}^k$ is computed using, for $k\geq 1$ and all $i$,
			\begin{equation}
				\bolds{U}^{k+1}_i = \Delta t^2 \frac{\bF(\bu^k_h)(\bx_i) + \bb(\bx_i, t^k)}{\rho} + 2 \bolds{U}^k_i - \bolds{U}^{k-1}_i
			\end{equation}
			and, for $k=0$ and all $i$,
			\begin{equation}
				\bolds{U}^{1}_i = \frac{\Delta t^2}{2} \frac{\bF(\bu^0_h)(\bx_i) + \bb(\bx_i, 0)}{\rho} + \Delta t \bv_0(\bx_i) + \bu_0(\bx_i) .
			\end{equation}
			In the above, the numerical evaluation of peridynamics force $\bF(\bu^k_h)(\bx_i)$ is nontrivial and, therefore, is detailed next. 
			
			\begin{algorithm}[H]
				\caption{NFEA implementation}
				\label{alg:centraldifference}
				\begin{algorithmic}[1]
					\STATE Read nodes and element-node connectivity from the mesh file
					\STATE \textcolor{mygray}{\it $\%\%$ Task: Create neighbor list and compute $V_{ij}$ using \eqref{eq:computeVij}}
					\FOR {Each integer $0\leq i \leq N-1$} \textcolor{mygray}{\it $\%$ $N$ is the total number of nodes}
					\FOR {Element $e \in E_j$} \textcolor{mygray}{\it $\%$ $E_j$ list of elements with node $j$ as vertex}
					\FOR {$1\leq q \leq Q^e$} \textcolor{mygray}{\it $\%$ Loop over quadrature points of element $T_e$}
					\IF {$\vert \bx_q - \bX[i] \vert \leq \epsilon$}
					\STATE Add $j$ to neighborList$[i]$
					\STATE Compute $V_{ij}$ using \eqref{eq:Vij}, add $V_{ij}$ to V$[i]$
					\ENDIF
					\ENDFOR
					\ENDFOR
					\ENDFOR \textcolor{mygray}{\it $\:\%$ End of loop over nodes for neighborlist}
					\STATE \textcolor{mygray}{\it $\%\%$ Task: Integrate in time}
					\FOR {Each integer $0\leq k \leq t_F/\Delta t$}
					\STATE \textcolor{mygray}{\it $\%$ $\bolds{U}, \bolds{V}$ are the displacement and velocity at step $k$}
					\STATE \textcolor{mygray}{\it $\%\%$ Task: Compute force $\bF$ using \eqref{eq:approxPDfinal}}
					\STATE Initialize vector $\bolds{F}$ with zeros
					\FOR {Each integer $0\leq i \leq N-1$}
					\FOR {Each integer $j\in$ neighborList$[i]$} 
					\STATE $S_-{ji} = \frac{\bolds{U}[j] - \bolds{U}[i]}{|\bolds{X}[j] - \bolds{X}[i]|} \cdot \frac{\bolds{X}[j] - \bolds{X}[i]}{|\bolds{X}[j] - \bolds{X}[i]|}$
					\STATE $\bolds{F}[i] = \bolds{F}[i] + \frac{\omega(\bolds{X}[i]) \omega(\bolds{X}[j])}{\mathrm{w}_d \epsilon^{d+1}}\, \psi'(|\bolds{X}[j] - \bolds{X}[i]|\,S_-{ji}^2) \, S_-{ji} \, \frac{\bolds{X}[j] - \bolds{X}[i]}{|\bolds{X}[j] - \bolds{X}[i]|} \,V[i][j]$
					\ENDFOR
					\ENDFOR	\textcolor{mygray}{\it $\:\%$ End of loop over nodes for $\bolds{F}$}
					\STATE \textcolor{mygray}{\it $\%\%$ Task: Update displacement $\bolds{U}$ and velocity $\bolds{V}$}						
					\FOR {Each integer $0\leq i \leq N-1$}
					\FOR {Each integer $0\leq l \leq d-i$} \textcolor{mygray}{\it $\%$ $d$ is the dimension of the problem}
					\STATE $\bolds{U}_-{temp} = \bolds{U}[i][l]$
					\IF {dof $l$ of node $i$ is free}
					\STATE $\bolds{U}[i][l] = \bolds{U}[i][l] + \Delta t \, \bolds{V}[i][l] + \Delta t^2 \frac{\bolds{F}[i][l] + \bb_-{l}(\bolds{X}[i], k\Delta t)}{\rho(\bolds{X}[i])}$						
					\ELSE
					\STATE Read $\bolds{U}[i][l]$ from boundary condition
					\ENDIF
					\STATE $\bolds{V}[i][l] = \frac{\bolds{U}[i][l] - \bolds{U}_-{temp}[i][l]}{\Delta t}$
					\ENDFOR
					\ENDFOR \textcolor{mygray}{\it $\:\%$ End of loop over nodes for $\bolds{U}$ and $\bolds{V}$ update}
					\ENDFOR \textcolor{mygray}{\it $\:\%$ End of loop over time}
				\end{algorithmic}
			\end{algorithm}
			
			From~\eqref{eq:rnptotalforce}, it holds that
			\begin{equation*}
				\bF(\bu^k_h)(\bx_i) = \sum_{T\in \mathcal{T}_h} \int_{H_\epsilon(\bx_i)\cap T} \frac{\omega(\bx_i) \omega(\by)}{\mathrm{w}_d \epsilon^{d+1}} J^\epsilon(|\by - \bx_i|) \psi'(|\by - \bx_i|S(\by, \bx_i;\bu^k_h)^2) S(\by , \bx_i;\bu^k_h) \frac{\by - \bx_i}{|\by - \bx_i|} d\by.
			\end{equation*}
			Let $N_{T}$ be the list of nodes that are vertices of element $T$. Recall that $\phi_i$ denotes the interpolation function of node $i$. For $\by \in T$, $\bu^k_h(\by) = \sum_{j\in N_{T}} \phi_j(\by) \bolds{U}^k_j$. Also, for any node $i$, $\bu^k_h(\bx_i) = \bolds{U}^k_i = \sum_{j\in N_{T}} \phi_j(\by) \bolds{U}^k_i$ for all $\by \in T$ (due to the partition of unity property, i.e., $\sum_{j \in N_{T}} \phi_j(\by) = 1$). Combining, it follows that
			\begin{equation*}
				\begin{split}
					\bF(\bu^k_h)(\bx_i) = \sum_{T\in \mathcal{T}_h} \bigg[ \sum_{j \in N_T} \int_{H_\epsilon(\bx_i)\cap T} & \frac{\omega(\bx_i) \omega(\by)}{\mathrm{w}_d \epsilon^{d+1}} J^\epsilon(|\by - \bx_i|) \psi'(|\by - \bx_i|S(\by, \bx_i;\bu^k_h)^2) \\ 
					&  \left( \phi_j(\by)\frac{\bolds{U}^k_j - \bolds{U}^k_i}{|\by - \bx_i|} \cdot \frac{\by - \bx_i}{|\by - \bx_i|} \right) \frac{\by - \bx_i}{|\by - \bx_i|} d\by \bigg].
				\end{split}
			\end{equation*}
			Motivated from the above, peridynamics force $\bF(\bu^k_h)(\bx_i)$ can be approximated as follows
			\begin{equation}\label{eq:approxPD}
				\begin{split}
					\bF(\bu^k_h)(\bx_i) \approx \sum_{j \in \textrm{PD}_i} &\frac{\omega(\bx_i) \omega(\bx_j)}{\mathrm{w}_d \epsilon^{d+1}} \psi'(|\bx_j - \bx_i|S(\bx_j, \bx_i;\bu^k_h)^2) \left(\frac{\bolds{U}^k_j - \bolds{U}^k_i}{|\bx_j - \bx_i|} \cdot \frac{\bx_j - \bx_i}{|\bx_j - \bx_i|} \right) \frac{\bx_j - \bx_i}{|\bx_j - \bx_i|} \\
					&\left[ \sum_{e\in E_j} \int_{H_\epsilon(\bx_i)\cap T_e} J^\epsilon(|\by - \bx_i|) \phi_j(\by) d\by\right]\,,
				\end{split}
			\end{equation}
			where, $E_j$ is the list of elements with node $j$ as its vertex, see \cref{fig:pdforcecalc}, and $\textrm{PD}_i$ is the list of neighboring nodes $j$ and it is defined as
			\begin{equation}
				\textrm{PD}_i := \{1\leq j \leq N: T_e \cap H_\epsilon(\bx_i) \neq \emptyset \text{ for at least one } e \in E_j\}\,.
			\end{equation}
			Note that $\{1\leq j \leq N: \bx_j \in H_\epsilon(\bx_i)\} \subset \textrm{PD}_i$, i.e., $\textrm{PD}_i$ consists of nodes that are in the horizon $H_\epsilon(\bx_i)$ and additional nodes outside $H_\epsilon(\bx_i)$ which may belong to element $T_e$ that intersects $H_\epsilon(\bx_i)$. 
			The above form of approximation is not unique, as one may also approximate the force as
			\begin{equation*}
				\begin{split}
					\bF(\bu^k_h)(\bx_i) \approx \sum_{j \in \textrm{PD}_i} &\frac{\omega(\bx_i) \omega(\bx_j)}{\mathrm{w}_d \epsilon^{d+1}} \left(\frac{\bolds{U}^k_j - \bolds{U}^k_i}{|\bx_j - \bx_i|} \cdot \frac{\bx_j - \bx_i}{|\bx_j - \bx_i|} \right) \frac{\bx_j - \bx_i}{|\bx_j - \bx_i|} \\
					&\left[ \sum_{e\in E_j} \int_{H_\epsilon(\bx_i)\cap T_e} \psi'(|\by - \bx_i|S(\by, \bx_i;\bu^k_h)^2) J^\epsilon(|\by - \bx_i|) \phi_j(\by) d\by\right].
				\end{split}
			\end{equation*}
			Similarly, other forms of approximation are possible by keeping some terms outside and some inside of the integration. In our implementation, the approximation~\eqref{eq:approxPD} is used for two reasons: \textit{1.} The term in the square bracket is independent of time and, therefore, can be computed only once in the beginning and stored, and \textit{2.} The choice of keeping nonlinear term outside the integral as well as the vector $(\bx_j - \bx_i)/|\bx_j - \bx_i|$ gives a stable simulation, and numerical results agree well with the benchmark problems. Proceeding further, let $V_{ij}$ be the weighted volume of a node $j$ for a pairwise force contribution to the node $i$. It is defined as
			\begin{equation}\label{eq:Vij}
				V_{ij} = \sum_{e\in E_j} \int_{H_\epsilon(\bx_i)\cap T_e} J^\epsilon(|\by - \bx_i|) \phi_j(\by) d\by .
			\end{equation}
			The above integration over an element is computed using the quadrature rule. In all the numerical results, the second-order quadrature rule is employed; higher-order schemes can be used as the above integration needs to be computed only once and stored in the memory. Let $Q^e$ be the number of quadrature points associated with the element $e$ ($e$ is the element number, and $T_e$ is the element). Further, let, for $q=1, ..., Q^e$, $(\bx_q, w_q)$ is the pair of quadrature points and weights. Then
			\begin{equation}\label{eq:computeVij}
				V_{ij} = \sum_{e\in E_j} \left[ \sum_{q=1}^{Q^e} \chi_{H_\epsilon(\bx_i)}(\bx_q) J^\epsilon(|\bx_q - \bx_i|) \phi_j(\bx_q) w_q \right],
			\end{equation}
			where $\chi_A(\bx)$ is the indicator function taking value $1$ if $\bx \in A$ and $0$ if $\bx\notin A$. Using the definition of $V_{ij}$, \eqref{eq:approxPD} can be written as
			\begin{equation}\label{eq:approxPDfinal}
				\bF(\bu^k_h)(\bx_i) \approx \sum_{j \in \textrm{PD}_i} \frac{\omega(\bx_i) \omega(\bx_j)}{\mathrm{w}_d \epsilon^{d+1}} \psi'(|\bx_j - \bx_i|S(\bx_j, \bx_i;\bu^k_h)^2) \left(\frac{\bolds{U}^k_j - \bolds{U}^k_i}{|\bx_j - \bx_i|} \cdot \frac{\bx_j - \bx_i}{|\bx_j - \bx_i|} \right) \frac{\bx_j - \bx_i}{|\bx_j - \bx_i|} V_{ij}.
			\end{equation}
			In \cref{fig:pdforcecalc}, one of the neighboring nodes $\bx_j$ contributing to the force at $\bx_i$ is shown on an example 2-d finite element mesh. \cref{alg:centraldifference} presents the implementation of NFEA.
			
			\subsection{Material properties}\label{sec:material:properties}
			Let $\rho$ denote the density, $E$ Young's modulus, $\nu$ Poisson ratio, and $G_c$ critical energy release rate. The bond-based peridynamics suffer from the restriction of Poisson ratio $\nu = 1/4$ in 3-d or 2-d plane strain and $\nu = 1/3$ in 2-d plane stress; see \cite{trageser2020bond}. All of the simulations are in 2-d, and plane strain is assumed. Therefore, $\nu$ is fixed to $1/4$. 
			
			To fix the parameters in the RNP model, see \eqref{eq:rnpbondforce}, the nonlinear potential function $\psi$ is set to $\psi(r) = c (1 - \exp[-\beta r])$, where $c$ and $\beta$ are two model parameters. The influence function is taken to be $J^\epsilon(r) = J(r/\epsilon)$, where $J(r) =  1 - r$ for $0\leq r< 1$ and $J(r) = 0$ for $r\geq 1$. The boundary function $\omega(\bx)$ is taken as $1$ for all points $\bx$ in the domain, i.e., $\omega(\bx) = 1$ for $\bx \in D$. Given $E$, Lam\'e parameter are $\lambda = \mu = 2E/5$ ($\nu = 1/4$ is assumed). The parameters $c$ and $\beta$, in 2-d, can be determined from (see \cite{CMPer-Lipton})
			\begin{align}
				c = \frac{G_c \pi}{4M_J}, \quad \beta = \frac{8E}{5cM_J},
			\end{align}
			where, $M_J = \int_0^1 J(r) r^2dr = 1/12$ for $J(r) = 1-r$. The inflection point of the potential function $\psi$ is given by $r^\ast = 1/\sqrt{2\beta}$ and the critical strain $S_c(\by, \bx) = \pm r^\ast/\sqrt{|\by - \bx|}$. 
			
			Let $c_L$, $c_S$, and $c_R$ are the longitudinal, shear, and Rayleigh wave speeds, respectively. Given elastic properties such as $E$ and $\nu$, wave speeds can be computed using the formulae:
			\begin{equation}\label{eq:wavespeeds}
				\begin{split}
					c_L &= \sqrt{\frac{\lambda + 2\mu}{\rho}} = \sqrt{\frac{1}{\rho}\frac{E(1-\nu)}{(1+
							\nu)(1-2\nu)}}, \quad c_S = \sqrt{\frac{\mu}{\rho}} = \sqrt{\frac{1}{\rho}\frac{E}{2(1+\nu)}}, \quad c_R \approx c_S \left(\frac{0.862 + 1.14\nu}{1+\nu}\right),
				\end{split}
			\end{equation}
			where the last formula to approximate Rayleigh wave speed can be found in \cite{royer2007improved}. Material properties employed in numerical experiments are listed in \cref{tab:matprops}. 
			
			\begin{table}
				\centering
				\addtolength{\tabcolsep}{+2pt}
				\renewcommand{\arraystretch}{1.3}
				
				\begin{tabular}{|l||r|c|l||r|}  
					\hline
					Properties & Values & & Properties & Values \\
					\hline\hline
					$\rho$ $(kg/m^3)$ & $1200$ && $c_L$ $(m/s)$ & $6123.7$ \\  
					\hline
					$E$  $(GPa)$ & $37.5$ && $c_S$ $(m/s)$ & $3535.5$ \\
					\hline
					$G_c$ $(J/m^2)$ & $500$ && $c_R$ $(m/s)$ & $3244.2$ \\
					\hline
				\end{tabular}
				
				\caption{Material properties. Formulae in \eqref{eq:wavespeeds} are used to compute longitudinal ($c_L$), shear ($c_S$), and Rayleigh ($c_R$) wave speeds.}
				\label{tab:matprops}
			\end{table}
			
			\begin{definition}[\textbf{Damage}]\label[definition]{def:damage}
				The damage at the material point $\bx$ is defined as
				\begin{align}\label{eq:defZ}
					Z(\bx) := \sup_{\by \in H_\epsilon(\bx) \cap D} \frac{|S(\by, \bx)|}{|S_c(\by, \bx)|}.
				\end{align}
				Based on the above, if $Z(\bx) \geq 1$, it follows that $\bx$ has at least one bond in the neighborhood with the bond strain above the critical bond strain. The damage zone of the material is given by the set $\{\bx\in D: Z(\bx) \geq 1\}$. Other measures of damage are also possible. For example, consider a function $\varphi$ given by
				\begin{equation}
					\varphi(\bx) = \frac{\int_{H_\epsilon(\bx)\cap D} \mu(\by, \bx; S) d\by}{\int_{H_\epsilon(\bx)\cap D} d\by},
				\end{equation}
				where, $\mu(\by, \bx; S)$ is a function that models the breakage of bond:
				\begin{equation*}
					\begin{split}
						\mu(\by, \bx; S) = \begin{cases}
							1, \qquad \text{if } S(\by, \bx) < S_c(\by,\bx), \\
							0, \qquad \text{otherwise}.
						\end{cases}
					\end{split}
				\end{equation*}
				Thus, $\varphi(\bx) \in [0, 1]$, and $\varphi(\bx) = 0$ implies all the bonds in the neighborhood of $\bx$ are stretched below the critical value, while $\varphi(\bx) = 1$ implies all the bonds in the neighborhood have strains above the critical value. 
			\end{definition}
			
			\begin{rem}\label{rem:plot}
				The plots of fields, such as damage, displacement, and linearized strain (symmetric gradient of displacement), are based on the finite element representation (surface plot in Paraview) and are in a deformed configuration unless otherwise stated in the figure's caption. The strain field is a piecewise constant over elements, and the strain value in each element is computed at the element's center. 
				
				The field $Z(\bx)$ is displayed using two colors (blue and red) to visualize damage. The blue color means $Z(\bx) \leq 1$ and the point $\bx$ is colored red if $Z(\bx) := \max_{\by\in H_\epsilon(\bx)} |S(\by, \bx)| / S_c(\by, \bx) > 1$, i.e., red colored point $\bx$ has at least one neighboring point $\by$ such that the bond $\bx$ -- $\by$ is stretched above the threshold stretch. Thus, red indicates that the point has at least one critically stretched bond.
			\end{rem}
			
			\subsection{Convergence test on square domain with displacement controlled loading}\label{ss:convgtest}
			A simple example of displacement-controlled loading is considered to test the convergence of the Nodal FE approximation. Consider a two-dimensional solid body $D = [0, 1 \text{m}]^2$ with density $\rho = 1$ kg/m$^3$, Young's modulus $E = 1$ Pa, and Poisson ratio $\nu = 0.25$. Thickness is one meter, and the horizon is fixed to $\epsilon = 0.05$ m. The left layer of specimens of $\epsilon$ thickness is clamped (zero displacements in x and y directions). The layer on the right of thickness $\epsilon$ is subjected to displacement $\bu_x(\bx, t) = 0.01 sin(2\pi t)$ in the x direction while the displacement in the y direction is kept free; see \cref{fig:convgtest}(a). To test the convergence, four different mesh sizes are considered: $h_i = \epsilon/m_i$ for $i=1,2,3,4$ with $m_1 = 4, m_2 = 8, m_3 = 12, m_4 = 16$. The domain is discretized into the uniform grid, and each grid is further divided into two triangles; see representative mesh in \cref{fig:convgtest}(b). The size of the grid is the mesh size. The final simulation time is $t_F = 0.5$ s, the time step size $\Delta t = 0.000125$ s, and the results are written to file every $\Delta t_{out} = 0.01$ s interval.

			\begin{figure}
				\centering
				\includegraphics[width=0.75\textwidth]{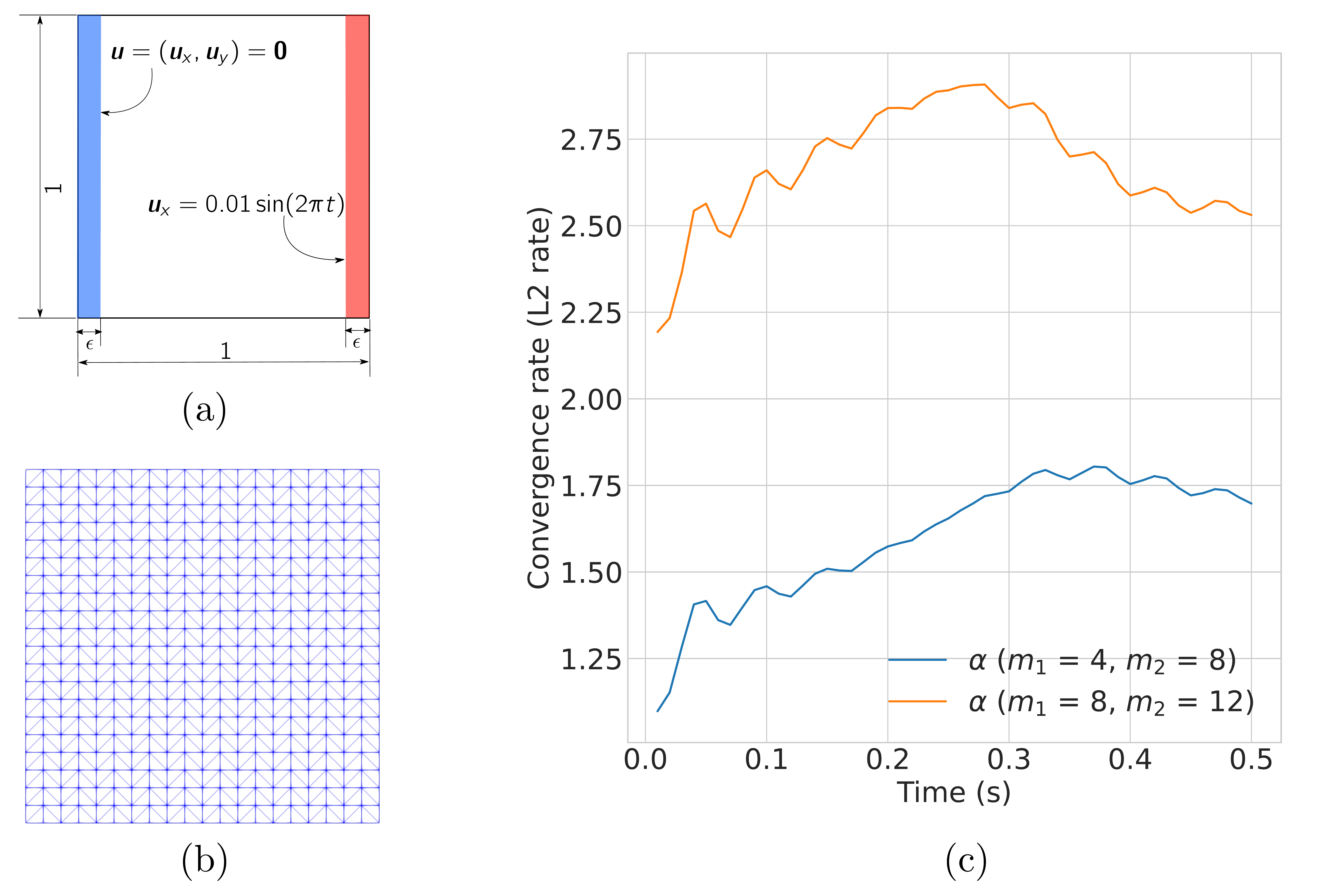}
				\caption{(a) {\bf Convergence test}: Setup. The horizon is fixed to $\epsilon = 0.05$ m. (b) A representative view of the mesh consisting of linear triangle elements. (c) Rate of convergence at discrete times using \eqref{eq:convergencerate}.}\label{fig:convgtest}
			\end{figure}
			
			\begin{figure}
				\centering
				\includegraphics[width=\textwidth]{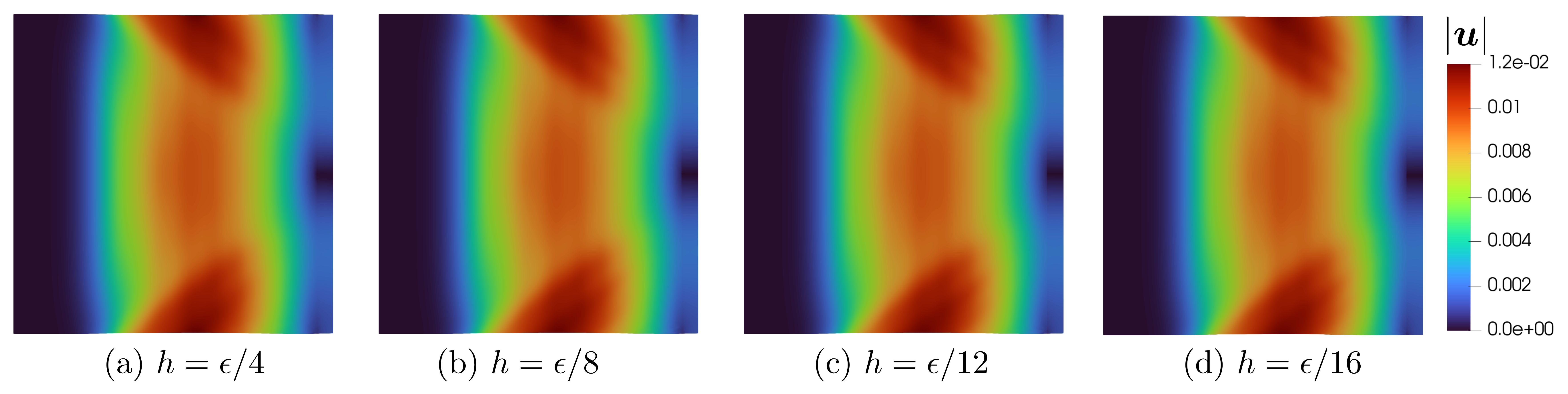}
				\caption{{\bf Convergence test}: Comparing the magnitude of displacement at the final simulation time $t_F$ for all four meshes. Unless otherwise stated, all the plots, including this figure, are in the deformed configuration and based on the finite element representation.}\label{fig:convgtestUMag}
			\end{figure}
			
			To calculate the convergence rate with mesh refinement, numerical solution $\bu_{h_4}(t_k)$ corresponding to the finest mesh is treated as an exact solution, where $t_k$, $k = 1, 2, ..., t_F/\Delta t_{out}$, is the output time (time at which solutions are written to a file). Setting $\bu = \bu_{h_4}$, the rate of convergence at a time $t_k$ from two solutions $\bu_{h_i}(t_k)$ and $\bu_{h_{i+1}}(t_k)$ can be computed as follows:
			\begin{equation}\label{eq:convergencerate}
				\alpha(t^k) = \frac{\log \left(||\bu_{h_i}(t_k) - \bu(t_k)||\right) - \log \left(||\bu_{h_{i+1}}(t_k) - \bu(t_k)||\right)}{\log(h_i) - \log(h_{i+1})}\,.
			\end{equation}
			Using the numerical solutions for three mesh sizes $h_1 = \epsilon/m_1, h_2 = \epsilon/m_2, h_3 = \epsilon/m_3$ ($h_4$ is used as an ``exact" solution), convergence rates $\alpha(t_k; m_1, m_2)$ and $\alpha(t_k; m_2, m_3)$ at output times $t_k$, $k=1,2, ..., t_F/\Delta t_{out}$, are computed and displayed in \cref{fig:convgtest}(c). Results show that the convergence rate is quite good and above 2 as the mesh size is reduced from $\epsilon/8$ to $\epsilon/12$. When reducing the mesh size from $\epsilon/4$ to $\epsilon/8$, the convergence rate is above $1.5$ and below the optimal rate of $2$ from the a-priori error analysis. These rates, however, are excellent considering that the NFEA is similar to the meshfree discretization and approximates the nodal peridynamics forces using a discrete summation of node-node interaction (integration through the quadrature method will increase the computational cost). In \cref{fig:convgtestUMag}, the magnitude of the displacement field is displayed for four simulations with different mesh sizes at the final time $t_F = 0.5$ s. Visually, all the results are in agreement.
			
			\subsection{Mode-I crack propagation}\label{ss:fractureproblem}
			Consider a square domain $D=[0,100\, \text{mm}]^2$ with a vertical pre-crack of length $l=20\,$mm located at the center; see \cref{fig:mode1crack}(a). If a specimen has a pre-crack line/curve (surface in 3-d), it means that in the peridynamics simulation, all the bonds intersecting the pre-crack line/curve (surface in 3-d) are initially broken and are kept broken during the course of the simulation. The constant velocity of $\pm 10^3$ mm/s is specified on the small area on the left and right sides to obtain the mode-I crack propagation. The simulation time and the size of the time step are $t_F = 40\, \mu$s and $\Delta t = 0.0008\,\mu$s, respectively, and the output is written to a file every $\Delta t_{out} = t_F / 50$ interval. The nonlocal length-scale, i.e., horizon, is fixed to $\epsilon = 2\,$mm. In the simulations, the RNP model with the material properties listed in \cref{tab:matprops} is employed.
			
			As in the previous example, to get the estimated convergence rate with mesh refinement, discretized solutions for mesh sizes $h_i = \epsilon / m_i$, $i=1,2,3,4$, with $m_1 = 4, m_2 = 8, m_3 = 12, m_4 = 16$ are computed. The finite element solution $\bu_{h_4}(t_k)$ corresponding to the finest mesh, i.e., $h = h_4 = \epsilon / 16$, is employed as an exact solution $\bu(t_k) = \bu_{h_4}(t_k)$, and using simulation results $\bu_{h_1}(t_k), \bu_{h_2}(t_k), \bu_{h_3}(t_k)$, convergence rates $\alpha(t_k; m_1, m_2)$ and $\alpha(t_k; m_2, m_3)$ at output times $t_k$, $k = 1, 2, ..., t_F/\Delta t_{out}$, are computed following \eqref{eq:convergencerate}. Note that $m_1$ and $m_2$ in $\alpha(\cdot; m_1, m_2)$ indicate that the rate is for meshes $h_1 = \epsilon/m_1$ and $h_2 = \epsilon/m_2$. 
			Two convergence rates, $\alpha(\cdot; m_1, m_2)$ and $\alpha(\cdot; m_2, m_3)$, are shown in \cref{fig:mode1crackrate}. 
			In \cref{fig:modeICompareMeshResults}, the fields $\bu_x$ and $Z$ are compared for four meshes at time $t = 28\,\mu$ s. The results show that the plots are visually indistinguishable.
			
			Next, the plot of the damage function $Z$ defined in \eqref{eq:defZ} is presented in the left column of \cref{fig:mode1plots}. In the right column, linearized strain is computed from the displacement field, and its magnitude (magnitude of the strain tensor $\bE = \frac{1}{2}\left[\nabla \bu + \nabla \bu^T\right]$ is taken as $||\bE|| = \sqrt{\bE\bcolon\bE}$, where $\bE\bcolon\bE = E_{ij}E_{ij}$ is the dot product) is depicted. In all the numerical results, it is found that the width of the process zone (damaged region) is approximately twice the horizon and envelopes the crack interface. Further, the strain tensor magnitude is unusually higher at the crack interface, as expected. To show crack opening, in \cref{fig:mode1plotsZcod}, the displacement is magnified by a factor of 100 and added to the reference configuration, and the damage $Z$ is displayed in the new artificial deformed configuration.
			
			\begin{figure}
				\centering
				
				\begin{subfigure}[t]{0.4\textwidth}
					\centering
					\includegraphics[width=0.9\textwidth]{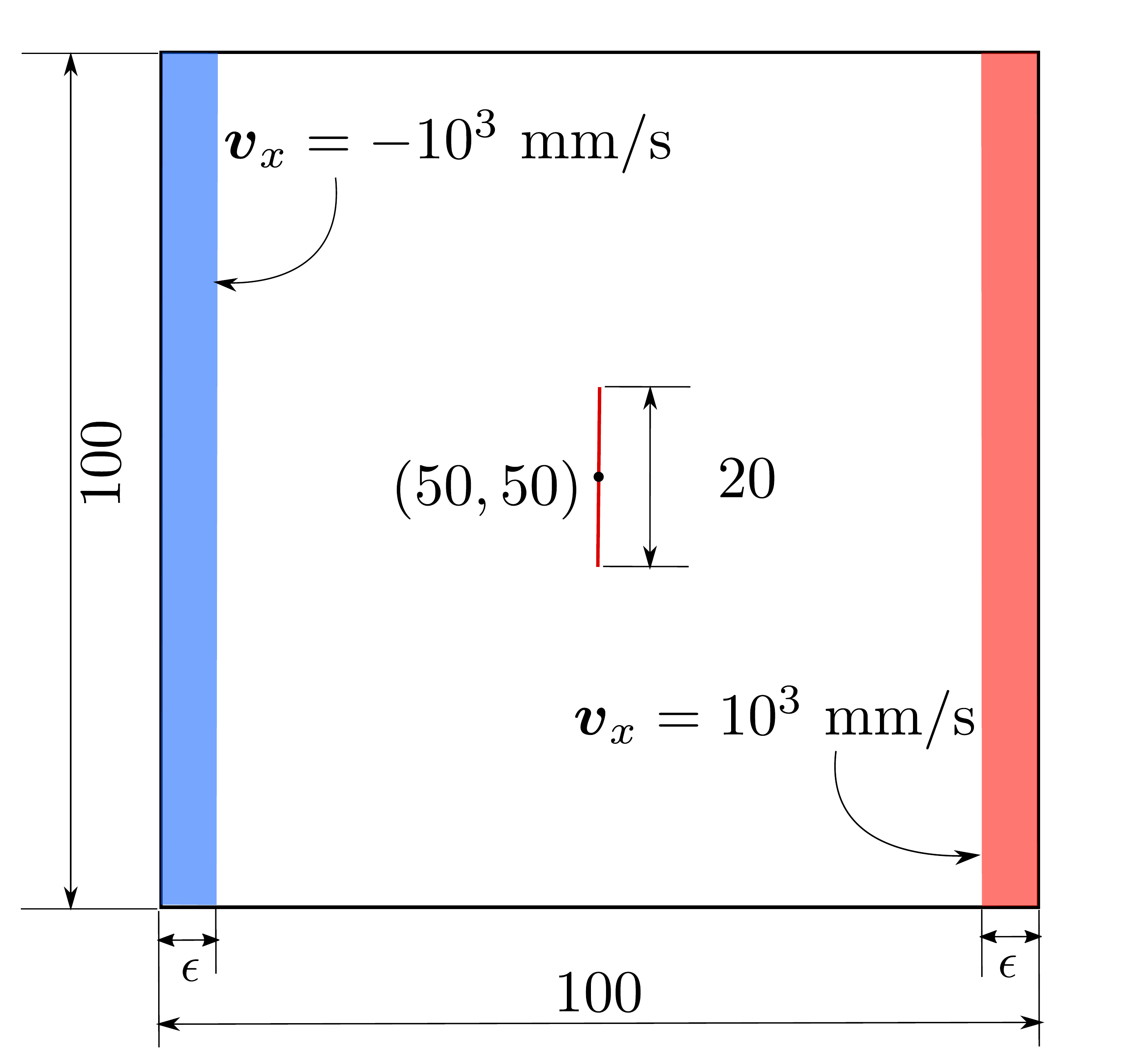}
					\caption{}
					\label{fig:mode1crack}
				\end{subfigure}  
				\begin{subfigure}[t]{0.55\textwidth}
					\centering
					\includegraphics[width=0.6\textwidth]{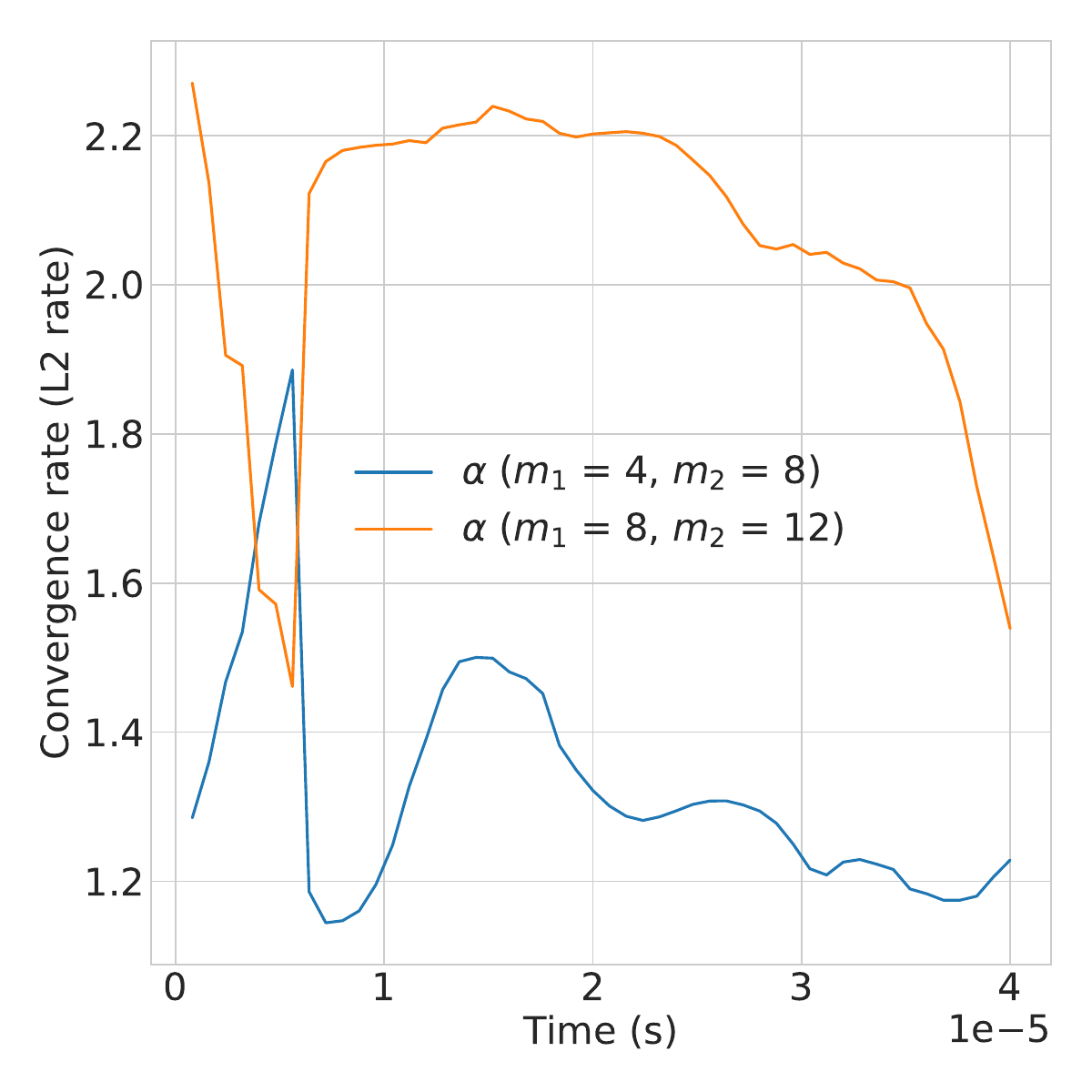}
					\caption{}
					\label{fig:mode1crackrate}
				\end{subfigure}
				\caption{(a) {\bf Mode-I crack problem}: Setup. The horizon is $\epsilon = 2$ mm. The center of the pre-crack coincides with the center of the domain. (b) {\bf Mode-I crack problem}: Convergence rate as the mesh is refined.}
			\end{figure}
			
			\begin{figure}
				\centering
				\includegraphics[width=0.9\textwidth]{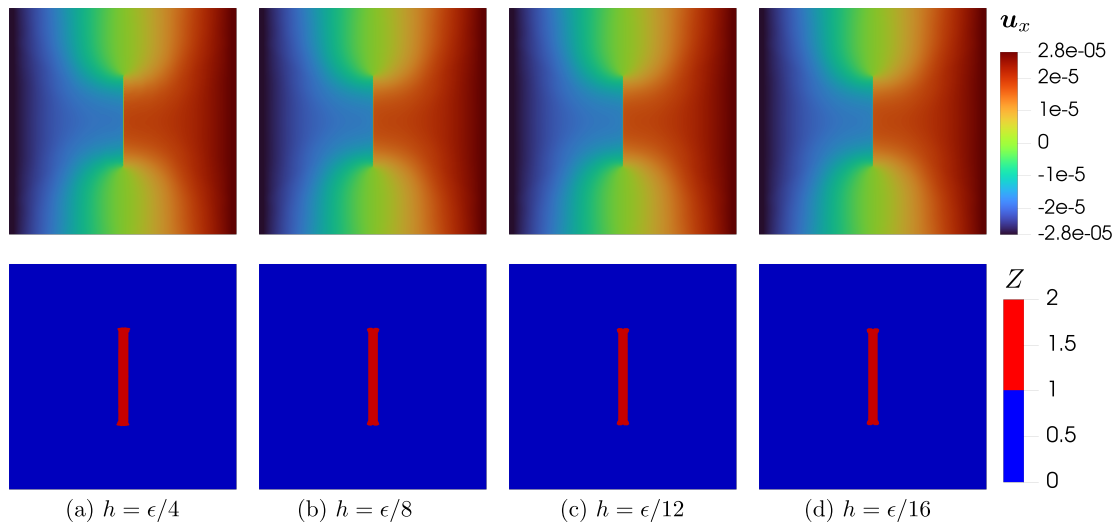}
				\caption{{\bf Mode-I crack problem}: Comparing $\bu_x$ and $Z$ fields for different mesh sizes at time $t = 28\,\mu$s. In the bottom row for the plots of damage field $Z$, note that the red indicates that the point has at least one critically stretched bond; see \cref{rem:plot}}\label{fig:modeICompareMeshResults}
			\end{figure}

			\begin{figure}
				\centering
				\includegraphics[width=\textwidth]{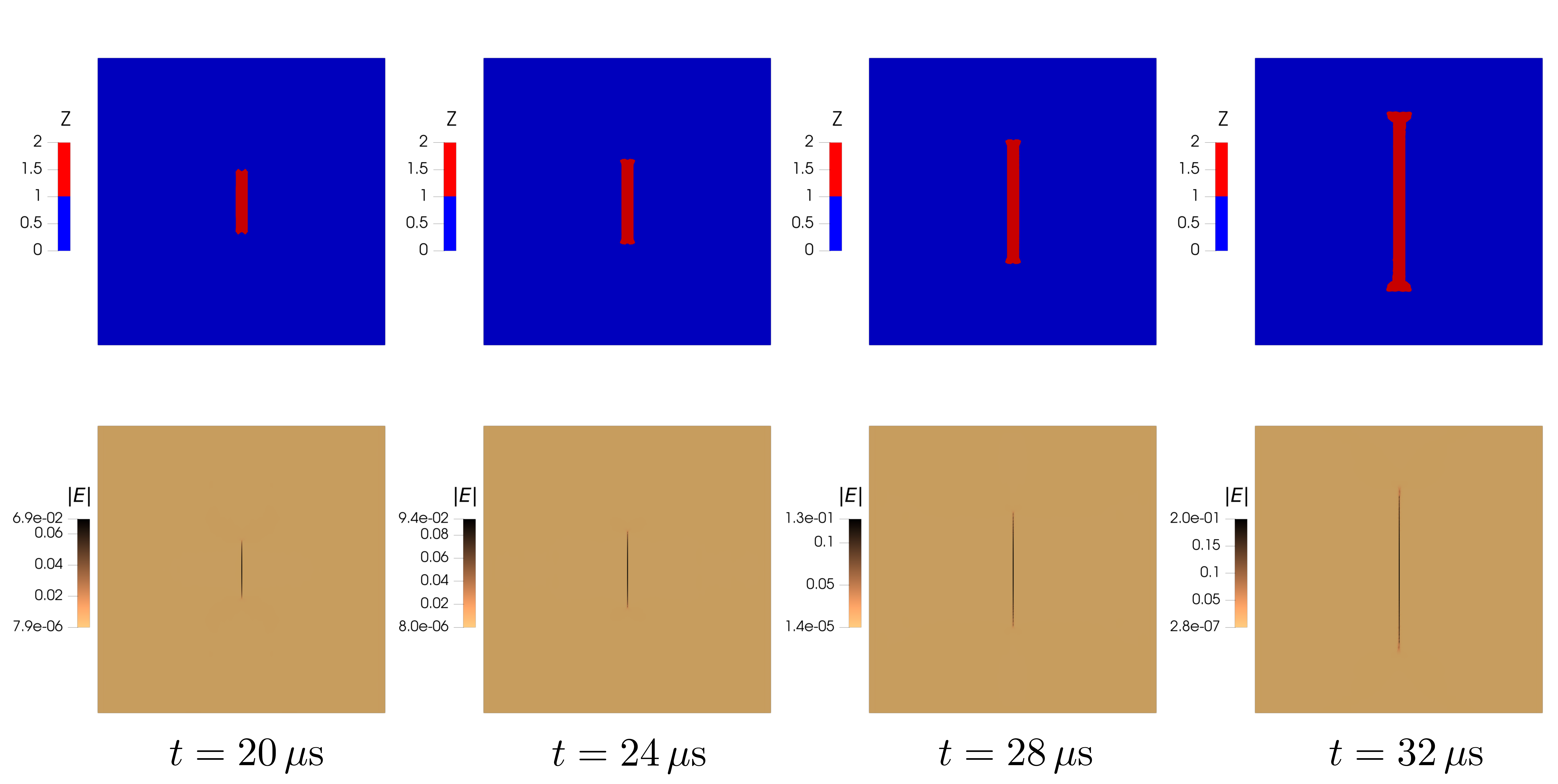}
				\caption{{\bf Mode-I crack problem}: Plot of damage ({\bf top row}) and the magnitude of the strain $\bE = \frac{1}{2}\left[\nabla \bu + \nabla \bu^T\right]$ ({\bf bottom row}). The damage and strains are localized near the crack interface. The thickness of the damage zone in the left column (red region) is approximately twice the horizon. Since the finite element displacement is continuous and piecewise linear, the strain tensor field over each element is constant, and the value is computed at the element's center.}
				\label{fig:mode1plots}
			\end{figure}
			
			\begin{figure}
				\centering
				\includegraphics[width=\textwidth]{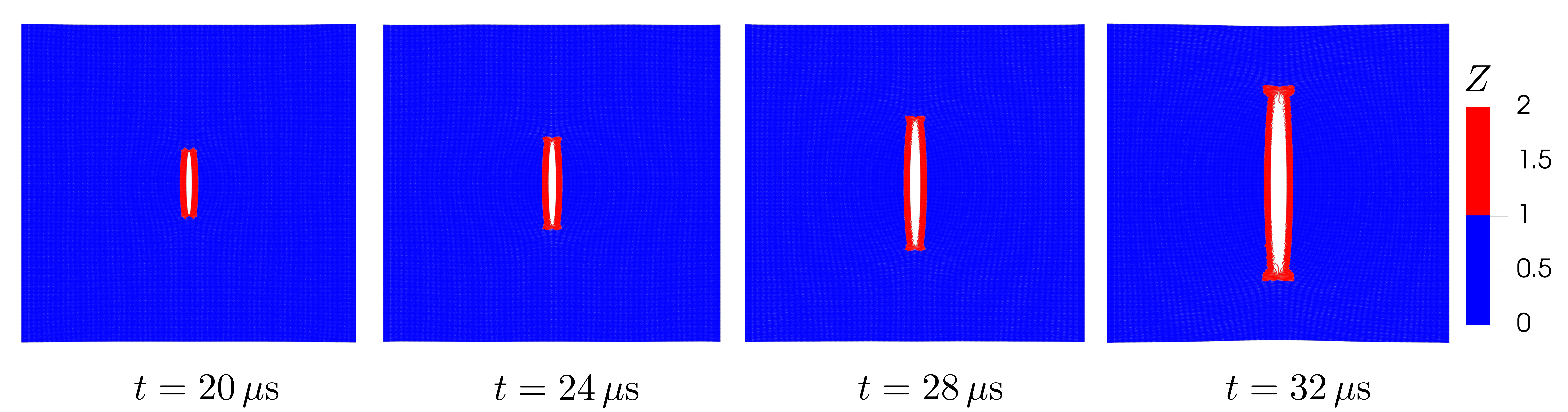}
				\caption{{\bf Mode-I crack problem}: Plot of damage in the deformed configuration after scaling the displacement field by $100$ to show the crack opening. The plots are obtained using the nodal data (without finite element representation) to show a clear opening. This procedure is repeated in other figures where the opening is visualized by scaling the displacement field.}
				\label{fig:mode1plotsZcod}
			\end{figure}
			
			\subsubsection{Comparing NFEA with meshfree method}
			To validate NFEA by comparing with the existing meshfree method employed in \cite{silling2005meshfree, lipton2019complex, jha2019numerical}, the present problem of mode-I fracture is used. Specifically, for the mesh size $h = \epsilon/8$, the peridynamics equation is solved using NFEA and the meshfree method. The rest of the parameters and the setup are the same as above. The plots of the damaged region ($Z(\bx)$) are presented in \cref{fig:compareWithMeshfree}. The figure shows that the NFEA produces results similar to those of the meshfree method. 
			
			\begin{figure}
				\centering
				\includegraphics[width=0.8\textwidth]{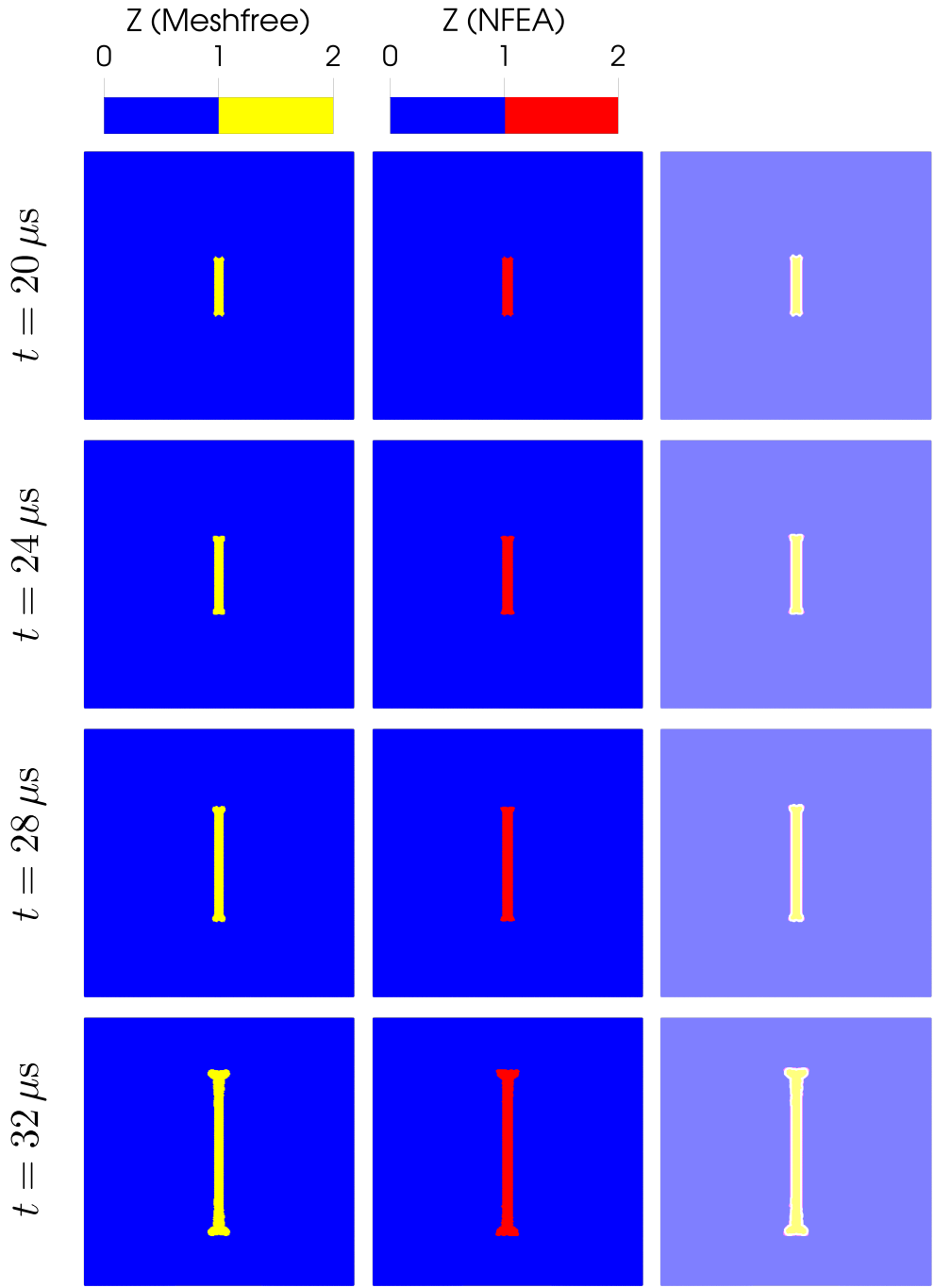}
				\caption{{\bf Mode-I crack problem}: Comparing NFEA and meshfree results for the mode-I fracture problem. Here, the first column plots the damaged zone from the meshfree method, the second column from NFEA, and the overlap of the damaged regions from the two methods is presented in the last column; the opacity of both plots is set to 50\% in the overlap plot. The rows correspond to the results at different times. From the plots in the third column, it is clear that the damaged regions from both methods are quite close.}\label{fig:compareWithMeshfree}
			\end{figure}

			\subsubsection{Localization of damage}
			The mode-I fracture problem from the above is employed to show that NFEA can capture the localization of damage and convergence with respect to the nonlocal length scale. To allow comparisons of results for different nonlocal length scales, in the setup of the mode-I fracture problem in \cref{fig:mode1crack}, the thickness of the left and right layers (where displacement boundary condition is applied) is fixed to $3$ mm, i.e., the thickness of the layers does not change with the nonlocal length scale. Keeping the rest of the setup details the same, including the final time and time step size, the mode-I fracture is simulated for three horizons $\epsilon \in {3, 2, 1}$ mm. Given a horizon $\epsilon$, the mesh size is fixed to $h = \epsilon/4$. In \cref{fig:compareHorizons}, the damaged region $Z(\bx) > 1$ is compared for three horizons at sample times. The results show damage localization as the nonlocal length scale is refined. Most importantly, points on the left of the localized damage zone do not interact with those on the right, and we have formed a piece of crack.
			
			\begin{figure}
				\centering
				\includegraphics[width=0.9\textwidth]{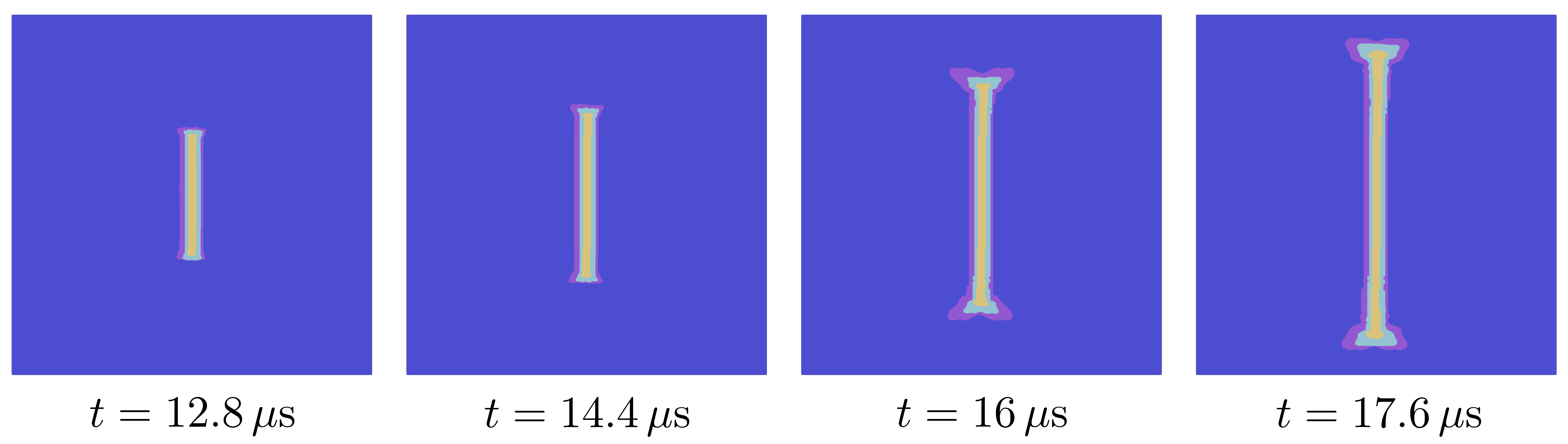}
				\caption{{\bf Mode-I crack problem with fixed thickness of boundary condition layers}: Comparing damaged zone for three different horizons at fixed times. The damaged region $Z(\bx) > 1$ shown in color other than blue is overlapped for three horizons, $\epsilon\in \{3, 2, 1\}$ mm. Light yellow, cyan, and purple regions correspond to $\epsilon = 3$ mm, $\epsilon = 2$ mm, and $\epsilon = 1$ mm, respectively. It is seen that the damaged region corresponding to the smallest horizon is contained inside the damaged region for the largest horizon. This demonstrates the localization of damage as $\epsilon$ is reduced.}
				\label{fig:compareHorizons}
			\end{figure}

			\subsection{Material with a circular hole subjected to an axial loading}\label{ss:voidproblem}
			A material with a hole, as shown in \cref{fig:voidproblem}, is subjected to displacement-controlled axial pulling. The details of the setup and boundary conditions are given in \cref{fig:voidproblem}. The remaining parameters are fixed as follows: horizon $\epsilon = 1$ mm, mesh size $h = 0.25$ mm, final time of the simulation $t_F = 160\, \mu$s, and the size of the time step $\Delta t = 0.0016\, \mu$s. Peridynamics force is computed using the RNP model. 
			
			The damage profile and the strains are presented in \cref{fig:circvoidplots}. Crack nucleation is seen when the internal stresses become large enough. It is also clear that the crack nucleates at the top and bottom edges of the void where the stresses are most significant; see \cref{fig:circvoidplots}(e). To visualize crack opening and branching, the displacement field is magnified by a factor of $50$ in \cref{fig:crackvoidZcod}. Branching of the cracks is also seen at later times.
			
			\begin{figure}
				\centering
				\includegraphics[scale=0.14]{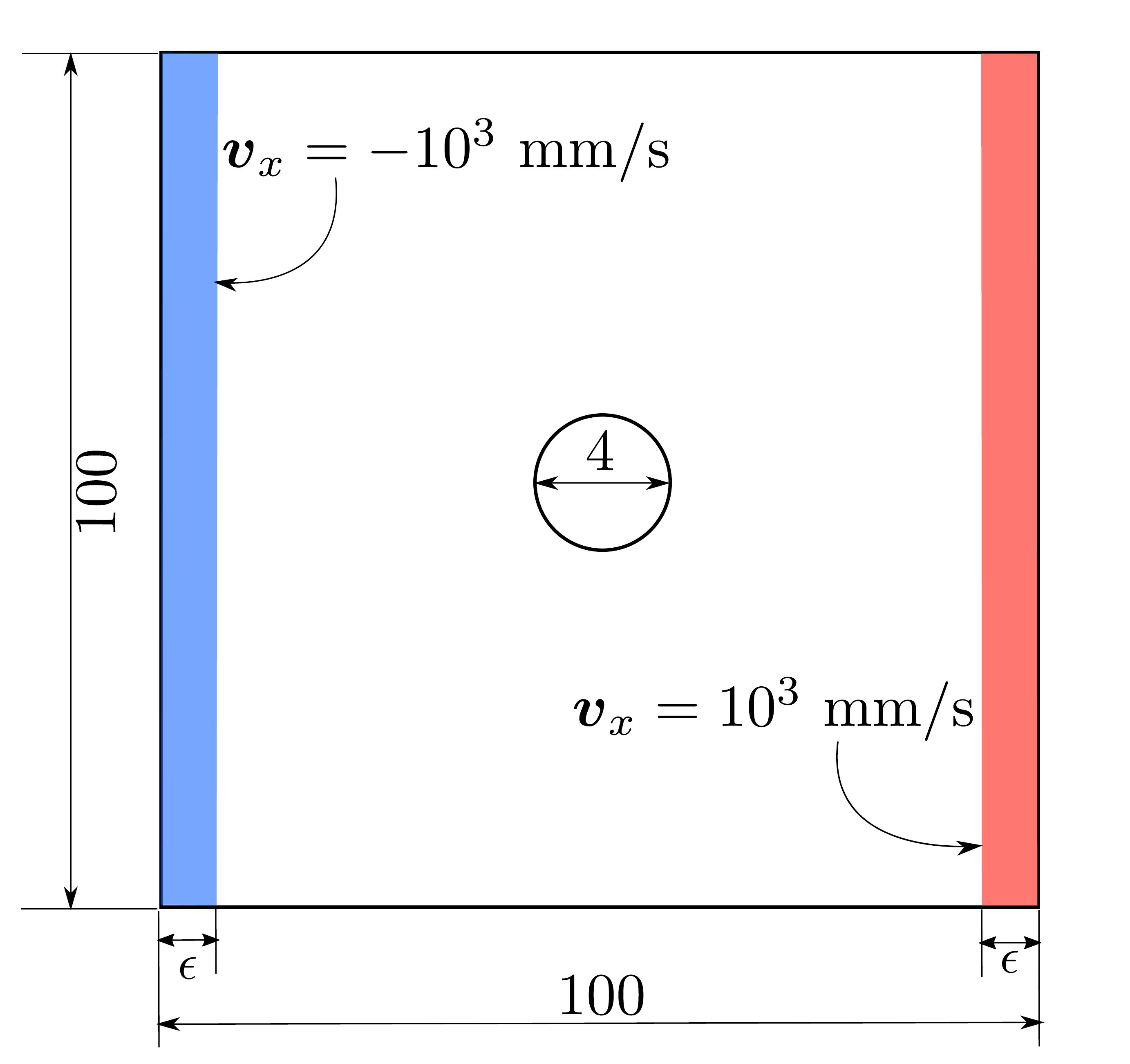}
				\caption{{\bf Circular hole problem}: Setup. The horizon is $\epsilon = 1$ mm. Constant velocity in the opposite direction is specified along the x-axis in the left (blue) and right (red) layers.}
				\label{fig:voidproblem}
			\end{figure}
			
			\begin{figure}
				\centering
				\includegraphics[width=\textwidth]{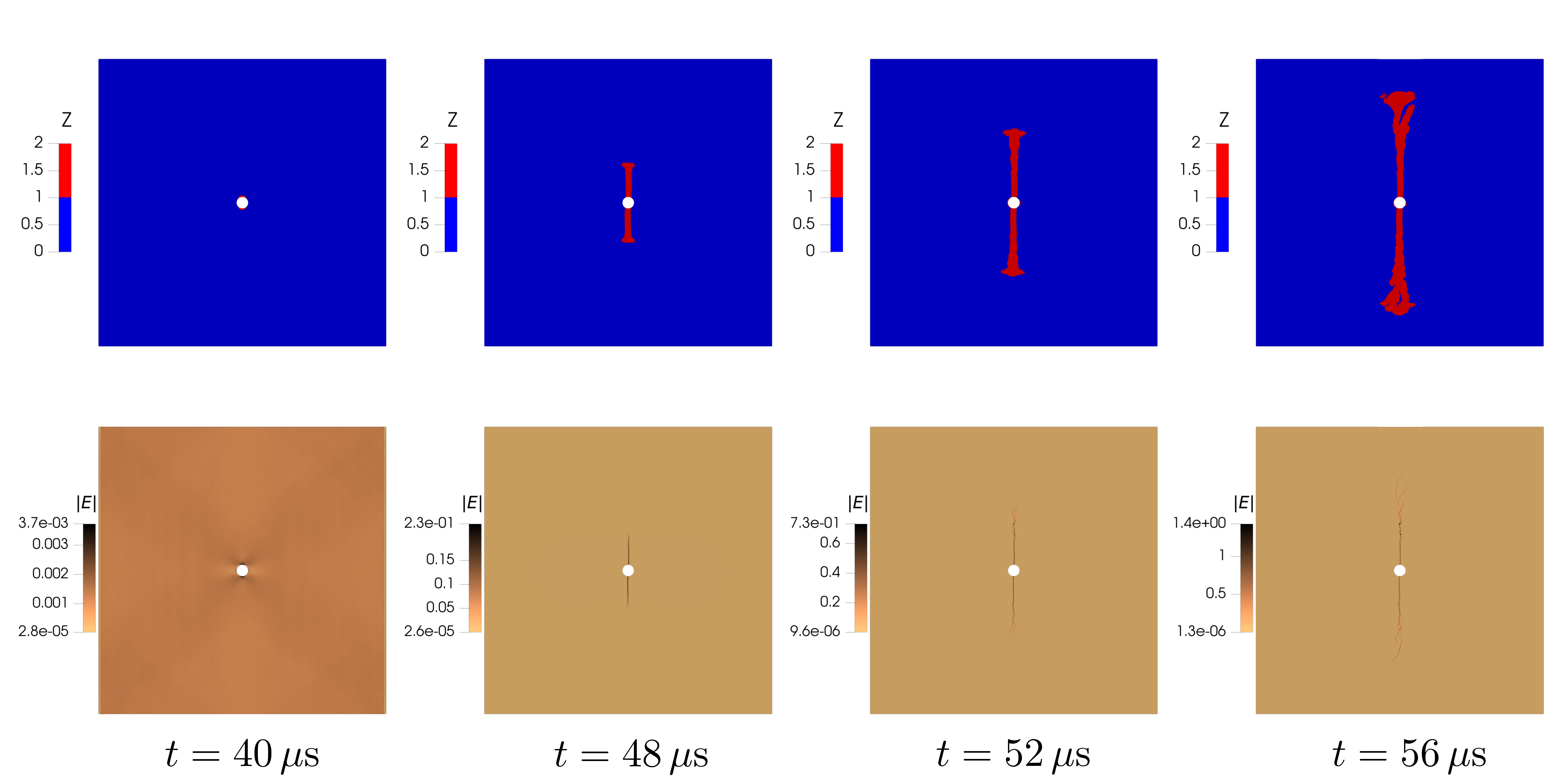}
				\caption{{\bf Circular hole problem}: Plot of damage ({\bf top row}) and the magnitude of the strain $\bE = \frac{1}{2}\left[\nabla \bu + \nabla \bu^T\right]$ ({\bf bottom row}).  A crack is seen nucleating at the two points on the hole's edge where the strain is maximum.}
				\label{fig:circvoidplots}
			\end{figure}
			
			\begin{figure}
				\centering
				\includegraphics[width=\textwidth]{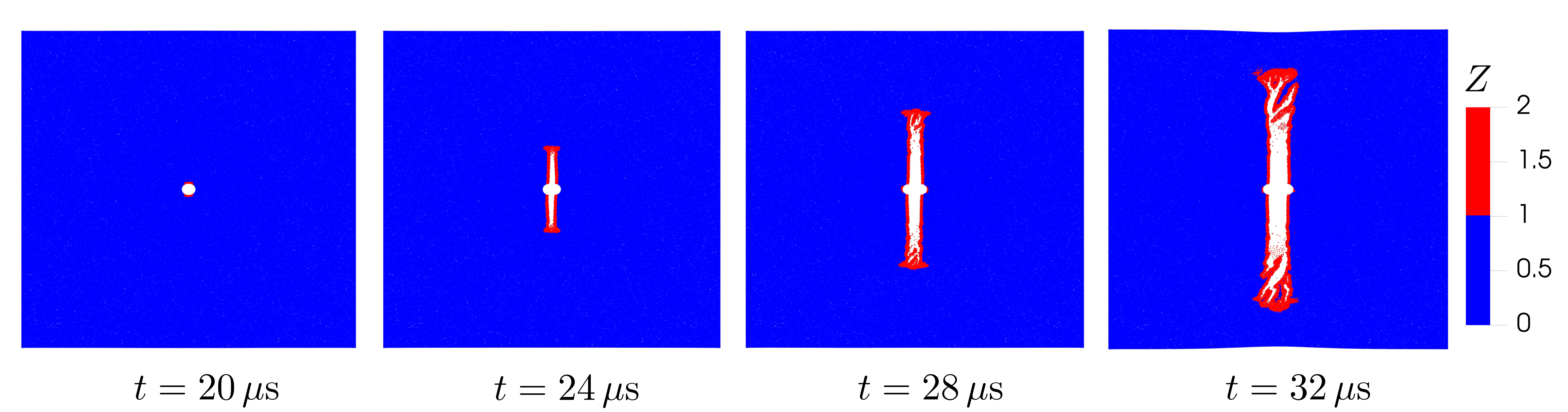}
				\caption{{\bf Circular hole problem}: Plot of damage in the deformed configuration after scaling the displacement field by $50$ to show the crack opening. }
				\label{fig:circvoidplotsZcod}
			\end{figure}
			
			\subsection{Material with a v-notch under bending load}\label{ss:bendingproblem}
			In this example, a rectangle beam with a v-notch is subjected to the bending load as shown in \cref{fig:vnotchsetup}. The horizon is fixed to $\epsilon = 1$ mm, mesh size $h = 0.25$ mm, final simulation time $t_F = 250\, \mu$s, and the size of the time step $\Delta t = 0.001667\, \mu$s. Peridynamics force is based on the RNP model. 
			
			The damage profile and the magnitude of the strain are presented in \cref{fig:vnotch}. As expected, the crack nucleates at the tip of the notch where the strain is most significant. Similar to the previous example, displacement is magnified by $50$ to highlight the separation of the structure in \cref{fig:vnotchZcod}.
			
			\begin{figure}
				\centering
				\includegraphics[scale=0.2]{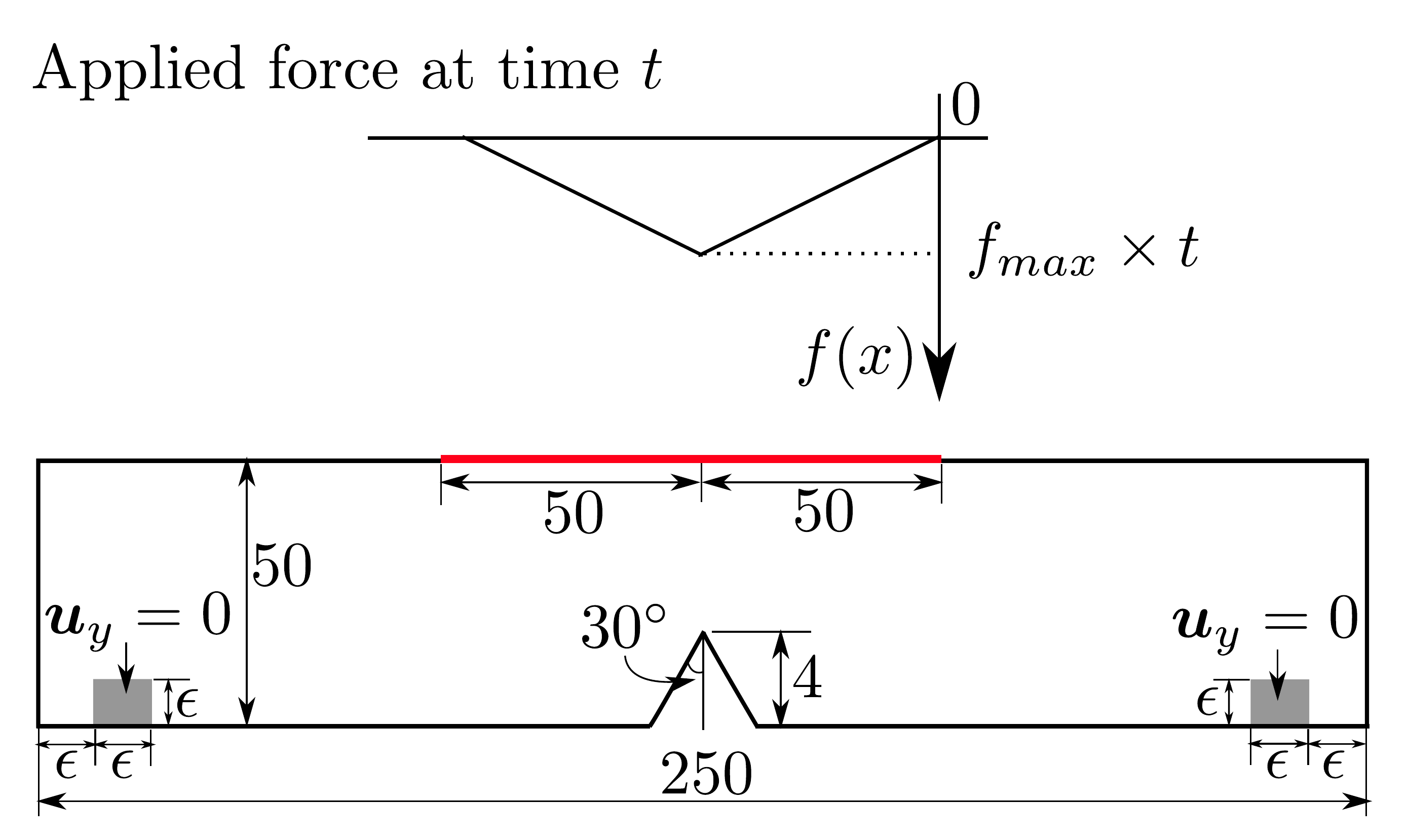}
				\caption{{\bf V-notch problem}: Setup. The horizon in this problem is $\epsilon = 1$ mm. Vertically downward distributed force is applied on the part of the top edge, as shown in red. The profile of distributed load is shown above the red line, where the loading parameter $f_{max}$ is set to $f_{max} = 2.5\times 10^{5}$ N/($\mu$s$\cdot$mm).}
				\label{fig:vnotchsetup}
			\end{figure}

			\begin{figure}
				\centering
				\includegraphics[width=\textwidth]{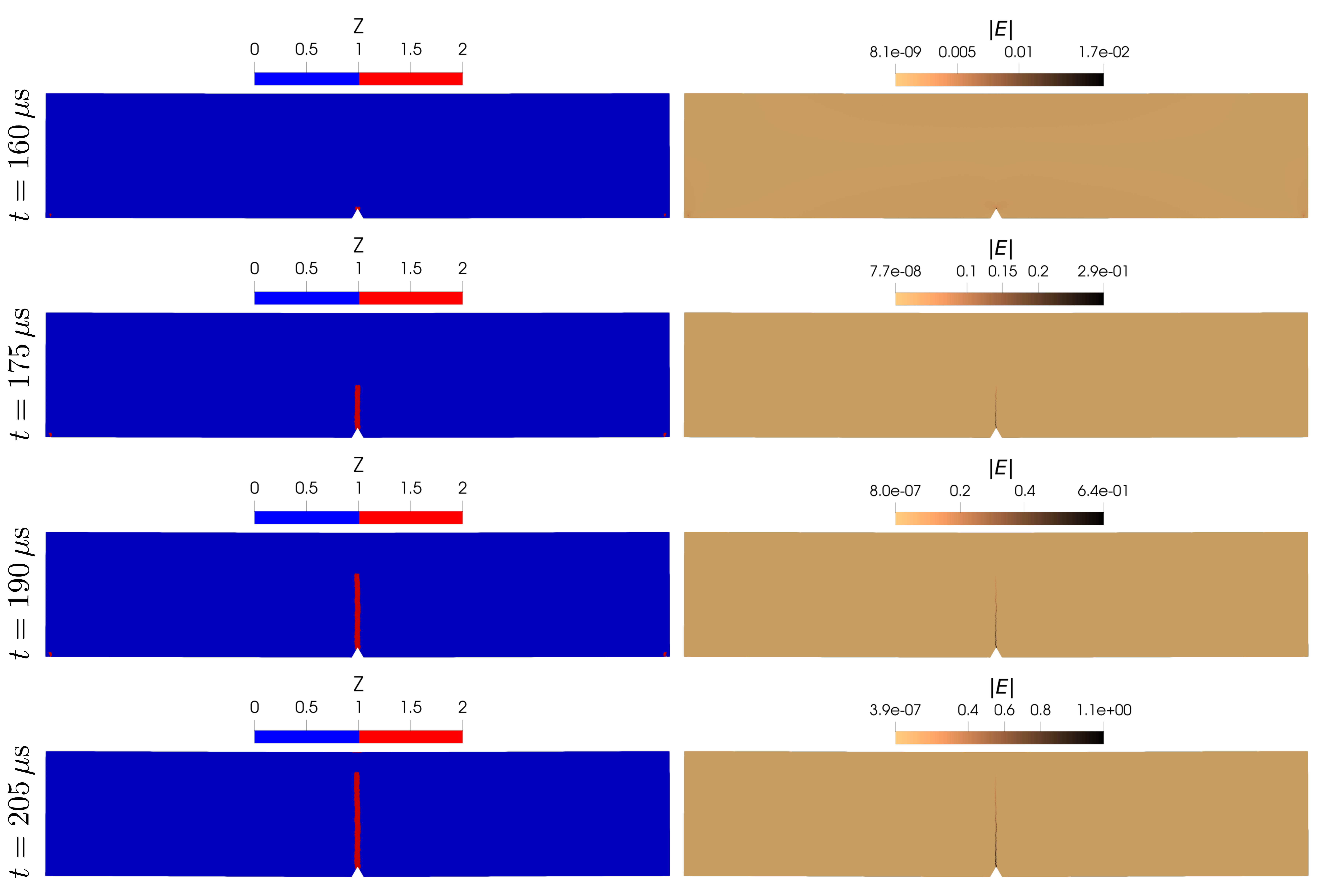}
				\caption{{\bf V-notch problem}: Plot of damage and the magnitude of the strain. Crack nucleates at the tip of the notch where the strain is maximum.}\label{fig:vnotch}
			\end{figure}
			
			\begin{figure}
				\centering
				\includegraphics[width=\textwidth]{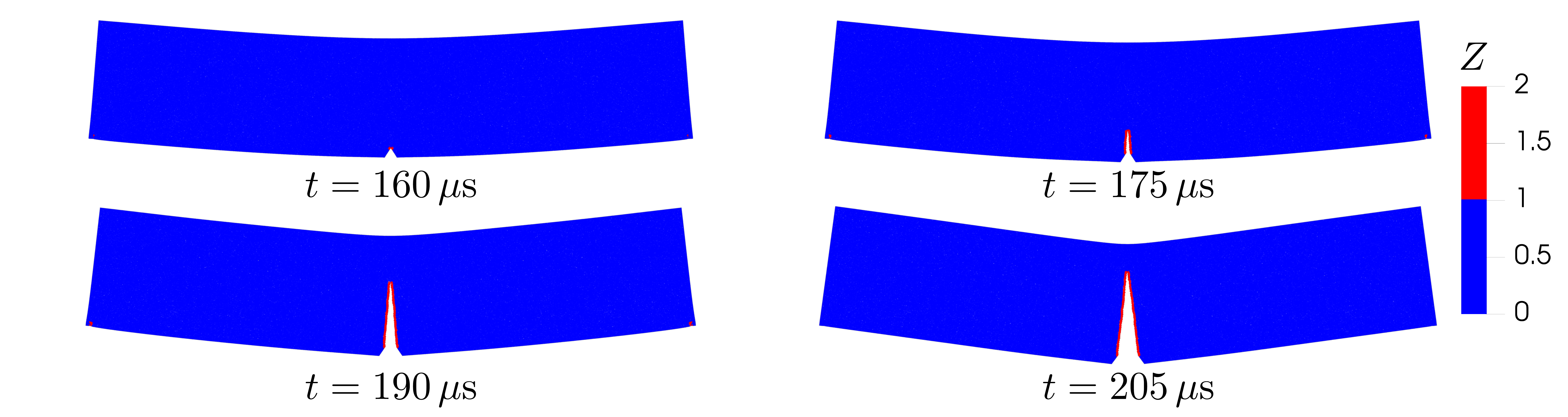}
				\caption{{\bf V-notch problem}: Plot of damage in the deformed configuration after scaling the displacement field by $50$ to highlight the separation of the structure.}
				\label{fig:vnotchZcod}
			\end{figure}
			
			\subsection{Material with a circular hole and pre-crack}\label{ss:crackvoidproblem}
			As a final example, a rectangular domain with an existing horizontal pre-crack and a circular hole in the neighborhood of a crack is considered. The setup, as shown in \cref{fig:crackvoidproblem}, is motivated by a similar example in \cite{CMPer-DaiBFEM}. The horizon is fixed to $\epsilon = 0.4$ mm, mesh size $h = 0.1$ mm, the final simulation time $t_F = 800\, \mu$s, and the size of the time step $\Delta t = 0.004\, \mu$s. The peridynamics force is computed using the RNP model. 
			
			Damage and the magnitude of the strain at different times are displayed in \cref{fig:crackvoidZ}. Initially, crack propagation is influenced by the hole nearby, and instead of growing horizontally, it is deflected. At later times, when the crack tip moves past the hole, the crack propagates horizontally. A similar problem was considered in \cite[Figure 18]{CMPer-DaiBFEM}, where results using different numerical methods were compared. The results of this work qualitatively agree with that in \cite{CMPer-DaiBFEM}. In \cref{fig:crackvoidStrain}, strain fields at different times are presented. Opening of the crack is visualized in \cref{fig:crackvoidZcod}.
			
			\begin{figure}
				\centering
				\includegraphics[scale=0.16]{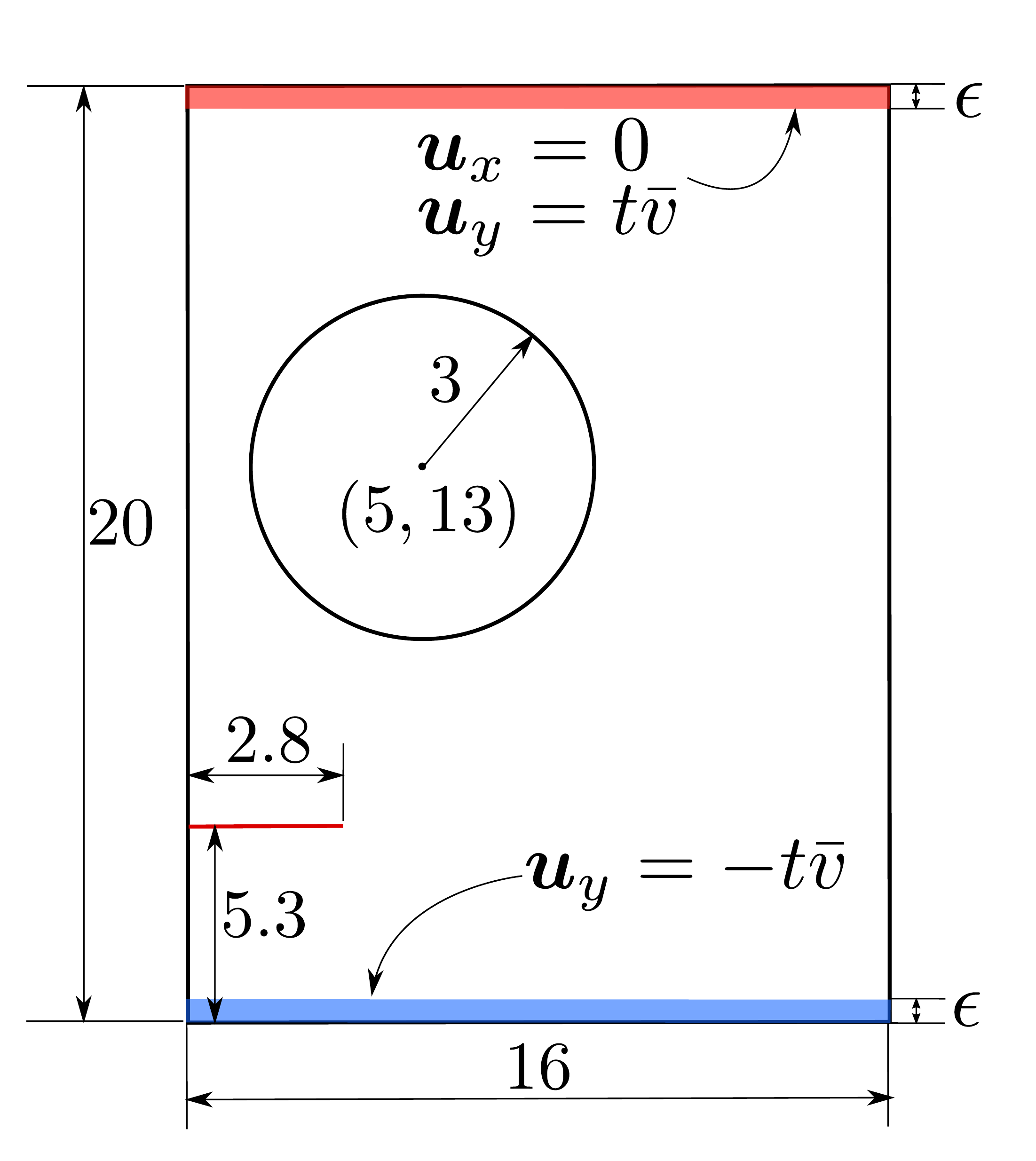}
				\caption{{\bf Circular hole and pre-crack problem}: Setup. The horizon is $\epsilon = 0.4$ mm. The magnitude of the prescribed vertical velocity on the top and bottom layers is $\bar{v} = 25$ mm/s.}
				\label{fig:crackvoidproblem}
			\end{figure}
			
			\begin{figure}
				\centering
				\includegraphics[scale=0.16]{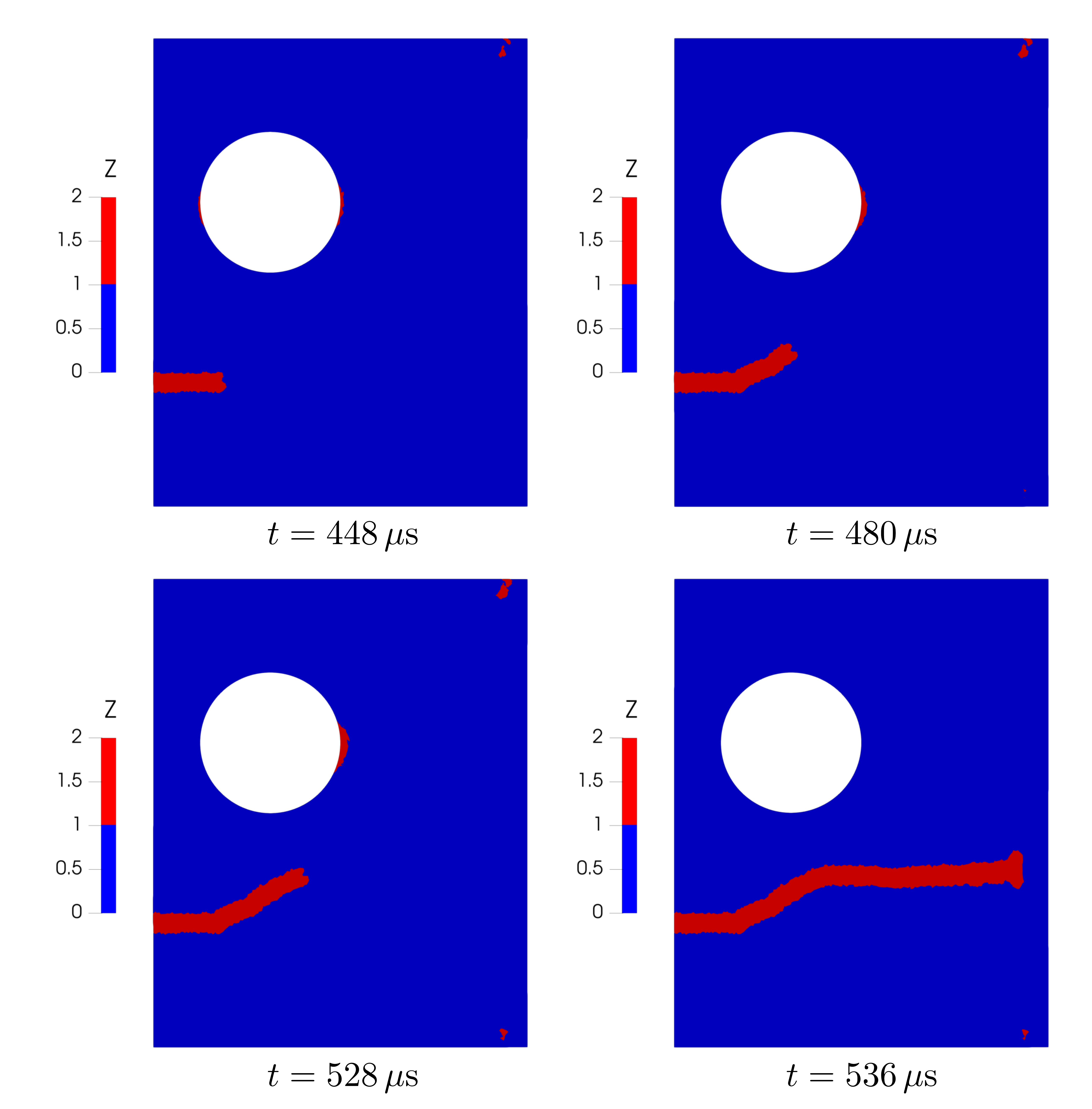}
				\caption{{\bf Circular hole and pre-crack problem}: Plot of damage.}\label{fig:crackvoidZ}
			\end{figure}
			
			\begin{figure}
				\centering
				\includegraphics[scale=0.16]{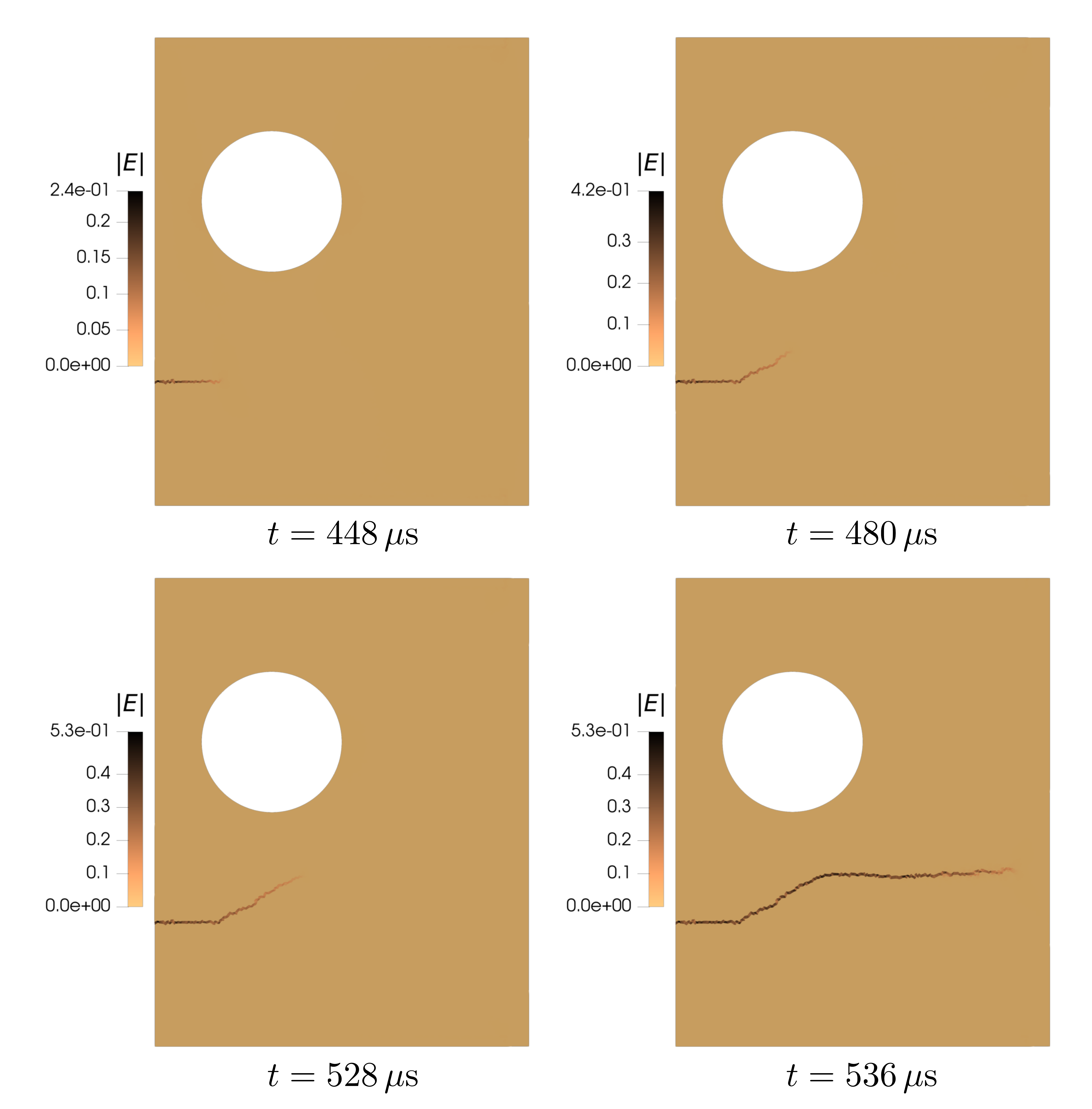}
				\caption{{\bf Circular hole and pre-crack problem}: Plot of the magnitude of strain $\bE = \frac{1}{2}\left[\nabla \bu + \nabla \bu^T\right]$.}\label{fig:crackvoidStrain}
			\end{figure}
			
			\begin{figure}
				\centering
				\includegraphics[scale=0.16]{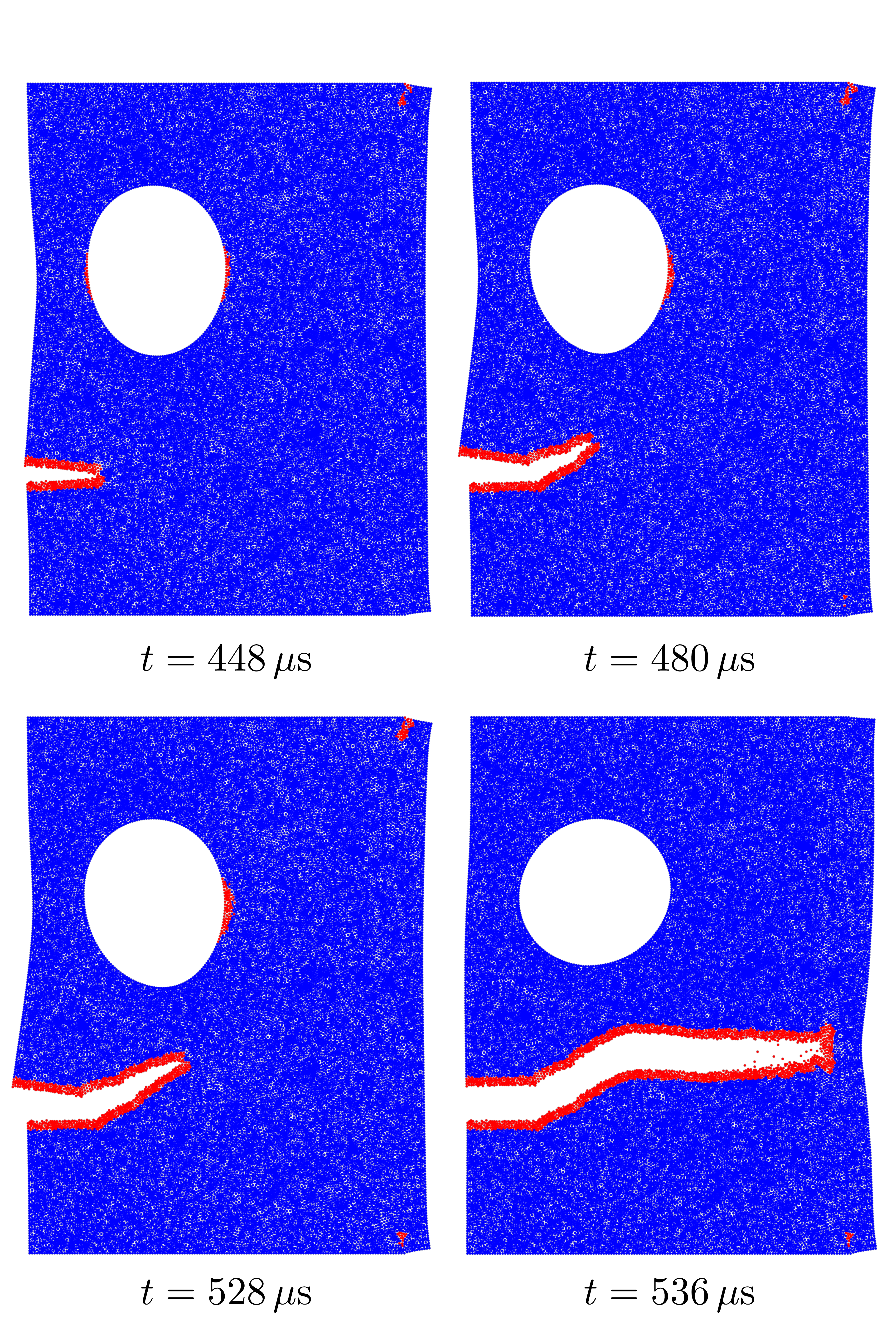}
				\caption{{\bf Circular hole and pre-crack problem}: Plot of damage in the deformed configuration after scaling the displacement field by $50$ to show the crack opening.}
				\label{fig:crackvoidZcod}
			\end{figure}
			
			\subsection{On the crack propagation speed}
			In this subsection, the crack propagation speeds from four problems are compared. Let $t_1$ and $t_2$ be times when the crack begins and stops propagating, respectively. Also, let $v(t)$, for $t\in [t_1, t_2]$, be the crack speed computed from the simulation at time $t$. To plot the crack speeds for all four examples in one plot, time $t\in [t_1, t_2]$ is transformed to $\bar{t} = (t - t_1)/(t_2 - t_1)$ so that $\bar{t} \in [0,1]$. Let $\bar{v}(\bar{t}) = v(t)$ be the crack speed as a function of normalized time $\bar{t}$. Next, crack speed is normalized by dividing the Rayleigh wave speed $c_R$; $c_R$ from \cref{tab:matprops} is $3244.2$ m/s.

			\begin{figure}
				\centering
				\includegraphics[scale=0.4]{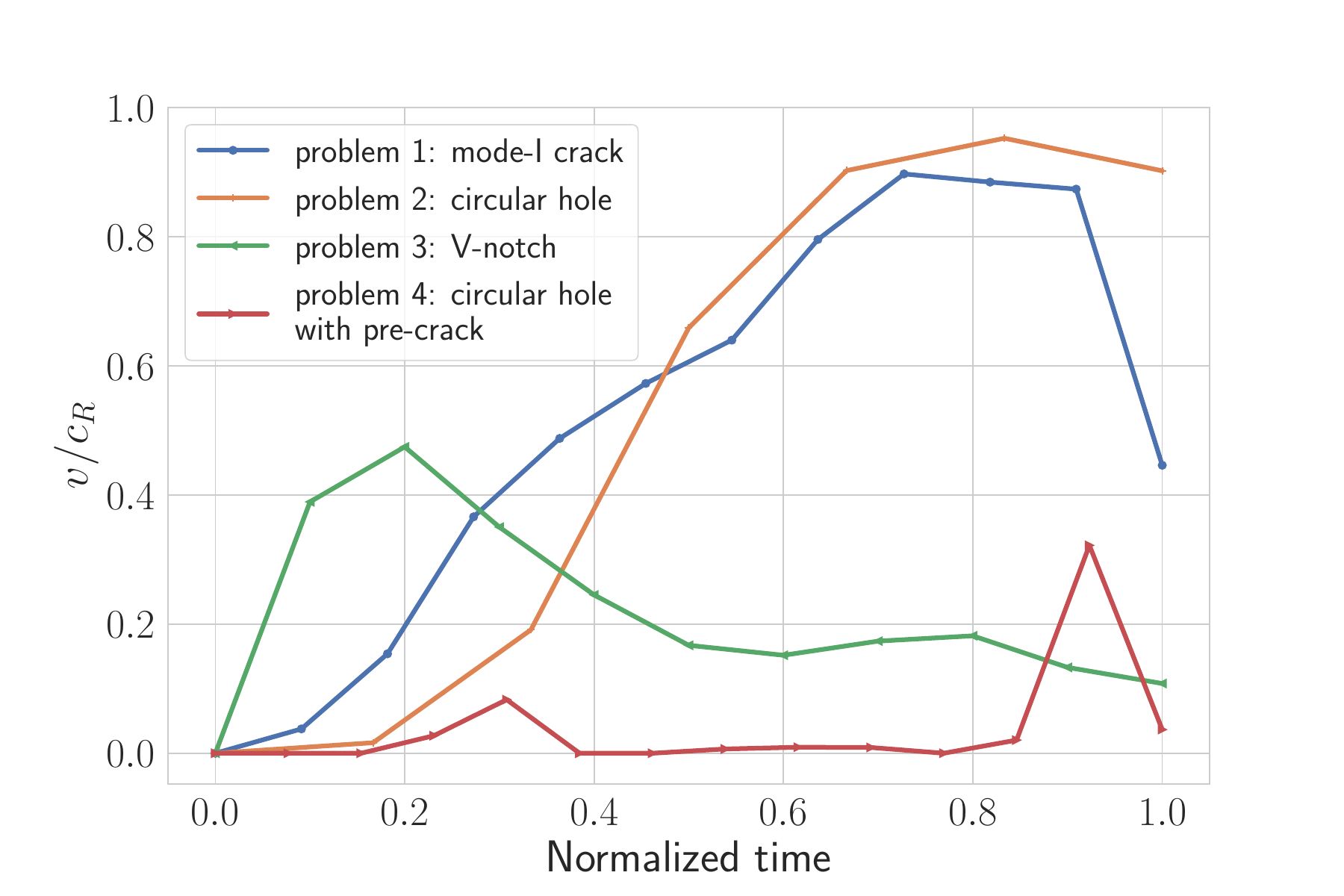}
				\caption{Comparison of the normalized crack speed for the four problems.}
				\label{fig:compareCrackSpeed}
			\end{figure}
			
			\begin{table}
				\centering
				\addtolength{\tabcolsep}{+1pt}
				\renewcommand{\arraystretch}{1.4}
				\begin{tabular}{|l||c|c|c|l||c|c|}  
					\hline
					Problem type & $\max \left(\frac{v}{c_R}\right)$ & $\mathrm{avg} \left(\frac{v}{c_R}\right)$ & & Problem type & $\max \left(\frac{v}{c_R}\right)$ & $\mathrm{avg} \left(\frac{v}{c_R}\right)$ \\
					\hline\hline
					Mode-I crack & 0.9 & 0.51 & & Circular hole & 0.95 & 0.52 \\
					\hline
					V-notch & 0.47 & 0.22 & & Circular hole with pre-crack & 0.32 & 0.04 \\
					\hline
				\end{tabular}
				
				\caption{Maximum and average normalized crack speeds for the four numerical problems. }
				\label{tab:compareCrackSpeeds}
			\end{table}
			
			\cref{fig:compareCrackSpeed} presents the normalized crack speed as a function of normalized time for the four problems. As expected, the normalized crack speeds are below 1, i.e., the crack propagates slower than the Rayleigh wave speed; see \cref{tab:compareCrackSpeeds}, which lists the maximum and average values of normalized crack speeds.

			%
			\section{Conclusion}\label{s:concl}
			This work analyzed the nodal finite element approximation for the peridynamics. Assuming exact solutions are in proper function spaces, consistency errors are shown to be bounded, and a-priori convergence of the discretization is established. The nodal finite element discretization implementation is discussed in detail, and a range of numerical experiments are performed using the method to show the utility of the approximation. The nodal finite element approximation is relatively straightforward and can be easily integrated with the standard finite element meshing libraries. Further, the method is computationally faster than the standard finite element approximation because the mass matrix is diagonal, and the nonlocal force calculation is similar to finite-difference/mesh-free approximation. Since NFEA is based on finite element representation and mesh, the coupling peridynamics with other PDE-based models for multiphysics simulation is straightforward.
			
			The work also presents an alternative NFEA method based on Cl'ement interpolation. This work considers Cl'ement interpolation only theoretically. The a-priori error estimates in the alternative NFEA improve the NFEA a-priori error estimates, and the regularity requirement on the exact solution in Cl\'ement NFEA is less restrictive. Future work will explore implementing the alternative NFEA method and developing a-posteriori error estimates. 
			
			%
			\section*{Acknowledgements}
			The majority of the work is done through the support of the U.S. Army Research Laboratory and the U.S. Army Research Office under contract/grant number W911NF1610456 RL. PKJ is also thankful to the Oden Institute for Computational Engineering and Sciences, The University of Texas at Austin, for providing resources to run some simulations. PD is thankful to the LSU Center of Computation \& Technology for supporting this work. 
			
			

\begin{thebibliography}{}
				
				\bibitem[Agwai et~al., 2011]{CMPer-Agwai}
				Agwai, A., Guven, I., and Madenci, E. (2011).
				\newblock Predicting crack propagation with peridynamics: a comparative study.
				\newblock {\em International journal of fracture}, 171(1):65--78.
				
				\bibitem[Ahrens et~al., 2005]{paraview}
				Ahrens, J., Geveci, B., and Law, C. (2005).
				\newblock Paraview: An end-user tool for large data visualization.
				\newblock {\em The visualization handbook}, 717.
				
				\bibitem[Aksoylu and Unlu, 2014]{AksoyluUnlu}
				Aksoylu, B. and Unlu, Z. (2014).
				\newblock Conditioning analysis of nonlocal integral operators in fractional
				sobolev spaces.
				\newblock {\em SIAM Journal on Numerical Analysis}, 52:653--677.
				
				\bibitem[Ambrosio and Brades, 1997]{AmbrosioBrades}
				Ambrosio, L. and Brades, A. (1997).
				\newblock {\em Energies in SBV and variational models in fracture mechanics:
					Homogenization and Applications to Materials Science}, volume~9.
				\newblock Gakuto, Gakkotosho, Tokyo, Japan.
				
				\bibitem[Ambrosio et~al., 1997]{Ambrosio}
				Ambrosio, L., Coscia, A., and Dal~Maso, G. (1997).
				\newblock Fine properties of functions with bounded deformation.
				\newblock {\em Archive for Rational Mechanics and Analysis}, 139(3):201--238.
				
				\bibitem[Anicode and Madenci, 2022]{anicode2022bond}
				Anicode, S. V.~K. and Madenci, E. (2022).
				\newblock Bond-and state-based peridynamic analysis in a commercial finite
				element framework with native elements.
				\newblock {\em Computer Methods in Applied Mechanics and Engineering},
				398:115208.
				
				\bibitem[Bobaru and Hu, 2012]{BobaruHu}
				Bobaru, F. and Hu, W. (2012).
				\newblock The meaning, selection, and use of the peridynamic horizon and its
				relation to crack branching in brittle materials.
				\newblock {\em International journal of fracture}, 176(2):215--222.
				
				\bibitem[Brenner and Scott, 2007]{MANa-Susanne}
				Brenner, S. and Scott, R. (2007).
				\newblock {\em The mathematical theory of finite element methods}, volume~15.
				\newblock Springer Science \& Business Media, 3 edition.
				
				\bibitem[Chen and Gunzburger, 2011]{chen2011continuous}
				Chen, X. and Gunzburger, M. (2011).
				\newblock Continuous and discontinuous finite element methods for a
				peridynamics model of mechanics.
				\newblock {\em Computer Methods in Applied Mechanics and Engineering},
				200(9-12):1237--1250.
				
				\bibitem[Cl\'ement, 1975]{clement1975}
				Cl\'ement, P. (1975).
				\newblock Approximation by finite element functions using local regularization.
				\newblock {\em RARO Anal. Num\'er}, 9:77--84.
				
				\bibitem[Dai et~al., 2015]{CMPer-DaiBFEM}
				Dai, S., Augarde, C., Du, C., and Chen, D. (2015).
				\newblock A fully automatic polygon scaled boundary finite element method for
				modelling crack propagation.
				\newblock {\em Engineering Fracture Mechanics}, 133:163--178.
				
				\bibitem[De~Meo and Oterkus, 2017]{de2017finite}
				De~Meo, D. and Oterkus, E. (2017).
				\newblock Finite element implementation of a peridynamic pitting corrosion
				damage model.
				\newblock {\em Ocean Engineering}, 135:76--83.
				
				\bibitem[Diehl et~al., 2020]{diehl2020asynchronous}
				Diehl, P., Jha, P.~K., Kaiser, H., Lipton, R., and L{\'e}vesque, M. (2020).
				\newblock An asynchronous and task-based implementation of peridynamics
				utilizing hpx—the c++ standard library for parallelism and concurrency.
				\newblock {\em SN Applied Sciences}, 2(12):1--21.
				
				\bibitem[Diehl et~al., 2019]{diehl2019review}
				Diehl, P., Prudhomme, S., and L{\'e}vesque, M. (2019).
				\newblock A review of benchmark experiments for the validation of peridynamics
				models.
				\newblock {\em Journal of Peridynamics and Nonlocal Modeling}, 1:14--35.
				
				\bibitem[Diyaroglu et~al., 2017]{diyaroglu2017peridynamic}
				Diyaroglu, C., Oterkus, S., Oterkus, E., and Madenci, E. (2017).
				\newblock Peridynamic modeling of diffusion by using finite-element analysis.
				\newblock {\em IEEE Transactions on Components, Packaging and Manufacturing
					Technology}, 7(11):1823--1831.
				
				\bibitem[Du et~al., 2018]{DuTaoTian}
				Du, Q., Tao, Y., and Tian, X. (2018).
				\newblock A peridynamic model of fracture mechanics with bond-breaking.
				\newblock {\em Journal of Elasticity}, 132(2):197--218.
				
				\bibitem[Emmrich et~al., 2013]{CMPer-Emmrich}
				Emmrich, E., Lehoucq, R.~B., and Puhst, D. (2013).
				\newblock Peridynamics: a nonlocal continuum theory.
				\newblock In {\em Meshfree Methods for Partial Differential Equations VI},
				pages 45--65. Springer.
				
				\bibitem[Foster et~al., 2011]{CMPer-Silling7}
				Foster, J.~T., Silling, S.~A., and Chen, W. (2011).
				\newblock An energy based failure criterion for use with peridynamic states.
				\newblock {\em International Journal for Multiscale Computational Engineering},
				9(6).
				
				\bibitem[Geuzaine and Remacle, 2009]{Gmsh}
				Geuzaine, C. and Remacle, J.-F. (2009).
				\newblock Gmsh: A 3-d finite element mesh generator with built-in pre-and
				post-processing facilities.
				\newblock {\em International journal for numerical methods in engineering},
				79(11):1309--1331.
				
				\bibitem[Ghajari et~al., 2014]{CMPer-Ghajari}
				Ghajari, M., Iannucci, L., and Curtis, P. (2014).
				\newblock A peridynamic material model for the analysis of dynamic crack
				propagation in orthotropic media.
				\newblock {\em Computer Methods in Applied Mechanics and Engineering},
				276:431--452.
				
				\bibitem[Ha and Bobaru, 2010]{HaBobaru}
				Ha, Y.~D. and Bobaru, F. (2010).
				\newblock Studies of dynamic crack propagation and crack branching with
				peridynamics.
				\newblock {\em International Journal of Fracture}, 162(1-2):229--244.
				
				\bibitem[Huang et~al., 2019]{huang2019finite}
				Huang, X., Bie, Z., Wang, L., Jin, Y., Liu, X., Su, G., and He, X. (2019).
				\newblock Finite element method of bond-based peridynamics and its abaqus
				implementation.
				\newblock {\em Engineering Fracture Mechanics}, 206:408--426.
				
				\bibitem[Jha and Lipton, 2021]{jha2021finite}
				Jha, P. and Lipton, R. (2021).
				\newblock Finite element approximation of nonlocal dynamic fracture models.
				\newblock {\em Discrete \& Continuous Dynamical Systems-B}, 26(3):1675.
				
				\bibitem[Jha et~al., 2021]{jha2020peridynamics}
				Jha, P.~K., Desai, P.~S., Bhattacharya, D., and Lipton, R. (2021).
				\newblock Peridynamics-based discrete element method (peridem) model of
				granular systems involving breakage of arbitrarily shaped particles.
				\newblock {\em Journal of the Mechanics and Physics of Solids}, 151:104376.
				
				\bibitem[Jha and Diehl, 2021]{jha2021nlmech}
				Jha, P.~K. and Diehl, P. (2021).
				\newblock Nlmech: Implementation of finite difference/meshfree discretization
				of nonlocal fracture models.
				\newblock {\em Journal of Open Source Software}, 6(65):3020.
				
				\bibitem[Jha and Lipton, 2018a]{jha2018numerical}
				Jha, P.~K. and Lipton, R. (2018a).
				\newblock Numerical analysis of nonlocal fracture models in holder space.
				\newblock {\em SIAM Journal on Numerical Analysis}, 56(2):906--941.
				
				\bibitem[Jha and Lipton, 2018b]{jha2018numerical2}
				Jha, P.~K. and Lipton, R. (2018b).
				\newblock Numerical convergence of nonlinear nonlocal continuum models to local
				elastodynamics.
				\newblock {\em International Journal for Numerical Methods in Engineering},
				114(13):1389--1410.
				
				\bibitem[Jha and Lipton, 2019]{jha2019numerical}
				Jha, P.~K. and Lipton, R. (2019).
				\newblock Numerical convergence of finite difference approximations for state
				based peridynamic fracture models.
				\newblock {\em Computer Methods in Applied Mechanics and Engineering},
				351:184--225.
				
				\bibitem[Jha and Lipton, 2020a]{jha2020finite}
				Jha, P.~K. and Lipton, R. (2020a).
				\newblock Finite element convergence for state-based peridynamic fracture
				models.
				\newblock {\em Communications on Applied Mathematics and Computation},
				2(1):93--128.
				
				\bibitem[Jha and Lipton, 2020b]{Jha2020peri}
				Jha, P.~K. and Lipton, R.~P. (2020b).
				\newblock Kinetic relations and local energy balance for lefm from a nonlocal
				peridynamic model.
				\newblock {\em International Journal of Fracture}.
				
				\bibitem[Lipton, 2014]{CMPer-Lipton3}
				Lipton, R. (2014).
				\newblock Dynamic brittle fracture as a small horizon limit of peridynamics.
				\newblock {\em Journal of Elasticity}, 117(1):21--50.
				
				\bibitem[Lipton, 2016]{CMPer-Lipton}
				Lipton, R. (2016).
				\newblock Cohesive dynamics and brittle fracture.
				\newblock {\em Journal of Elasticity}, 124(2):143--191.
				
				\bibitem[Lipton et~al., 2016]{CMPer-Lipton2}
				Lipton, R., Silling, S., and Lehoucq, R. (2016).
				\newblock Complex fracture nucleation and evolution with nonlocal
				elastodynamics.
				\newblock {\em arXiv preprint arXiv:1602.00247}.
				
				\bibitem[Lipton et~al., 2019]{lipton2019complex}
				Lipton, R.~P., Lehoucq, R.~B., and Jha, P.~K. (2019).
				\newblock Complex fracture nucleation and evolution with nonlocal
				elastodynamics.
				\newblock {\em Journal of Peridynamics and Nonlocal Modeling}, 1(2):122--130.
				
				\bibitem[Liu and Hong, 2012]{liu2012coupling}
				Liu, W. and Hong, J.-W. (2012).
				\newblock A coupling approach of discretized peridynamics with finite element
				method.
				\newblock {\em Computer methods in applied mechanics and engineering},
				245:163--175.
				
				\bibitem[Macek and Silling, 2007]{CMPer-Richard}
				Macek, R.~W. and Silling, S.~A. (2007).
				\newblock Peridynamics via finite element analysis.
				\newblock {\em Finite Elements in Analysis and Design}, 43(15):1169--1178.
				
				\bibitem[Madenci et~al., 2018]{CMPer-Madenci}
				Madenci, E., Dorduncu, M., Barut, A., and Phan, N. (2018).
				\newblock A state-based peridynamic analysis in a finite element framework.
				\newblock {\em Engineering Fracture Mechanics}, 195:104--128.
				
				\bibitem[Mengesha and Du, 2015]{CMPer-Mengesha2}
				Mengesha, T. and Du, Q. (2015).
				\newblock On the variational limit of a class of nonlocal functionals related
				to peridynamics.
				\newblock {\em Nonlinearity}, 28(11):3999.
				
				\bibitem[Ni et~al., 2018]{ni2018peridynamic}
				Ni, T., Zhu, Q.-z., Zhao, L.-Y., and Li, P.-F. (2018).
				\newblock Peridynamic simulation of fracture in quasi brittle solids using
				irregular finite element mesh.
				\newblock {\em Engineering Fracture Mechanics}, 188:320--343.
				
				\bibitem[Royer and Clorennec, 2007]{royer2007improved}
				Royer, D. and Clorennec, D. (2007).
				\newblock An improved approximation for the rayleigh wave equation.
				\newblock {\em Ultrasonics}, 46(1):23--24.
				
				\bibitem[Seleson et~al., 2016]{seleson2016consistency}
				Seleson, P., Du, Q., and Parks, M.~L. (2016).
				\newblock On the consistency between nearest-neighbor peridynamic
				discretizations and discretized classical elasticity models.
				\newblock {\em Computer Methods in Applied Mechanics and Engineering},
				311:698--722.
				
				\bibitem[Shojaei et~al., 2017]{shojaei2017coupling}
				Shojaei, A., Zaccariotto, M., and Galvanetto, U. (2017).
				\newblock Coupling of 2d discretized peridynamics with a meshless method based
				on classical elasticity using switching of nodal behaviour.
				\newblock {\em Engineering Computations}, 34(5):1334--1366.
				
				\bibitem[Silling, 2003]{silling2003dynamic}
				Silling, S. (2003).
				\newblock Dynamic fracture modeling with a meshfree peridynamic code.
				\newblock In {\em Computational fluid and solid mechanics 2003}, pages
				641--644. Elsevier.
				
				\bibitem[Silling et~al., 2010]{CMPer-Silling5}
				Silling, S., Weckner, O., Askari, E., and Bobaru, F. (2010).
				\newblock Crack nucleation in a peridynamic solid.
				\newblock {\em International Journal of Fracture}, 162(1-2):219--227.
				
				\bibitem[Silling, 2000]{CMPer-Silling}
				Silling, S.~A. (2000).
				\newblock Reformulation of elasticity theory for discontinuities and long-range
				forces.
				\newblock {\em Journal of the Mechanics and Physics of Solids}, 48(1):175--209.
				
				\bibitem[Silling and Askari, 2005]{silling2005meshfree}
				Silling, S.~A. and Askari, E. (2005).
				\newblock A meshfree method based on the peridynamic model of solid mechanics.
				\newblock {\em Computers \& structures}, 83(17):1526--1535.
				
				\bibitem[Silling and Bobaru, 2005]{SillBob}
				Silling, S.~A. and Bobaru, F. (2005).
				\newblock Peridynamic modeling of membranes and fibers.
				\newblock {\em International Journal of Non-Linear Mechanics}, 40(2):395--409.
				
				\bibitem[Silling et~al., 2007]{States}
				Silling, S.~A., Epton, M., Weckner, O., Xu, J., and Askari, E. (2007).
				\newblock Peridynamic states and constitutive modeling.
				\newblock {\em Journal of Elasticity}, 88(2):151--184.
				
				\bibitem[Silling and Lehoucq, 2008]{CMPer-Silling4}
				Silling, S.~A. and Lehoucq, R.~B. (2008).
				\newblock Convergence of peridynamics to classical elasticity theory.
				\newblock {\em Journal of Elasticity}, 93(1):13--37.
				
				\bibitem[Silling and Lehoucq, 2010]{silling2010peridynamic}
				Silling, S.~A. and Lehoucq, R.~B. (2010).
				\newblock Peridynamic theory of solid mechanics.
				\newblock {\em Advances in applied mechanics}, 44:73--168.
				
				\bibitem[Trageser and Seleson, 2020]{trageser2020bond}
				Trageser, J. and Seleson, P. (2020).
				\newblock Bond-based peridynamics: A tale of two poisson’s ratios.
				\newblock {\em Journal of Peridynamics and Nonlocal Modeling}, 2(3):278--288.
				
				\bibitem[Trask et~al., 2019]{trask2019asymptotically}
				Trask, N., You, H., Yu, Y., and Parks, M.~L. (2019).
				\newblock An asymptotically compatible meshfree quadrature rule for nonlocal
				problems with applications to peridynamics.
				\newblock {\em Computer Methods in Applied Mechanics and Engineering},
				343:151--165.
				
				\bibitem[Weckner and Abeyaratne, 2005]{WeckAbe}
				Weckner, O. and Abeyaratne, R. (2005).
				\newblock The effect of long-range forces on the dynamics of a bar.
				\newblock {\em Journal of the Mechanics and Physics of Solids}, 53(3):705--728.
				
				\bibitem[Wildman et~al., 2017]{CMPer-Wildman2}
				Wildman, R.~A., O’Grady, J.~T., and Gazonas, G.~A. (2017).
				\newblock A hybrid multiscale finite element/peridynamics method.
				\newblock {\em International Journal of Fracture}, 207(1):41--53.
				
				\bibitem[Yang et~al., 2019]{yang2019implementation}
				Yang, Z., Oterkus, E., Nguyen, C.~T., and Oterkus, S. (2019).
				\newblock Implementation of peridynamic beam and plate formulations in finite
				element framework.
				\newblock {\em Continuum Mechanics and Thermodynamics}, 31:301--315.
				
			\end{thebibliography}

		\end{document}